\documentclass[12pt,a4paper,twoside]{article}
\usepackage{latexsym,amsmath,amssymb,amsfonts,amsthm}
\usepackage[mathscr]{eucal}
\usepackage{makeidx}

\newtheorem{definition}{Definition}[section]

\newtheorem{theorem}{Theorem}[section]
\def\nabla{\bigtriangledown}
\makeindex

\begin{document}

\title{ Nonlinear Connections and Spinor Geometry}
\author{Sergiu I.\ Vacaru \thanks{%
e--mail: vacaru@fisica.ist.utl.pt, sergiu$_{-}$vacaru@yahoo.com,\ $%
^\diamondsuit$ nadejda$_{-}$vicol@yahoo.com }, and Nadejda A. Vicol $%
^\diamondsuit$ \\
\\
{\small * Centro Multidisciplinar de Astrofisica - CENTRA,}\\
{\small Departamento de Fisica, Instituto Superior Tecnico,}\\
{\small Av. Rovisco Pais 1, Lisboa, 1049-001, Portugal }\\
{} \\
{\small $\diamondsuit$ Faculty of Mathematics and Informatics, gr. 33 MI,}\\
{\small State University of Moldova, Mateevici str. 60,}\\
{\small Chi\c sin\v au MD2009, Republic of Moldova} }
\date{June 28, 2004}
\maketitle

\begin{abstract}
We present an introduction to the geometry of higher order vector and
co--vector bundles (including higher order generalizations of the Finsler
geometry and Kaluza--Klein gravity) and review the basic results on Clifford
and spinor structures on spaces with generic local anisotropy modeled by
anholonomic frames with associated nonlinear connection structures. We
emphasize strong arguments for application of Finsler like geometries in
modern string and gravity theory and noncommutative geometry and
noncommutative field theory and gravity.%


\end{abstract}

\newpage

\section{Introduction}

Nowadays, it has been established an interest to non--Riemannian geometries
derived in the low energy string theory \cite{sgr}, noncommutative geometry %
\cite{ncg,vsv} and quantum groups \cite{majid}. Various types of Finsler
like structures can be parametrized by generic off--diagonal metrics, which
can not be diagonalized by coordinate transforms but only by anholonomic
maps with associated nonlinear connection (in brief, N--connection). Such
structures may be defined as exact solutions of gravitational field
equations in the Einstein gravity and its generalizations \cite{vb,vsv,vrec}%
, for instance, in the metric--affine \cite{mag} Riemann--Cartan gravity %
\cite{rcg}. Finsler like configurations are considered in locally
anisotropic thermodynamics and kinetics and related stochastic processes %
\cite{vtherm} and (super) string theory \cite{vst96,vstr2,vbook}.

The following natural step in these lines is to elucidate the theory of
spinors in effectively derived Finsler geometries and to relate this
formalism of Clifford structures to non--commutative Finsler geometry. It
should be noted that the rigorous definition of spinors for Finsler spaces
and generalizations was not a trivial task because (on such spaces) there
are not defined even local groups of authomorphisms. The problem was solved
in Refs. \cite{vjmp,vsp1,vhsp} by adapting the geometric constructions with
respect to anholonomic frames with associated N--connection structure. The
aim of this work is to outline the geometry of generalized Finsler spinors
in a form more oriented to applications in modern mathematical physics.

We start with some historical remarks:\ The spinors studied by
mathematicians and physicists are connected with the general theory of
Clifford spaces introduced in 1876 \cite{clifford}. The theory of spinors
and Clifford algebras play a major role in contemporary physics and
mathematics. The spinors were discovered by \`{E}lie Cartan in 1913 in
mathematical form in his researches on representation group theory \cite%
{car38}; he showed that spinors furnish a linear representation of the
groups of rotations of a space of arbitrary dimensions. Physicists Pauli %
\cite{pauli} and Dirac \cite{dirac} (in 1927, respectively, for the
three--dimensional and four--dimensional space--time) introduced spinors for
the representation of the wave functions. In general relativity theory
spinors and the Dirac equations on (pseudo) Riemannian spaces were defined
in 1929 by H. Weyl \cite{weyl}, V. Fock \cite{foc} and E. Schr\"{o}dinger %
\cite{schr}. The books \cite{pen} by R. Penrose and W. Rindler monograph
summarize the spinor and twistor methods in space--time geometry (see
additional references \cite{hladik} on Clifford structures and spinor
theory).

Spinor variables were introduced in Finsler geometries by Y. Takano in 1983 %
\cite{t1} where he dismissed anisotropic dependencies not only on vectors on
the tangent bundle but on some spinor variables in a spinor bundle on a
space--time manifold. Then generalized Finsler geometries, with spinor
variables, were developed by T. Ono and Y. Takano in a series of
publications during 1990--1993 \cite{ono}. The next steps were
investigations of anisotropic and deformed geometries with spinor and vector
variables and applications in gauge and gravity theories elaborated by P.
Stavrinos and his students, S. Koutroubis, P. Manouselis, and V. Balan
starting from 1994 \cite{sk.sk1,sk2,sk3}. In those works the authors assumed
that some spinor variables may be introduced in a Finsler-like way but they
did not relate the Finlser metric to a Clifford structure and restricted the
spinor--gauge Finsler constructions only for antisymmetric spinor metrics on
two--spinor fibers with possible generalizations to four dimensional Dirac
spinors.

Isotopic spinors, related with $SU(2)$ internal structural groups, were
considered in generalized Finsler gravity and gauge theories also by G.
Asanov and S. Ponomarenko \cite{asa88}, in 1988. In that book, and in other
papers on Finsler geometry with spinor variables, the authors did not
investigate the possibility of introducing a rigorous mathematical
definition of spinors on spaces with generic local anisotropy.

An alternative approach to spinor differential geometry and generalized
Finsler spaces was elaborated, starting from 1994, in a series of papers and
communications by S. Vacaru and co--authors \cite{vsp0,vsp1,vsp2,vsp3,vdeb}.
This direction originates from Clifford algebras and Clifford bundles \cite%
{kar} and Penrose's spinor and twistor space--time geometry \cite{pen},
which were re--considered for the case of nearly autoparallel maps
(generalized conformal transforms) in Refs. \cite{v87}. In the works \cite%
{vsp0,vjmp,vsp1,vsp2,vsp3}, a rigorous definition of spinors for Finsler
spaces, and their generalizations, was given. It was proven that a Finsler,
or Lagrange, metric (in a tangent, or, more generally, in a vector bundle)
induces naturally a distinguished Clifford (spinor) structure which is
locally adapted to the nonlinear connection structure. Such spinor spaces
could be defined for arbitrary dimensions of base and fiber subspaces, their
spinor metrics are symmetric, antisymmetric or nonsymmetric, depending on
the corresponding base and fiber dimensions. That work resulted in formation
of the spinor differential geometry of generalized Finsler spaces and
developed a number of geometric applications to the theory of gravitational
and matter field interactions with generic local anisotropy.

The geometry of anisotropic spinors and (distinguished by nonlinear
connections) Clifford structures was elaborated for higher order anisotropic
spaces \cite{vsp1,vhsp,vbook} and, more recently, for Hamilton and Lagrange
spaces \cite{vsv}.

We emphasize that the theory of anisotropic spinors may be related not only
to generalized Finsler, Lagrange, Cartan and Hamilton spaces or their higher
order generalizations, but also to anholonomic frames with associated
nonlinear connections which appear naturally even in (pseudo) Riemannian and
Riemann--Cartan geometries if off--diagonal metrics are considered \cite{vb}%
. In order to construct exact solutions of the Einstein equations in general
relativity and extra dimension gravity (for lower dimensions see \cite%
{vtherm}), it is more convenient to diagonalize space--time metrics by using
some anholonomic transforms. As a result one induces locally anisotropic
structures on space--time which are related to anholonomic (anisotropic)
spinor structures.

The main purpose of the present review is to present a detailed summary and
new results on spinor differential geometry for generalized Finsler spaces
and (pseudo) Riemannian space--times provided with anholonomic frame and
associated nonlinear connection structure, to discuss and compare the
existing approaches and to consider applications to modern gravity and gauge
theories.

\section{(Co) Vector Bundles and N--Connections}

We outline the basic definitions and denotations for the vector and tangent
(and theirs dual spaces) bundles and higher order vector/covector bundle
geometry. In this work, we consider that the space--time geometry can be
modeled both on a (pseudo) Riemannian manifold $V^{[n+m]}$ of dimension $n+m$
and/or on a vector bundle (or its dual, covector bundle) being, for
simplicity, locally trivial with a base space $M$ of dimension $n$ and a
typical fiber $F$ (cofiber $F^{\ast }$) of dimension $m,$ or as a higher
order extended vector/covector bundle (we follow the geometric constructions
and definitions of monographs \cite{ma87,mhss} which were generalized for
vector superbundles in Refs. \cite{vstr2,vbook}). \ Such (pseudo) Riemanian
spaces and/or vector/covector bundles enabled with compatible fibered and/or
anholonomic structures are called \textbf{anisotropic space--times}. If the
anholonomic structure with associated nonlinear connection is modeled on
higher order vector/covector bundles we use the term of \textbf{higher order
anisotropic space--time.} In this section, we usually shall omit proofs
which can be found in the mentioned monographs \cite{ma87,mhss,vbook}.

\subsection{(Co) vector and tangent bundles}

A locally trivial \textbf{vector bundle}, in brief, \textbf{v--bundle}, $%
\mathcal{E}=\left( E,\pi ,M,Gr,F\right) $ is introduced as a set of spaces
and surjective map with the properties that a real vector space $F=\mathcal{R%
}^{m}$ of dimension $m$ ($\dim F=m,$ $\mathcal{R\ }$ denotes the real number
field) defines the typical fiber, the structural group is chosen to be the
group of automorphisms of $\mathcal{R}^{m},$ i. e. $Gr=GL\left( m,\mathcal{R}%
\right) ,\ $ and $\pi :E\rightarrow M$ is a differentiable surjection of a
differentiable manifold $E$ (total space, $\dim E=n+m)$ to a differentiable
manifold $M$ $\left( \mbox{base
space, }\dim M=n\right) .$ The local coordinates on $\mathcal{E}$ are
denoted $u^{{\alpha }}=\left( x^{{i}},y^{{a}}\right) ,$ or in brief ${%
u=\left( x,y\right) }$ (the Latin indices ${i,j,k,...}=1,2,...,n$ define
coordinates of geometrical objects with respect to a local frame on base
space $M;$\ the Latin indices ${a,b,c,...=}1,2,...,m$ define fiber
coordinates of geometrical objects and the Greek indices ${\alpha ,\beta
,\gamma ,...}$ are considered as cumulative ones for coordinates of objects
defined on the total space of a v-bundle).

Coordinate transforms $u^{{\ \alpha ^{\prime }\ }} = u^{{\ \alpha ^{\prime
}\ }}\left( u^{{\ \alpha }}\right) $ on a v--bundle $\mathcal{E}$ are
defined $(x^{{\ i}}, y^{{\ a}}) \rightarrow ( x^{{\ i^{\prime }\ }}, y^{{\
a^{\prime }}}) ,$ where
\begin{equation}
x^{{i^{\prime }\ }} = x^{{\ i^{\prime }\ }}(x^{{\ i}}),\qquad y^{{a^{\prime
}\ }} = K_{{\ a\ }}^{{\ a^{\prime }}}(x^{{i\ }})y^{{a}}  \label{coordtr}
\end{equation}
and matrix $K_{{\ a\ }}^{{\ a^{\prime }}}(x^{{\ i\ }}) \in GL\left( m,%
\mathcal{R}\right) $ are functions of necessary smoothness class.

A local coordinate parametrization of v--bundle $\mathcal{E}$ naturally
defines a coordinate basis
\begin{equation}
\partial _{\alpha }=\frac{\partial }{\partial u^{\alpha }}=\left( \partial
_{i}=\frac{\partial }{\partial x^{i}},\ \partial _{a}=\frac{\partial }{%
\partial y^{a}}\right) ,  \label{pder}
\end{equation}%
and the reciprocal to (\ref{pder}) coordinate basis
\begin{equation}
d^{\alpha }=du^{\alpha }=(d^{i}=dx^{i},\ d^{a}=dy^{a})  \label{pdif}
\end{equation}%
which is uniquely defined from the equations $d^{\alpha }\circ \partial
_{\beta }=\delta _{\beta }^{\alpha },$ where $\delta _{\beta }^{\alpha }$ is
the Kronecher symbol and by ''$\circ $'' we denote the inner (scalar)
product in the tangent bundle $\mathcal{TE}.$

A \textbf{tangent bundle} (in brief, \textbf{t--bundle}) $(TM,\pi ,M)$ to a
manifold $M$ can be defined as a particular case of a v--bundle when the
dimension of the base and fiber spaces (the last one considered as the
tangent subspace) are identic, $n=m.$ In this case both type of indices $%
i,k,...$ and $a,b,...$ take the same values $1,2,...n$. For t--bundles the
matrices of fiber coordinates transforms from (\ref{coordtr}) can be written
$K_{{\ i\ }}^{{\ i^{\prime }}}={\partial x^{i^{\prime }}}/{\partial x^{i}}.$

We shall also use the concept of \textbf{covector bundle}, (in brief,
\textbf{cv--bundles)} \newline
$\breve{\mathcal{E}}=\left( {\breve{E}},\pi ^{\ast },M,Gr,F^{\ast }\right) $%
, which is introduced as a dual vector bundle for which the typical fiber $%
F^{\ast }$ (cofiber) is considered to be the dual vector space (covector
space) to the vector space $F.$ The fiber coordinates $p_{a}$ of $\breve{E}$
are dual to $y^{a}$ in $E.$ The local coordinates on total space $\breve{E}$
are denoted $\breve{u}=(x,p)=(x^{i},p_{a}).$ The coordinate transform on $%
\breve{E},$ ${\breve{u}}=(x^{i},p_{a})\rightarrow {\breve{u}}^{\prime
}=(x^{i^{\prime }},p_{a^{\prime }}), $ are written
\begin{equation}
x^{{i^{\prime }\ }}=x^{{\ i^{\prime }\ }}(x^{{\ i}}),\qquad p_{{a^{\prime }\
}}=K_{{\ a^{\prime }}}^{{\ a\ }}(x^{{i\ }})p_{{a}}.  \label{coordtrd}
\end{equation}%
The coordinate bases on $E^{\ast }$ are denoted
\begin{equation}
{\breve{\partial}}_{\alpha }=\frac{\breve{\partial}}{\partial u^{\alpha }}%
=\left( \partial _{i}=\frac{\partial }{\partial x^{i}},{{\breve{\partial}}%
^{a}}=\frac{\breve{\partial}}{\partial p_{a}}\right)  \label{pderct}
\end{equation}%
and
\begin{equation}
{\breve{d}}^{\alpha }={\breve{d}}u^{\alpha }=\left( d^{i}=dx^{i},{\breve{d}}%
_{a}=dp_{a}\right) .  \label{pdifct}
\end{equation}%
We use ''breve'' symbols in order to distinguish the geometrical objects on
a cv--bundle $\mathcal{E}^{\ast }$ from those on a v--bundle $\mathcal{E}$.

As a particular case with the same dimension of base space and cofiber one
obtains the \textbf{cotangent bundle} $(T^{\ast }M,\pi ^{\ast },M)$ , in
brief, \textbf{ct--bundle,} being dual to $TM.$ The fibre coordinates $p_{i}
$ of $T^{\ast }M$ are dual to $y^{i}$ in $TM.$ The coordinate transforms (%
\ref{coordtrd}) on $T^{\ast }M$ are stated by some matrices $K_{{\ k^{\prime
}}}^{k}(x^{i})={\partial x^{k}}/{\partial x^{k^{\prime }}}.$

In our further considerations we shall distinguish the base and cofiber
indices.

\subsection{Higher order (co) vector bundles}

The geometry of higher order tangent and cotangent bundles provided with
nonlinear connection structure was elaborated in Refs. \cite{mhss} in order
to geometrize the higher order Lagrange and Hamilton mechanics. In this case
we have base spaces and fibers of the same dimension. To develop the
approach to modern high energy physics (in superstring and Kaluza--Klein
theories), we introduced (in Refs \cite{vsp1,vhsp,vbook,vstr2}) the concept
of higher order vector bundle with the fibers consisting from finite
'shells'' of vector, or covector, spaces of different dimensions not
obligatory coinciding with the base space dimension.

\begin{definition}
A distinguished vector/covector space, in brief dvc--space, of type
\begin{equation}
{\tilde{F}}=F[v(1),v(2),cv(3),...,cv(z-1),v(z)]  \label{orient}
\end{equation}%
is a vector space decomposed into an invariant oriented direct sum
\begin{equation*}
{\tilde{F}}=F_{(1)}\oplus F_{(2)}\oplus F_{(3)}^{\ast }\oplus ...\oplus
F_{(z-1)}^{\ast }\oplus F_{(z)}
\end{equation*}%
of vector spaces $F_{(1)},F_{(2)},...,F_{(z)}$ of respective dimensions
\begin{equation*}
dimF_{(1)}=m_{1},dimF_{(2)}=m_{2},...,dimF_{(z)}=m_{z}
\end{equation*}%
and of covector spaces $F_{(3)}^{\ast },...,F_{(z-1)}^{\ast }$ of respective
dimensions
\begin{equation*}
dimF_{(3)}^{\ast }=m_{3}^{\ast },...,dimF_{(z-1)}^{\ast }=m_{(z-1)}^{\ast }.
\end{equation*}
\end{definition}

As a particular case we obtain a distinguished vector space, in brief
dv--space (a distinguished covector space, in brief dcv--space), if all
components of the sum are vector (covector) spaces. We note that we have
fixed, for simplicity, an orientation of (co) vector subspaces like in (\ref%
{orient}).

Coordinates on ${\tilde{F}}$ are denoted
\begin{equation*}
\tilde{y}=(y_{(1)},y_{(2)},p_{(3)},...,p_{(z-1)},y_{(z)})=\{y^{<\alpha
_{z}>}\}=(y^{a_{1}},y^{a_{2}},p_{a_{3}},...,p_{a_{z-1}},y^{a_{z}}),
\end{equation*}%
where indices run correspondingly the values $a_{1}=1,2,...,m_{1};\
a_{2}=1,2,...,m_{2},\ ...,a_{z}=1,2,...,m_{z}.$

\begin{definition}
A higher order vector/covector bundle (in brief, hvc-\--bund\-le) of type ${%
\tilde{\mathcal{E}}}={\tilde{\mathcal{E}}}[v(1),v(2),cv(3),...,cv(z-1),v(z)]$
is a vector bundle ${\tilde{\mathcal{E}}}=({\tilde{E}},p^{<d>},{\tilde{F}}%
,M) $ with corresponding total, ${\tilde{E}}$, and base, $M,$ spaces,
surjective projection $p^{<d>}:\ {\tilde{E}}\rightarrow M$ and typical fiber
${\tilde{F}}.$
\end{definition}

We define the higher order vector (covector) bundles, in brief, hv--bundles
(in brief, hcv--bundles), if the typical fibre is a dv--space (dcv--space)
as particular cases of the hvc--bundles.

A hvc--bundle is constructed as an oriented set of enveloping 'shell by
shell' v--bundles and/or cv--bundles,
\begin{equation*}
p^{<s>}:\ {\tilde{E}}^{<s>}\rightarrow {\tilde{E}}^{<s-1>},
\end{equation*}%
where we use the index $<s>=0,1,2,...,z$ in order to enumerate the shells,
when ${\tilde{E}}^{<0>}=M.$ Local coordinates on ${\tilde{E}}^{<s>}$ are
denoted
\begin{equation}
{\tilde{u}}_{(s)}=(x,{\tilde{y}}%
_{<s>})=(x,y_{(1)},y_{(2)},p_{(3)},...,y_{(s)})=(x^{i},y^{a_{1}},y^{a_{2}},p_{a_{3}},...,y^{a_{s}}).
\notag
\end{equation}%
If $<s>=<z>$ we obtain a complete coordinate system on ${\tilde{\mathcal{E}}}
$ denoted in brief
\begin{equation*}
\tilde{u}=(x,{\tilde{y}})=\tilde{u}^{\alpha
}=(x^{i}=y^{a_{0}},y^{a_{1}},y^{a_{2}},p_{a_{3}},...,p_{a_{z-1}},y^{a_{z}}).
\end{equation*}%
We shall use the general commutative indices $\alpha ,\beta ,...$ for
objects on hvc---bundles which are marked by tilde, like $\tilde{u},\tilde{u}%
^{\alpha },...,$ ${\tilde{E}}^{<s>},....$

The coordinate bases on ${\tilde{\mathcal{E}}}$ are denoted
\begin{eqnarray}
{\tilde{\partial}}_{\alpha } &=&\frac{{\tilde{\partial}}}{\partial u^{\alpha
}}=\left( \partial _{i}=\frac{\partial }{\partial x^{i}},\partial _{a_{1}}=%
\frac{\partial }{\partial y^{a_{1}}},\partial _{a_{2}}=\frac{\partial }{%
\partial y^{a_{2}}},{{\breve{\partial}}}^{a_{3}}=\frac{\breve{\partial}}{%
\partial p_{a_{3}}},...,\partial _{a_{z}}=\frac{\partial }{\partial y^{a_{z}}%
}\right)  \label{pderho} \\
&&  \notag
\end{eqnarray}%
and
\begin{eqnarray}
{\tilde{d}}^{\alpha } &=&{\tilde{d}}u^{\alpha }=\left(
d^{i}=dx^{i},d^{a_{1}}=dy^{a_{1}},d^{a_{2}}=dy^{a_{2}},{\breve{d}}%
_{a_{3}}=dp_{a_{3}},...,d^{a_{z}}=dy^{a_{z}}\right) .  \label{pdifho} \\
&&  \notag
\end{eqnarray}

We give two examples of higher order tangent / co\-tan\-gent bundles (when
the dimensions of fibers/cofibers coincide with the dimension of bundle
space, see Refs. \cite{mhss}).

\subsubsection{Osculator bundle}

The $k$--osculator bundle is identified with the $k$--tangent bundle $\left(
T^{k}M,p^{(k)},M\right) $ of a $n$--dimensional manifold $M.$ We denote the
local coordinates ${\tilde{u}}^{\alpha }=\left(
x^{i},y_{(1)}^{i},...,y_{(k)}^{i}\right) ,$ where we have identified $%
y_{(1)}^{i}\simeq y^{a_{1}},...,y_{(k)}^{i}\simeq y^{a_{k}},k=z,$ in order
to to have similarity with denotations from \cite{mhss}. The coordinate
transforms ${\tilde{u}}^{\alpha ^{\prime }}\rightarrow {\tilde{u}}^{\alpha
^{\prime }}\left( {\tilde{u}}^{\alpha }\right)$ preserving the structure of
such higher order vector bundles are parametrized
\begin{eqnarray*}
x^{i^{\prime }} &=&x^{i^{\prime }}\left( x^{i}\right) ,\det \left( \frac{%
\partial x^{i^{\prime }}}{\partial x^{i}}\right) \neq 0, \\
y_{(1)}^{i^{\prime }} &=&\frac{\partial x^{i^{\prime }}}{\partial x^{i}}%
y_{(1)}^{i},\ 2y_{(2)}^{i^{\prime }}=\frac{\partial y_{(1)}^{i^{\prime }}}{%
\partial x^{i}}y_{(1)}^{i}+2\frac{\partial y_{(1)}^{i^{\prime }}}{\partial
y^{i}}y_{(2)}^{i}, \\
&& \\
&&................... \\
ky_{(k)}^{i^{\prime }} &=&\frac{\partial y_{(1)}^{i^{\prime }}}{\partial
x^{i}}y_{(1)}^{i}+...+k\frac{\partial y_{(k-1)}^{i^{\prime }}}{\partial
y_{(k-1)}^{i}}y_{(k)}^{i},
\end{eqnarray*}%
where the equalities
\begin{equation*}
\frac{\partial y_{(s)}^{i^{\prime }}}{\partial x^{i}}=\frac{\partial
y_{(s+1)}^{i^{\prime }}}{\partial y_{(1)}^{i}}=...=\frac{\partial
y_{(k)}^{i^{\prime }}}{\partial y_{(k-s)}^{i}}
\end{equation*}%
hold for $s=0,...,k-1$ and $y_{(0)}^{i}=x^{i}.$

The natural coordinate frame on $\left( T^{k}M,p^{(k)},M\right) $ is defined
${\tilde{\partial}}_{\alpha }= (\frac{\partial }{\partial x^{i}},\frac{%
\partial }{\partial y_{(1)}^{i}},...,\frac{\partial }{\partial y_{(k)}^{i}})$
and the coframe is ${\tilde{d}}_{\alpha}
=(dx^{i},dy_{(1)}^{i},...,dy_{(k)}^{i}).$ These formulas are respectively
some particular cases of$\ $(\ref{pderho}) and (\ref{pdifho}).

\subsubsection{The dual bundle of k--osculator bundle}

This higher order vector/covector bundle, denoted as $\left( T^{\ast
k}M,p^{\ast k},M\right) ,$ is defined as the dual bundle to the k--tangent
bundle $\left( T^{k}M,p^{k},M\right) .$ The local coordinates (parametrized
as in the previous paragraph) are
\begin{equation*}
\tilde{u}=\left( x,y_{(1)},...,y_{(k-1)},p\right) =\left(
x^{i},y_{(1)}^{i},...,y_{(k-1)}^{i},p_{i}\right) \in T^{\ast k}M.
\end{equation*}%
The coordinate transforms on $\left( T^{\ast k}M,p^{\ast k},M\right) $ are
\begin{eqnarray*}
x^{i^{\prime }} &=&x^{i^{\prime }}\left( x^{i}\right) ,\det \left( \frac{%
\partial x^{i^{\prime }}}{\partial x^{i}}\right) \neq 0, \\
y_{(1)}^{i^{\prime }} &=&\frac{\partial x^{i^{\prime }}}{\partial x^{i}}%
y_{(1)}^{i},\ 2y_{(2)}^{i^{\prime }}=\frac{\partial y_{(1)}^{i^{\prime }}}{%
\partial x^{i}}y_{(1)}^{i}+2\frac{\partial y_{(1)}^{i^{\prime }}}{\partial
y^{i}}y_{(2)}^{i}, \\
&&................... \\
(k-1)y_{(k-1)}^{i^{\prime }} &=&\frac{\partial y_{(k-2)}^{i^{\prime }}}{%
\partial x^{i}}y_{(1)}^{i}+...+k\frac{\partial y_{(k-1)}^{i^{\prime }}}{%
\partial y_{(k-2)}^{i}}y_{(k-1)}^{i},\ p_{i^{\prime }}=\frac{\partial x^{i}}{%
\partial x^{i^{\prime }}}p_{i},
\end{eqnarray*}%
where the equalities
\begin{equation*}
\frac{\partial y_{(s)}^{i^{\prime }}}{\partial x^{i}}=\frac{\partial
y_{(s+1)}^{i^{\prime }}}{\partial y_{(1)}^{i}}=...=\frac{\partial
y_{(k-1)}^{i^{\prime }}}{\partial y_{(k-1-s)}^{i}}
\end{equation*}%
hold for $s=0,...,k-2$ and $y_{(0)}^{i}=x^{i}.$

The natural coordinate frame on $\left( T^{\ast k}M,p^{\ast (k)},M\right) $
is written in the form\newline
${\tilde{\partial}}_{\alpha }= ( \frac{\partial }{\partial x^{i}},\frac{%
\partial }{\partial y_{(1)}^{i}},...,\frac{\partial }{\partial y_{(k-1)}^{i}}%
,\frac{\partial }{\partial p_{i}})$ and the coframe is written ${\tilde{d}}%
_{\alpha }=\left( dx^{i},dy_{(1)}^{i},...,dy_{(k-1)}^{i},dp_{i}\right) .$
These formulas are, respespectively, certain particular cases of (\ref%
{pderho}) and (\ref{pdifho}).

\subsection{Nonlinear Connections}

The concept of \textbf{nonlinear connection,} in brief, N-connection, is
fundamental in the geometry of vector bundles and anisotropic spaces (see a
detailed study and basic references in \cite{ma87} and, for supersymmetric
and/or spinor bundles in \cite{vbook,vp,vstr2}). A rigorous mathematical
definition is possible by using the formalism of exact sequences of vector
bundles.

\subsubsection{N--connections in vector bundles}

Let $\mathcal{E=}=(E,p,M)$ be a v--bundle with typical fiber $\mathcal{R}^m$
and $\pi ^T:\ TE\to TM$ being the differential of the map $P$ which is a
fibre--preserving morphism of the tangent bundle $TE,\tau _E,E)\to E$ and of
tangent bundle $(TM,\tau ,M)\to M.$ The kernel of the vector bundle
morphism, denoted as $(VE,\tau _V,E),$ is called the \textbf{vertical
subbundle} over $E,$ which is a vector subbundle of the vector bundle $%
(TE,\tau _E,E).$

A vector $X_u$ tangent to a point $u\in E$ is locally written as $%
(x,y,X,Y)=(x^i,y^a,X^i,Y^a),$ where the coordinates $(X^i,Y^a)$ are defined
by the equality $X_u=X^i\partial _i+Y^a\partial _a.$ We have\ $\pi
^T(x,y,X,Y)=(x,X).$ Thus the submanifold $VE$ contains the elements which
are locally represented as $(x,y,0,Y).$

\begin{definition}
\label{ncon} A nonlinear connection $\mathbf{N}$ in a vector bundle $%
\mathcal{E}=(E,\pi ,M)$ is the splitting on the left of the exact sequence
\begin{equation*}
0\mapsto VE\mapsto TE\mapsto TE/VE\mapsto 0
\end{equation*}
where $TE/VE$ is the factor bundle.
\end{definition}

By definition (\ref{ncon}) it is defined a morphism of vector bundles $C:\
TE\to VE$ such the superposition of maps $C\circ i$ is the identity on $VE,$
where $i:\ VE\mapsto VE.$ The kernel of the morphism $C$ is a vector
subbundle of $(TE,\tau _E,E)$ which is the horizontal subbundle, denoted by $%
(HE,\tau _H,E).$ Consequently, we can prove that in a v-bundle $\mathcal{E}$
a N--connection can be introduced as a distribution
\begin{equation*}
\{N:\ E_u\rightarrow H_uE,T_uE=H_uE\oplus V_uE\}
\end{equation*}
for every point $u\in E$ defining a global decomposition, as a Whitney sum,
into horizontal, $H\mathcal{E},$ and vertical, $V\mathcal{E},$ subbundles of
the tangent bundle $T\mathcal{E}$
\begin{equation}
T\mathcal{E}=H\mathcal{E}\oplus V\mathcal{E}.  \label{whitney}
\end{equation}

Locally a N-connection in a v--bundle $\mathcal{E}$ is given by its
coefficients $N_{{\ i}}^{{a}}({\ u})=N_{{i}}^{{a}}({x,y})$ with respect to
bases (\ref{pder}) and (\ref{pdif}), $\mathbf{N}=N_{i}^{~a}(u)d^{i}\otimes
\partial _{a}.$ We note that a linear connection in a v--bundle $\mathcal{E}$
can be considered as a particular case of a N--connection when $%
N_i^{~a}(x,y)=K_{bi}^a\left( x\right) y^b,$ where functions $K_{ai}^b\left(
x\right) $ on the base $M$ are called the Christoffel coefficients.

\subsubsection{N--connections in covector bundles}

A nonlinear connection in a cv--bundle ${\breve{\mathcal{E}}}$ (in brief a
\v{N}--connection) can be introduced in a similar fashion as for v--bundles
by reconsidering the corresponding definitions for cv--bundles. For
instance, it is stated by a Whitney sum, into horizontal, $H{\breve{\mathcal{%
E}}},\ $ and vertical, $V{\breve{\mathcal{E}}},$ subbundles of the tangent
bundle $T{\breve{\mathcal{E}}}:$
\begin{equation}
T{\breve{\mathcal{E}}}=H{\breve{\mathcal{E}}}\oplus V{\breve{\mathcal{E}}}.
\label{whitneyc}
\end{equation}

Hereafter, for the sake of brevity, we shall omit details on definition of
geometrical objects on cv--bundles if they are very similar to those for
v--bundles:\ we shall present only the basic formulas by emphasizing the
most important common points and differences.

\begin{definition}
\label{ctvn} A \v{N}--connection on ${\breve{\mathcal{E}}}$ is a
differentiable distribution
\begin{equation*}
\breve{N}:\ {\breve{\mathcal{E}}}\rightarrow {\breve{N}}_{u}\in T_{u}^{\ast }%
{\breve{\mathcal{E}}}
\end{equation*}%
which is suplimentary to the vertical distribution $V,$ i. e. $T_{u}{\breve{%
\mathcal{E}}}={\breve{N}}_{u}\oplus {\breve{V}}_{u},\forall {\breve{\mathcal{%
E}}}.$
\end{definition}

The same definition is true for \v N--connections in ct--bundles, we have to
change in the definition (\ref{ctvn}) the symbol ${\breve{\mathcal{E}}}$
into $T^*M.$

A N--connection in a cv--bundle ${\breve{\mathcal{E}}}$ is given locally by
its coefficients ${\breve N}_{{\ ia}}({\ u})= {\breve N}_{{ia}}({x,p})$ with
respect to bases (\ref{pder}) and (\ref{pdif}), $\mathbf{\breve N}={\breve N}%
_{ia}(u)d^i\otimes {{\breve{\partial}}^a}.$

We emphasize that if a N--connection is introduced in a v--bundle
(cv--bundle), we have to adapt the geometric constructions to the
N--connection structure.

\subsubsection{N--connections in higher order bundles}

The concept of N--connection can be defined for higher order vector /
covec\-tor bundle in a standard manner like in the usual vector bundles:

\begin{definition}
A nonlinear connection ${\tilde{\mathbf{N}}}$ in hvc--bundle
\index{N--connection}
\begin{equation*}
{%
\tilde{\mathcal{E}}}={\tilde{\mathcal{E}}}[v(1),v(2),cv(3),...,cv(z-1),v(z)]
\end{equation*}
is a splitting of the left of the exact sequence
\begin{equation}
0\to V{\tilde{\mathcal{E}}}\to T{\tilde{\mathcal{E}}}\to T{\tilde{\mathcal{E}%
}}/V{\tilde{\mathcal{E}}}\to 0  \label{exacts}
\end{equation}
\end{definition}

We can associate sequences of type (\ref{exacts}) to every mappings of
intermediary subbundles. For simplicity, we present here the Whitney
decomposition
\begin{equation*}
T{\tilde{\mathcal{E}}}=H{\tilde{\mathcal{E}}}\oplus V_{v(1)}{\tilde{\mathcal{%
E}}}\oplus V_{v(2)}{\tilde{\mathcal{E}}}\oplus V_{cv(3)}^{\ast }{\tilde{%
\mathcal{E}}}\oplus ....\oplus V_{cv(z-1)}^{\ast }{\tilde{\mathcal{E}}}%
\oplus V_{v(z)}{\tilde{\mathcal{E}}}.
\end{equation*}%
Locally a N--connection ${\tilde{\mathbf{N}}}$ in ${\tilde{\mathcal{E}}}$ is
given by its coefficients
\begin{equation}
\begin{array}{llllll}
N_{i}^{~a_{1}}, & N_{i}^{~a_{2}}, & N_{ia_{3}}, & ..., & N_{ia_{z-1}}, &
N_{i}^{~a_{z}}, \\
0, & N_{a_{1}}^{~a_{2}}, & N_{a_{1}a_{3}}, & ..., & N_{a_{1}a_{z-1}}, &
N_{a_{1}}^{~a_{z}}, \\
0, & 0, & N_{a_{2}a_{3}}, & ..., & N_{a_{2}a_{z-1}}, & N_{a_{2}}^{~a_{z}},
\\
..., & ..., & ..., & ..., & ..., & ..., \\
0, & 0, & 0, & ..., & N_{a_{z-2}~a_{z-1}}, & N_{a_{z-2}}^{~a_{z}}, \\
0, & 0, & 0, & ..., & 0, & N^{a_{z-1}a_{z}},%
\end{array}
\label{nconho}
\end{equation}%
which are given with respect to the components of bases (\ref{pderho}) and (%
\ref{pdifho}).

\subsubsection{Anholonomic frames and N--connections}

Having defined a N--connection structure in a (vector, covector, or higher
order vector / covector) bundle, we can adapt with respect to this structure
(by 'N--elonga\-ti\-on') the operators of partial derivatives and
differentials and to consider decompositions of geometrical objects with
respect to adapted bases and cobases.

\paragraph{Anholonomic frames in v--bundles}

In a v--bundle $\mathcal{E}$ provided with a N-connection we can adapt to
this structure the geometric constructions by introducing locally adapted
basis (N--frame, or N--basis):
\begin{equation}
\delta _\alpha =\frac \delta {\delta u^\alpha }=\left( \delta _i=\frac
\delta {\delta x^i}=\partial _i-N_i^{~a}\left( u\right) \partial _a,\partial
_a=\frac \partial {\partial y^a}\right) ,  \label{dder}
\end{equation}
and its dual N--basis, (N--coframe, or N--cobasis),
\begin{equation}
\delta \ ^\alpha =\delta u^\alpha =\left( d^i=\delta x^i=dx^i,\delta
^a=\delta y^a+N_i^{~a}\left( u\right) dx^i\right) .  \label{ddif}
\end{equation}

The\textbf{\ anholonomic coefficients, } $\mathbf{w}=\{w_{\beta \gamma
}^{\alpha }\left( u\right) \},$ of N--frames are defined to satisfy the
relations
\begin{equation}
\left[ \delta _{\alpha },\delta _{\beta }\right] =\delta _{\alpha }\delta
_{\beta }-\delta _{\beta }\delta _{\alpha }=w_{\beta \gamma }^{\alpha
}\left( u\right) \delta _{\alpha }.  \label{anhol}
\end{equation}

A frame bases is holonomic is all anholonomy coefficients vanish (like for
usual coordinate bases (\ref{pdif})), or anholonomic if there are nonzero
values of $w_{\beta \gamma }^{\alpha }.$

The operators (\ref{dder}) and (\ref{ddif}) on a v--bundle $\mathcal{E}$
en\-abled with a N--connecti\-on can be considered as respective equivalents
of the operators of partial derivations and differentials:\ the existence of
a N--connection structure results in 'elongation' of partial derivations on $%
x$--variables and in 'elongation' of differentials on $y$--variables.

The algebra of tensorial distinguished fields $DT( \mathcal{E})$ (d--fields,
d--ten\-sors, d--objects) on $\mathcal{E}$ is introduced as the tensor
algebra $\mathcal{T}=\{\mathcal{T}_{qs}^{pr}\}$ of the v--bundle $\mathcal{E}%
_{\left( d\right) }=\left( H\mathcal{E}\oplus V\mathcal{E},p_{d},\mathcal{E}%
\right),$ where $p_{d}:\ H\mathcal{E}\oplus V\mathcal{E}\rightarrow \mathcal{%
E}.$

\paragraph{Anholonomic frames in cv--bundles}

The anholnomic frames adapted to the \v{N}--con\-nec\-ti\-on structure are
introduced similarly to (\ref{dder}) and (\ref{ddif}):

the locally adapted basis (\v N--basis, or \v N--frame):
\begin{equation}
{\breve \delta}_\alpha = \frac {\breve \delta}{\delta u^\alpha}= \left(
\delta _i=\frac \delta {\delta x^i}= \partial _i + {\breve N}_{ia}\left({%
\breve u}\right) {\breve \partial}^a, {\breve \partial}^a = \frac \partial
{\partial p_a}\right) ,  \label{ddercv}
\end{equation}

and its dual (\v{N}--cobasis, or \v{N}--coframe):
\begin{equation}
{\breve{\delta}}^{\alpha }={\breve{\delta}}u^{\alpha }=\left( d^{i}=\delta
x^{i}=dx^{i},\ {\breve{\delta}}_{a}={\breve{\delta}}p_{a}=dp_{a}-{\breve{N}}%
_{ia}\left( {\breve{u}}\right) dx^{i}\right) .  \label{ddifcv}
\end{equation}

We note that for the singes of \v N--elongations are inverse to those for
N--elongations.

The\textbf{\ anholonomic coefficients, } $\mathbf{\breve w}= \{{\breve w}%
_{\beta \gamma }^\alpha \left({\breve u}\right) \},$ of \v N--frames are
defined by the relations
\begin{equation}  \label{anhola}
\left[{\breve \delta} _\alpha ,{\breve \delta} _\beta \right] ={\breve \delta%
}_\alpha {\breve \delta}_\beta - {\breve \delta}_\beta {\breve \delta}%
_\alpha = {\breve w}_{\beta \gamma }^\alpha \left({\breve u}\right) {\breve
\delta}_\alpha .
\end{equation}

The \textbf{algebra of tensorial distinguished fields} $DT\left({\breve {%
\mathcal{E}}}\right) $ (d--fields, d--tensors, d--objects) on ${\breve {%
\mathcal{E}}}$ is introduced as the tensor algebra ${\breve {\mathcal{T}}}
=\{ {\breve {\mathcal{T}}}_{qs}^{pr}\}$ of the cv--bundle ${\breve {\mathcal{%
E}}}_{\left( d\right) }= \left(H{\breve {\mathcal{E}}}\oplus V{\breve {%
\mathcal{E}}},{\breve p}_d, {\breve {\mathcal{E}}} \right),$ where ${\breve p%
}_d:\ H{\breve {\mathcal{E}}}\oplus V{\breve {\mathcal{E}}}\rightarrow {%
\breve {\mathcal{E}}}.$

An element ${\breve{\mathbf{t}}}\in {\breve{\mathcal{T}}}_{qs}^{pr},$
d--tensor field of type $\left(
\begin{array}{cc}
p & r \\
q & s%
\end{array}%
\right) ,$ can be written in local form as
\begin{equation}
{\breve{\mathbf{t}}}={\breve{t}}%
_{j_{1}...j_{q}b_{1}...b_{r}}^{i_{1}...i_{p}a_{1}...a_{r}}\left( {\breve{u}}%
\right) {\breve{\delta}}_{i_{1}}\otimes ...\otimes {\breve{\delta}}%
_{i_{p}}\otimes {\breve{\partial}}_{a_{1}}\otimes ...\otimes {\breve{\partial%
}}_{a_{r}}\otimes {\breve{d}}^{j_{1}}\otimes ...\otimes {\breve{d}}%
^{j_{q}}\otimes {\breve{\delta}}^{b_{1}}...\otimes {\breve{\delta}}^{b_{r}}.
\notag
\end{equation}

We shall respectively use the denotations $\mathcal{X}\left({\breve E}%
\right) $ (or $\mathcal{X} \left( M\right) ),\ \Lambda ^p\left({\breve {%
\mathcal{E}}}\right) $ or $\left(\Lambda ^p\left( M\right) \right) $ and $%
\mathcal{F}\left( {\breve E} \right)$ (or $\mathcal{F}$ $\left( M\right) $)
for the module of d--vector fields on ${\breve {\mathcal{E}}}$ (or $M$), the
exterior algebra of p--forms on ${\breve {\mathcal{E}}}$\ (or $M)$ and the
set of real functions on ${\breve {\mathcal{E}}}$ (or $M).$

\paragraph{Anholonomic frames in hvc--bundles}

The anholnomic frames adapted to a N--connec\-ti\-on in hvc--bundle $\tilde{%
\mathcal{E}}$ are defined by the set of coefficients (\ref{nconho}); having
restricted the constructions to a vector (covector) shell we obtain some
generalizations of the formulas for corresponding N(or \v{N})--connection
elongation of partial derivatives defined by (\ref{dder}) (or (\ref{ddercv}%
)) and (\ref{ddif}) (or (\ref{ddifcv})).

We introduce the adapted partial derivatives (anholonomic N--frames, or
N--bases) in $\tilde{\mathcal{E}}$ by applying the coefficients (\ref{nconho}%
)
\begin{equation*}
{\tilde{\delta}}_{\alpha }=\frac{{\tilde{\delta}}}{\delta \tilde{u}^{\alpha }%
}=\left( \delta _{i},\delta _{a_{1}},\delta _{a_{2}},{\breve{\delta}}%
^{a_{3}},...,{\breve{\delta}}^{a_{z-1}},\partial _{a_{z}}\right) ,
\end{equation*}%
where
\begin{eqnarray*}
&&\delta _{i}=\partial _{i}-N_{i}^{~a_{1}}\partial
_{a_{1}}-N_{i}^{~a_{2}}\partial _{a_{2}}+N_{ia_{3}}{\breve{\partial}}%
^{a_{3}}-...+N_{ia_{z-1}}{\breve{\partial}}^{a_{z-1}}-N_{i}^{~a_{z}}\partial
_{a_{z}}, \\
&&\delta _{a_{1}}=\partial _{a_{1}}-N_{a_{1}}^{~a_{2}}\partial
_{a_{2}}+N_{a_{1}a_{3}}{\breve{\partial}}^{a_{3}}-...+N_{a_{1}a_{z-1}}{%
\breve{\partial}}^{a_{z-1}}-N_{a_{1}}^{~a_{z}}\partial _{a_{z}}, \\
&&\delta _{a_{2}}=\partial _{a_{2}}+N_{a_{2}a_{3}}{\breve{\partial}}%
^{a_{3}}-...+N_{a_{2}a_{z-1}}{\breve{\partial}}^{a_{z-1}}-N_{a_{2}}^{~a_{z}}%
\partial _{a_{z}}, \\
&&{\breve{\delta}}^{a_{3}}={\tilde{\partial}}^{a_{3}}-N^{a_{3}a_{4}}\partial
_{a_{4}}-...+N_{~a_{z-1}}^{a_{3}}{\breve{\partial}}%
^{a_{z-1}}-N^{~a_{3}a_{z}}\partial _{a_{z}}, \\
&&................. \\
&&{\breve{\delta}}^{a_{z-1}}={\tilde{\partial}}^{a_{z-1}}-N^{~a_{z-1}a_{z}}%
\partial _{a_{z}}, \\
&&\partial _{a_{z}}=\partial /\partial y^{a_{z}}.
\end{eqnarray*}%
These formulas can be written in the matrix form:
\begin{equation}
{\tilde{\delta}}_{_{\bullet }}=\widehat{\mathbf{N}}(u)\times {\tilde{\partial%
}}_{_{\bullet }}  \label{dderho}
\end{equation}%
where
\begin{equation*}
{\tilde{\delta}}_{_{\bullet }}=\left(
\begin{array}{l}
\delta _{i} \\
\delta _{a_{1}} \\
\delta _{a_{2}} \\
{\breve{\delta}}^{a_{3}} \\
... \\
{\breve{\delta}}^{a_{z-1}} \\
\partial _{a_{z}}%
\end{array}%
\right) ,\quad {\tilde{\partial}}_{_{\bullet }}=\left(
\begin{array}{l}
\partial _{i} \\
\partial _{a_{1}} \\
\partial _{a_{2}} \\
{\tilde{\partial}}^{a_{3}} \\
... \\
{\tilde{\partial}}^{a_{z-1}} \\
\partial _{a_{z}}%
\end{array}%
\right) ,\quad
\end{equation*}%
and
\begin{equation*}
\widehat{\mathbf{N}}=\left(
\begin{array}{llllllll}
1 & -N_{i}^{~a_{1}} & -N_{i}^{~a_{2}} & N_{ia_{3}} & -N_{i}^{~a_{4}} & ... &
N_{ia_{z-1}} & -N_{i}^{~a_{z}} \\
0 & 1 & -N_{a_{1}}^{~a_{2}} & N_{a_{1}a_{3}} & -N_{a_{1}}^{~a_{4}} & ... &
N_{a_{1}a_{z-1}} & -N_{a_{1}}^{~a_{z}} \\
0 & 0 & 1 & N_{a_{2}a_{3}} & -N_{a_{2}}^{~a_{4}} & ... & N_{a_{2}a_{z-1}} &
-N_{a_{2}}^{~a_{z}} \\
0 & 0 & 0 & 1 & -N^{a_{3}a_{4}} & ... & N_{~a_{z-1}}^{a_{3}} &
-N^{~a_{3}a_{z}} \\
... & ... & ... & ... & ... & ... & ... & ... \\
0 & 0 & 0 & 0 & 0 & ... & 1 & -N^{~a_{z-1}a_{z}} \\
0 & 0 & 0 & 0 & 0 & ... & 0 & 1%
\end{array}%
\right) .
\end{equation*}

The adapted differentials (anholonomic N--coframes, or N--cobases) in $%
\tilde{\mathcal{E}}$ are introduced in the simplest form by using the matrix
formalism: The respective dual matrices
\begin{eqnarray*}
{\tilde{\delta}}^{\bullet } &=&\{{\tilde{\delta}}^\alpha \}=\left(
\begin{array}{lllllll}
d^i & \delta ^{a_1} & \delta ^{a_2} & {\breve{\delta}}_{a_3} & ... & {\breve{%
\delta}}_{a_{z-1}} & \delta ^{a_z}%
\end{array}
\right) , \\
{\tilde{d}}^{\bullet } &=&\{{\tilde{\partial}}^\alpha \}=\left(
\begin{array}{lllllll}
d^i & d^{a_1} & d^{a_2} & d_{a_3} & ... & {d}_{a_{z-1}} & d^{a_z}%
\end{array}
\right)
\end{eqnarray*}
are related via a matrix relation
\begin{equation}
{\tilde{\delta}}^{\bullet }={\tilde{d}}^{\bullet }\widehat{\mathbf{M}}
\label{ddifho}
\end{equation}
which defines the formulas for anholonomic N--coframes. The matrix $\widehat{%
\mathbf{M}}$ from (\ref{ddifho}) is the inverse to $\widehat{\mathbf{N}},$
i. e. satisfies the condition
\begin{equation}
\widehat{\mathbf{M}}\times \widehat{\mathbf{N}}=I.  \label{invmatr}
\end{equation}

The\textbf{\ anholonomic coefficients, } $\widetilde{\mathbf{w}}=\{%
\widetilde{w}_{\beta \gamma }^\alpha \left( \widetilde{u}\right) \},$ on
hcv--bundle $\tilde{\mathcal{E}}$ are expressed via coefficients of the
matrix $\widehat{\mathbf{N}}$ and their partial derivatives following the
relations
\begin{equation}
\left[ \widetilde{\delta }_\alpha ,\widetilde{\delta }_\beta \right] =%
\widetilde{\delta }_\alpha \widetilde{\delta }_\beta -\widetilde{\delta }%
_\beta \widetilde{\delta }_\alpha =\widetilde{w}_{\beta \gamma }^\alpha
\left( \widetilde{u}\right) \widetilde{\delta }_\alpha .  \label{anholho}
\end{equation}
We omit the explicit formulas on shells.

A d--tensor formalism can be also developed on the space $\tilde{\mathcal{E}}%
.$ In this case the indices have to be stipulated for every shell
separately, like for v--bundles or cv--bundles.

\section{Distinguished connections and metrics}

In general, distinguished objects (d--objects) on a v--bundle $\mathcal{E}$
(or cv--bundle ${\breve{\mathcal{E}}}$) are introduced as geometric objects
with various group and coordinate transforms coordinated with the
N--connection structure on $\mathcal{E}$ (or ${\breve{\mathcal{E}}}$). For
example, a distinguished connection (in brief, \textbf{d--connection)} $D$
on $\mathcal{E}$ (or ${\breve{\mathcal{E}}}$) is defined as a linear
connection $D$ on $E$ (or $\breve{E}$) conserving under a parallelism the
global decomposition (\ref{whitney}) (or (\ref{whitneyc})) into horizontal
and vertical subbundles of $T\mathcal{E}$ (or $T{\breve{\mathcal{E}}}).$ A
covariant derivation associated to a d--connection becomes d--covariant. We
shall give necessary formulas for cv--bundles in round brackets.

\subsection{D--connections}

\subsubsection{D--connections in v--bundles (cv--bundles)}

A N--connection in a v--bundle $\mathcal{E}$ (cv--bundle ${\breve{\mathcal{E}%
}}$) induces a corresponding decomposition of d--tensors into sums of
horizontal and vertical parts, for example, for every d--vector $X\in
\mathcal{X\left( E\right) }$ (${\breve X}\in \mathcal{X}\left( {\breve{%
\mathcal{E}}}\right) $ ) and 1--form $A\in \Lambda ^1\left( \mathcal{E}%
\right) $ ($\breve{A}\in \Lambda ^1\left( {\breve{\mathcal{E}}}\right) $) we
have respectively
\begin{eqnarray}
X &=&hX+vX\mathbf{\ \quad }\mbox{and \quad }A=hA+vA,  \label{vdecomp} \\
(\breve{X} &=&h\breve{X}+vX\mathbf{\ \quad }\mbox{and \quad }\breve{A}=h%
\breve{A}+v\breve{A})  \notag
\end{eqnarray}
where
\begin{equation*}
hX = X^i\delta _i,vX=X^a\partial _a \ (h\breve{X} = \breve{X}^i{\tilde \delta%
}_i,v\breve{X}=\breve{X}_a{\breve \partial}^a)
\end{equation*}
and
\begin{equation*}
hA =A_i\delta ^i,vA=A_ad^a \ (h\breve{A} = \breve{A}_i{\breve \delta}^i,v%
\breve{A}= \breve{A}^a{\breve d}_a).
\end{equation*}

In consequence, we can associate to every d--covariant derivation along the%
\newline
d--vector (\ref{vdecomp}), $D_{X}=X\circ D$ ($D_{\breve{X}}=\breve{X}\circ D$%
) two new operators of h- and v--covariant derivations
\begin{eqnarray*}
D_{X}^{(h)}Y &=&D_{hX}Y\quad \mbox{ and \quad }D_{X}^{\left( v\right)
}Y=D_{vX}Y,\quad \forall Y{\ \in }\mathcal{X\left( E\right) } \\
(D_{\breve{X}}^{(h)}\breve{Y} &=&D_{h\breve{X}}\breve{Y}\quad
\mbox{ and
\quad }D_{\breve{X}}^{\left( v\right) }\breve{Y}=D_{v\breve{X}}\breve{Y}%
,\quad \forall \breve{Y}{\ \in }\mathcal{X}\left( {\breve{\mathcal{E}}}%
\right) )
\end{eqnarray*}%
for which the following conditions hold:
\begin{eqnarray}
D_{X}Y &=&D_{X}^{(h)}Y{\ +}D_{X}^{(v)}Y  \label{dcovpr} \\
(D_{\breve{X}}\breve{Y} &=&D_{\breve{X}}^{(h)}\breve{Y}{\ +}D_{\breve{X}%
}^{(v)}\breve{Y}),  \notag
\end{eqnarray}%
where
\begin{eqnarray*}
D_{X}^{(h)}f &=&(hX\mathbf{)}f\mbox{ \quad and\quad }D_{X}^{(v)}f=(vX\mathbf{%
)}f,\quad X,Y\mathbf{\in }\mathcal{X\left( E\right) },f\in \mathcal{F}\left(
M\right) \\
({\breve{D}}_{\breve{X}}^{(h)}f &=&(h\breve{X}\mathbf{)}f%
\mbox{ \quad
and\quad }{\breve{D}}_{\breve{X}}^{(v)}f=(v\breve{X}\mathbf{)}f,\quad \breve{%
X},\breve{Y}{\in }\mathcal{X}\left( {\breve{\mathcal{E}}}\right) ,f\in
\mathcal{F}\left( M\right) ).
\end{eqnarray*}

The components $\Gamma _{\beta \gamma }^\alpha $ ( ${\breve{\Gamma}}_{\beta
\gamma }^\alpha )$of a d--connection ${\breve{D}}_\alpha =({\breve{\delta}}%
_\alpha \circ D),$ locally adapted to the N---connection structure with
respect to the frames (\ref{dder}) and (\ref{ddif}) ((\ref{ddercv}) and (\ref%
{ddifcv})), are defined by the equations
\begin{equation*}
D_\alpha \delta _\beta =\Gamma _{\alpha \beta }^\gamma \delta _\gamma ~({%
\breve{D}}_\alpha {\breve{\delta}}_\beta ={\breve{\Gamma}}_{\alpha \beta
}^\gamma {\breve{\delta}}_\gamma ~),
\end{equation*}
from which one immediately follows
\begin{equation}
\Gamma _{\alpha \beta }^\gamma \left( u\right) =\left( D_\alpha \delta
_\beta \right) \circ \delta ^\gamma \quad ~({\breve{\Gamma}}_{\alpha \beta
}^\gamma \left( {\breve{u}}\right) =\left( {\breve{D}}_\alpha {\breve{\delta}%
}_\beta \right) \circ {\breve{\delta}}^\gamma ).  \label{gamma}
\end{equation}

The coefficients of operators of h- and v--covariant derivations,
\begin{eqnarray*}
D_k^{(h)} &=&\{L_{jk}^i,L^a_{bk\;}\}\mbox{ and }D_c^{(v)}=%
\{C_{jk}^i,C_{bc}^a\} \\
({\breve{D}}_k^{(h)} &=&\{{\breve{L}}_{jk}^i,{\breve{L}}_{a k}^{~b}\}%
\mbox{
and }{\breve{D}}^{(v)c}=\{{\breve C}_{~j}^{i~c},{\breve C}_a^{~bc}\})
\end{eqnarray*}
(see (\ref{dcovpr})), are introduced as corresponding h- and
v--paramet\-ri\-za\-ti\-ons of (\ref{gamma})
\begin{eqnarray}
L_{jk}^i &=&\left( D_k\delta _j\right) \circ d^i,\quad L_{bk}^a=\left(
D_k\partial _b\right) \circ \delta ^a  \label{hgamma} \\
({\breve{L}}_{jk}^i &=&\left( {\breve{D}}_k{\breve{\delta}}_j\right) \circ
d^i,\quad {\breve{L}}_{a k}^{~b}=\left( {\breve{D}}_k{\breve \partial}%
^b\right) \circ {\breve{\delta}}_a)  \notag
\end{eqnarray}
and
\begin{eqnarray}
C_{jc}^i &=&\left( D_c\delta _j\right) \circ d^i,\quad C_{bc}^a=\left(
D_c\partial _b\right) \circ \delta ^a  \label{vgamma} \\
({\breve C}_{~j}^{i~c} &=&\left( {\breve{D}}^c{\breve{\delta}}_j\right)
\circ d^i,\quad {\breve C}_a^{~bc}=\left( {\breve{D}}^c{\breve \partial}%
^b\right) \circ {\breve{\delta}}_a).  \notag
\end{eqnarray}

A set of components (\ref{hgamma}) and (\ref{vgamma}) \
\begin{equation*}
\Gamma _{\alpha \beta }^\gamma =[L_{jk}^i,L_{bk}^a,C_{jc}^i,C_{bc}^a]~\left(
{\breve{\Gamma}}_{\alpha \beta }^\gamma =[{\breve{L}}_{jk}^i,{\breve{L}}%
_{ak}^{~b},{\breve{C}}_{~j}^{i~c},{\breve{C}}_a^{~bc}]\right)
\end{equation*}
completely defines the local action of a d---connection $D$ in $\mathcal{E}$
(${\breve{D}}$ in ${\breve{\mathcal{E}}).}$

For instance, having taken on $\mathcal{E}$ (${\breve{\mathcal{E}})}$ a
d---tensor field of type $\left(
\begin{array}{cc}
1 & 1 \\
1 & 1%
\end{array}%
\right) ,$
\begin{equation*}
\mathbf{t}=t_{jb}^{ia}\delta _{i}\otimes \partial _{a}\otimes d^{j}\otimes
\delta ^{b},{{\tilde{\mathbf{t}}}}={\breve{t}}_{ja}^{ib}{\breve{\delta}}%
_{i}\otimes {\breve{\partial}}^{a}\otimes d^{j}\otimes {\breve{\delta}}_{b},
\end{equation*}%
and a d--vector $\mathbf{X}$ ($\mathbf{{\breve{X}}})$ we obtain
\begin{eqnarray*}
D_{X}\mathbf{t} &=&D_{X}^{(h)}\mathbf{t+}D_{X}^{(v)}\mathbf{t=}\left( X^{k}{%
\breve{t}}_{jb|k}^{ia}+X^{c}t_{jb\perp c}^{ia}\right) \delta _{i}\otimes
\partial _{a}\otimes d^{j}\otimes \delta ^{b}, \\
({\breve{D}}_{\breve{X}}{{\tilde{\mathbf{t}}}} &=&\breve{D}_{\breve{X}}^{(h)}%
{\tilde{\mathbf{t}}}+\breve{D}_{\breve{X}}^{(v)}{\tilde{\mathbf{t}}}=\left(
\breve{X}^{k}{\breve{t}}_{ja|k}^{ib}+\breve{X}_{c}{\breve{t}}_{ja}^{ib\perp
c}\right) {\breve{\delta}}_{i}\otimes {\breve{\partial}}^{a}\otimes
d^{j}\otimes {\breve{\delta}}_{b})
\end{eqnarray*}%
where the h--covariant derivative is written
\begin{eqnarray*}
t_{jb|k}^{ia} &=&\delta
_{k}t_{jb}^{ia}+L_{hk}^{i}t_{jb}^{ha}+L_{ck}^{a}t_{jb}^{ic}
-L_{jk}^{h}t_{hb}^{ia}-L_{bk}^{c}t_{jc}^{ia} \\
({\breve{t}}_{ja|k}^{ib} &=&{\breve{\delta}}_{k}{\breve{t}}_{ja}^{ib}+{%
\breve{L}}_{hk}^{i}{\breve{t}}_{ja}^{hb}+{\breve{L}}_{ck}^{~b}{\breve{t}}%
_{ja}^{ic}-{\breve{L}}_{jk}^{h}{\breve{t}}_{ha}^{ib}-{\breve{L}}_{ck}^{~b}{%
\breve{t}}_{ja}^{ic})
\end{eqnarray*}%
and the v-covariant derivative is written
\begin{eqnarray}
t_{jb\perp c}^{ia} &=&\partial
_{c}t_{jb}^{ia}+C_{hc}^{i}t_{jb}^{ha}+C_{dc}^{a}t_{jb}^{id}
-C_{jc}^{h}t_{hb}^{ia}-C_{bc}^{d}t_{jd}^{ia}  \label{covder1} \\
({\breve{t}}_{ja}^{ib\perp c} &=&{\breve{\partial}}^{c}{\breve{t}}_{ja}^{ib}+%
{\breve{C}}_{~j}^{i~c}{\breve{t}}_{ja}^{hb}+{\breve{C}}_{a}^{~dc}{\breve{t}}%
_{jd}^{ib}-{\breve{C}}_{~j}^{i~c}{\breve{t}}_{ha}^{ib}-{\breve{C}}_{d}^{~bc}{%
\breve{t}}_{ja}^{id}).  \label{covder2}
\end{eqnarray}%
For a scalar function $f\in \mathcal{F}\left( \mathcal{E}\right) $ ( $f\in
\mathcal{F}\left( {\breve{\mathcal{E}}}\right) $) we have
\begin{eqnarray*}
D_{k}^{(h)} &=&\frac{\delta f}{\delta x^{k}}=\frac{\partial f}{\partial x^{k}%
}-N_{k}^{a}\frac{\partial f}{\partial y^{a}}\mbox{ and }D_{c}^{(v)}f=\frac{%
\partial f}{\partial y^{c}} \\
(\breve{D}_{k}^{(h)} &=&\frac{{\breve{\delta}}f}{\delta x^{k}}=\frac{%
\partial f}{\partial x^{k}}+N_{ka}\frac{\partial f}{\partial p_{a}}%
\mbox{
and }\breve{D}^{(v)c}f=\frac{\partial f}{\partial p_{c}}).
\end{eqnarray*}

\subsubsection{D--connections in hvc--bundles}

The theory of connections in higher order anisotropic vector superbundles
and vector bundles was elaborated in Refs. \cite{vstr2,vhsp,vbook}. Here we
re--formulate that formalism for the case when some shells of higher order
anisotropy could be covector spaces by stating the general rules of
covariant derivation compatible with the N--connection structure in
hvc--bundle $\tilde{\mathcal{E}}$ and omit details and cumbersome formulas.

For a hvc--bundle of type ${\tilde{\mathcal{E}}}={\tilde{\mathcal{E}}}%
[v(1),v(2),cv(3),...,cv(z-1),v(z)]$ a d--connection ${\tilde{\Gamma}}%
_{\alpha \beta }^{\gamma }$ has the next shell decomposition of components
(on induction being on the $p$-th shell, considered as the base space, which
in this case a hvc--bundle, we introduce in a usual manner, like a vector or
covector fiber, the $(p+1)$-th shell)%
\begin{eqnarray*}
{\tilde{\Gamma}}_{\alpha \beta }^{\gamma } &=&\{\Gamma _{\alpha _{1}\beta
_{1}}^{\gamma _{1}}=[L_{j_{1}k_{1}}^{i_{1}},L_{b_{1}k_{1}}^{a_{1}},
C_{j_{1}c_{1}}^{i_{1}},C_{b_{1}c_{1}}^{a_{1}}],\Gamma _{\alpha _{2}\beta
_{2}}^{\gamma _{2}}=[L_{j_{2}k_{2}}^{i_{2}},L_{b_{2}k_{2}}^{a_{2}},
C_{j_{2}c_{2}}^{i_{2}},C_{b_{2}c_{2}}^{a_{2}}], \\
\ {\breve{\Gamma}}_{\alpha _{3}\beta _{3}}^{\gamma _{3}} &=&[{\breve{L}}%
_{j_{3}k_{3}}^{i_{3}},{\breve{L}}_{a_{3}k_{3}}^{~b_{3}},{\breve{C}}%
_{~j_{3}}^{i_{3}~c_{3}},{\breve{C}}_{a_{3}}^{~b_{3}c_{3}}],..., \\
{\breve{\Gamma}}_{\alpha _{z-1}\beta _{z-1}}^{\gamma _{z-1}} &=&[{\breve{L}}%
_{j_{z-1}k_{z-1}}^{i_{z-1}},{\breve{L}}_{a_{z-1}k_{z-1}}^{~b_{z-1}},{\breve{C%
}}_{~j_{z-1}}^{i_{z-1}~c_{z-1}},{\breve{C}}_{a_{z-1}}^{~b_{z-1}c_{z-1}}],
\Gamma _{\alpha _{z}\beta _{z}}^{\gamma
_{z}}=[L_{j_{z}k_{z}}^{i_{z}},L_{b_{z}k_{z}}^{a_{z}},
C_{j_{z}c_{z}}^{i_{z}},C_{b_{z}c_{z}}^{a_{z}}]\}.
\end{eqnarray*}%
These coefficients determine the rules of a covariant derivation $\tilde{D}$
on ${\tilde{\mathcal{E}}}.$

For example, let us consider a d--tensor ${\tilde{\mathbf{t}}}$ of type $%
\left(
\begin{array}{llllll}
1 & 1_{1} & 1_{2} & {\breve{1}}_{3} & ... & 1_{z} \\
1 & 1_{1} & 1_{2} & {\breve{1}}_{3} & ... & 1_{z}%
\end{array}%
\right)$ with corresponding tensor product of components of anholonomic
N--frames (\ref{dderho}) and (\ref{ddifho})
\begin{equation*}
{\tilde{\mathbf{t}}}={\tilde{t}}_{jb_{1}b_{2}{\breve{a}}_{3}...{\breve{a}}%
_{z-1}b_{z}}^{ia_{1}a_{2}{\breve{b}}_{3}...{\breve{b}}_{z-1}a_{z}}\delta
_{i}\otimes \partial _{a_{1}}\otimes d^{j}\otimes \delta ^{b_{1}}\otimes
\partial _{a_{2}}\otimes \delta ^{b_{2}}\otimes {\breve{\partial}}%
^{a_{3}}\otimes {\breve{\delta}}_{b_{3}}{...}\otimes {\breve{\partial}}%
^{a_{z-1}}\otimes {\breve{\delta}}_{bz-1}\otimes \partial _{a_{z}}\otimes
\delta ^{b_{z}}.
\end{equation*}%
The d--covariant derivation $\tilde{D}$ of ${\tilde{\mathbf{t}}}$ is to be
performed separately for every shall according the rule (\ref{covder1}) if a
shell is defined by a vector subspace, or according the rule (\ref{covder2})
if the shell is defined by a covector subspace.

\subsection{Metric structure}

\subsubsection{D--metrics in v--bundles}

We define a \textbf{metric structure }$\mathbf{G\ }$
\index{metric} in the total space $E$ of a v--bundle $\mathcal{E=}$ $\left(
E,p,M\right) $ over a connected and paracompact base $M$ as a symmetric
covariant tensor field of type $\left( 0,2\right) $,
\begin{equation*}
\mathbf{G} = G_{\alpha \beta } du^{\alpha}\otimes du^\beta
\end{equation*}
being non degenerate and of constant signature on $E.$

Nonlinear connection $\mathbf{N}$ and metric $\mathbf{G}$ structures on $%
\mathcal{E}$ are mutually compatible it there are satisfied the conditions:
\begin{equation}  \label{comp}
\mathbf{G}\left( \delta _i,\partial _a\right) =0,\mbox{or equivalently, }%
G_{ia}\left( u\right) -N_i^b\left( u\right) h_{ab}\left( u\right) =0,
\end{equation}
where $h_{ab}=\mathbf{G}\left( \partial _a,\partial _b\right) $ and $G_{ia}=%
\mathbf{G}\left( \partial _i,\partial _a\right),$ which gives
\begin{equation}  \label{ncon1}
N_i^b\left( u\right) = h^{ab}\left( u\right) G_{ia}\left( u\right)
\end{equation}
( the matrix $h^{ab}$ is inverse to $h_{ab}).$ One obtains the following
decomposition of metric:
\begin{equation}  \label{metrdec}
\mathbf{G}(X,Y)\mathbf{=hG}(X,Y)+\mathbf{vG}(X,Y),
\end{equation}
where the d--tensor $\mathbf{hG}(X,Y)$ = $\mathbf{\ G}(hX,hY)$ is of type $%
\left(
\begin{array}{cc}
0 & 0 \\
2 & 0%
\end{array}
\right) $ and the d--tensor\newline
$\mathbf{vG}(X,Y) = \mathbf{G}(vX,vY)$ is of type $\left(
\begin{array}{cc}
0 & 0 \\
0 & 2%
\end{array}
\right) .$ With respect to the anholonomic basis (\ref{dder}) the d--metric (%
\ref{metrdec}) is written
\begin{equation}  \label{dmetric}
\mathbf{G}=g_{\alpha \beta }\left( u\right) \delta ^\alpha \otimes \delta
^\beta =g_{ij}\left( u\right) d^i\otimes d^j+h_{ab}\left( u\right) \delta
^a\otimes \delta ^b,
\end{equation}
where $g_{ij}=\mathbf{G}\left( \delta _i,\delta _j\right) .$

A metric structure of type (\ref{metrdec}) (equivalently, of type (\ref%
{dmetric})) or a metric on $E$ with components satisfying the constraints (%
\ref{comp}), (equivalently (\ref{ncon1})) defines an adapted to the given
N--connection inner (d--scalar) product on the tangent bundle $\mathcal{TE}$.

A d--connection $D_{X}$ is \textbf{metric} (or \textbf{compatible } with
metric $\mathbf{G}$) on $\mathcal{E}$ if $D_{X}\mathbf{G}=0,\forall X\mathbf{%
\in }\mathcal{X\left( E\right)}.$ With respect to anholonomic frames these
conditions are written
\begin{equation}
D_{\alpha }g_{\beta \gamma }=0,  \label{comatib}
\end{equation}%
where by $g_{\beta \gamma }$ we denote the coefficients in the block form (%
\ref{dmetric}).

\subsubsection{D--metrics in cv-- and hvc--bundles}

The presented considerations on self--consistent definition of
N--connection, d--connection and metric structures in v--bundles can
re--formulated in a similar fashion for another types of anisotropic
space--times, on cv--bundles and on shells of hvc--bundles. For simplicity,
we give here only the analogous formulas for the metric d--tensor (\ref%
{dmetric}):

\begin{itemize}
\item On cv--bundle ${%
\breve{\mathcal{E}}}$ we write
\begin{equation}
{\breve{\mathbf{G}}}={\breve{g}}_{\alpha \beta }\left( {\breve{u}}\right) {%
\breve{\delta}}^\alpha \otimes {\breve{\delta}}^\beta ={\breve{g}}%
_{ij}\left( {\breve{u}}\right) d^i\otimes d^j+{\breve{h}}^{ab}\left( {%
\breve{u}}\right) {\breve{\delta}}_a\otimes {\breve{\delta}}_b,
\label{dmetricvc}
\end{equation}
where ${\breve{g}}_{ij}={\breve{\mathbf{G}}}\left( {\breve{\delta}}_i,{%
\breve{\delta}}_j\right) $ and ${\breve{h}}^{ab}={\breve{\mathbf{G}}}\left( {%
\breve{\partial}}^a,{\breve{\partial}}^b\right) $ and the N--coframes are
given by formulas (\ref{ddifcv}).

For simplicity, we shall consider that the metricity conditions are
satisfied, ${\breve{D}}_\gamma {\breve{g}}_{\alpha \beta}=0.$

\item On hvc--bundle ${\tilde{\mathcal{E}}}$ we write
\begin{eqnarray}
{\tilde{\mathbf{G}}} &=&{\tilde{g}}_{\alpha \beta }\left( {\tilde{u}}\right)
{\tilde{\delta}}^{\alpha }\otimes {\tilde{\delta}}^{\beta }={\tilde{g}}%
_{ij}\left( {\tilde{u}}\right) d^{i}\otimes d^{j}+{\tilde{h}}%
_{a_{1}b_{1}}\left( {\tilde{u}}\right) {\delta }^{a_{1}}\otimes {\delta }%
^{b_{1}}+{\tilde{h}}_{a_{2}b_{2}}\left( {\tilde{u}}\right) {\delta }%
^{a_{2}}\otimes {\delta }^{b_{2}}  \label{dmetrichcv} \\
&&+{\tilde{h}}^{a_{3}b_{3}}\left( {\tilde{u}}\right) {\breve{\delta}}%
_{a_{3}}\otimes {\breve{\delta}}_{b_{3}}+...+{\tilde{h}}^{a_{z-1}b_{z-1}}%
\left( {\tilde{u}}\right) {\breve{\delta}}_{a_{z-1}}\otimes {\breve{\delta}}%
_{b_{z-1}}+{\tilde{h}}_{a_{z}b_{z}}\left( {\tilde{u}}\right) {\delta }%
^{a_{z}}\otimes {\delta }^{b_{z}},  \notag
\end{eqnarray}%
where ${\tilde{g}}_{ij}={\tilde{\mathbf{G}}}\left( {\tilde{\delta}}_{i},{%
\tilde{\delta}}_{j}\right) $ and ${\tilde{h}}_{a_{1}b_{1}}={\tilde{\mathbf{G}%
}}\left( \partial _{a_{1}},\partial _{b_{1}}\right) ,$ ${\tilde{h}}%
_{a_{2}b_{2}}={\tilde{\mathbf{G}}}\left( \partial _{a_{2}},\partial
_{b_{2}}\right) ,$\newline
${\tilde{h}}^{a_{3}b_{3}}={\tilde{\mathbf{G}}}({\breve{\partial}}^{a_{3}}, {%
\breve{\partial}}^{b_{3}}),$ ... and the N--coframes are given by formulas (%
\ref{ddifho}).

The metricity conditions are ${\tilde{D}}_{\gamma }{\tilde{g}}_{\alpha \beta
}=0.$

\item On osculator bundle $T^{2}M=Osc^{2}M$, we have a particular case of (%
\ref{dmetrichcv}) when
\begin{eqnarray}
{\tilde{\mathbf{G}}} &=&{\tilde{g}}_{\alpha \beta }\left( {\tilde{u}}\right)
{\tilde{\delta}}^{\alpha }\otimes {\tilde{\delta}}^{\beta }
\label{dmetrichosc2} \\
&=&{\tilde{g}}_{ij}\left( {\tilde{u}}\right) d^{i}\otimes d^{j}+{\tilde{h}}%
_{ij}\left( {\tilde{u}}\right) {\delta y}_{(1)}^{i}\otimes {\delta y}%
_{(1)}^{i}+{\tilde{h}}_{ij}\left( {\tilde{u}}\right) {\delta y}%
_{(2)}^{i}\otimes {\delta y}_{(2)}^{i}  \notag
\end{eqnarray}
with respect to N--coframes.

\item On dual osculator bundle $\left( T^{\ast 2}M,p^{\ast 2},M\right) $ we
have another particular case of (\ref{dmetrichcv}) when
\begin{eqnarray}
{\tilde{\mathbf{G}}} &=&{\tilde{g}}_{\alpha \beta }\left( {\tilde{u}}\right)
{\tilde{\delta}}^{\alpha }\otimes {\tilde{\delta}}^{\beta }
\label{dmetrichosc2d} \\
&=&{\tilde{g}}_{ij}\left( {\tilde{u}}\right) d^{i}\otimes d^{j}+{\tilde{h}}%
_{ij}\left( {\tilde{u}}\right) {\delta y}_{(1)}^{i}\otimes {\delta y}%
_{(1)}^{i}+{\tilde{h}}^{ij}\left( {\tilde{u}}\right) {\delta p}%
_{i}^{(2)}\otimes {\delta p}_{i}^{(2)}  \notag
\end{eqnarray}
with respect to N--coframes.
\end{itemize}

\subsection{Some examples of d--connections}

We emphasize that the geometry of connections in a v--bundle $\mathcal{E}$
is very reach. If a triple of fundamental geometric objects $\left(
N_i^a\left( u\right) ,\Gamma _{\beta \gamma }^\alpha \left( u\right)
,g_{\alpha \beta }\left( u\right) \right) $ is fixed on $\mathcal{E}$, a
multi--connection structure (with corresponding different rules of covariant
derivation, which are, or not, mutually compatible and with the same, or
not, induced d--scalar products in $\mathcal{TE)}$ is defined on this
v--bundle. We can give a priority to a connection structure following some
physical arguments, like the reduction to the Christoffel symbols in the
holonomic case, mutual compatibility between metric and N--connection and
d--connection structures and so on.

In this subsection we enumerate some of the connections and covariant
derivations in v--bundle $\mathcal{E}$, cv--bundle ${\breve{\mathcal{E}}}$
and in some hvc--bundles which can present interest in investigation of
locally anisotropic gravitational and matter field interactions :

\begin{enumerate}
\item Every N--connection in $\mathcal{E}$ with coefficients $N_i^a\left(
x,y\right) $ being differentiable on y--variables, induces a structure of
linear connection $N_{\beta \gamma }^\alpha ,$ where
\begin{equation}
N_{bi}^a=\frac{\partial N_i^a}{\partial y^b}\mbox{ and
}N_{bc}^a\left( x,y\right) =0.  \label{nlinearized}
\end{equation}
For some $Y\left( u\right) =Y^i\left( u\right) \partial _i+Y^a\left(
u\right) \partial _a$ and $B\left( u\right) =B^a\left( u\right) \partial _a$
one introduces a covariant derivation as
\begin{equation*}
D_Y^{(\widetilde{N})}B=\left[ Y^i\left( \frac{\partial B^a}{\partial x^i}%
+N_{bi}^aB^b\right) +Y^b\frac{\partial B^a}{\partial y^b}\right] \frac
\partial {\partial y^a}.
\end{equation*}

\item The d--connection of Berwald type on v--bundle $\mathcal{E} $
(cv--bundle ${\breve{\mathcal{E}})}$
\begin{eqnarray}
\Gamma _{\beta \gamma }^{(B)\alpha } &=&\left( L_{jk}^i,\frac{\partial N_k^a%
}{\partial y^b},0,C_{bc}^a\right) ,  \label{berwald} \\
({\breve{\Gamma}}_{\beta \gamma }^{(B)\alpha } &=&\left( \breve{L}_{jk}^i,-%
\frac{\partial \breve{N}_{ka}}{\partial p_b},0,{\breve{C}}_a^{~bc}\right) )
\notag
\end{eqnarray}
where 
\begin{eqnarray}
L_{.jk}^i\left( x,y\right) &=&\frac 12g^{ir}\left( \frac{\delta g_{jk}}{%
\delta x^k}+\frac{\delta g_{kr}}{\delta x^j}-\frac{\delta g_{jk}}{\delta x^r}%
\right) ,  \label{lccoef} \\
C_{.bc}^a\left( x,y\right) &=&\frac 12h^{ad}\left( \frac{\partial h_{bd}}{%
\partial y^c}+\frac{\partial h_{cd}}{\partial y^b}-\frac{\partial h_{bc}}{%
\partial y^d}\right)  \notag \\
(\breve{L}_{.jk}^i\left( x,p\right) &=&\frac 12\breve{g}^{ir}\left( \frac{{%
\breve{\delta}}\breve{g}_{jk}}{\delta x^k}+\frac{{\breve{\delta}}\breve{g}%
_{kr}}{\delta x^j}-\frac{{\breve{\delta}}\breve{g}_{jk}}{\delta x^r}\right) ,
\notag \\
{\breve{C}}_a^{~bc}\left( x,p\right) &=&\frac 12\breve{h}_{ad}\left( \frac{%
\partial \breve{h}^{bd}}{\partial p_c}+\frac{\partial \breve{h}^{cd}}{%
\partial p_b}-\frac{\partial \breve{h}^{bc}}{\partial p_d}\right) ),  \notag
\end{eqnarray}
which is hv---metric, i.e. there are satisfied the conditions $%
D_k^{(B)}g_{ij}=0$ and $D_c^{(B)}h_{ab}=0$ ($\breve{D}_k^{(B)}\breve{g}%
_{ij}=0$ and $\breve{D}^{(B)c}\breve{h}^{ab}=0).$

\item The canonical d--connection $\mathbf{\Gamma ^{(c)}}$ (or $\mathbf{%
\breve{\Gamma}^{(c)})}$ on a v--bundle (or cv--bundle) is associated to a
metric $\mathbf{G}$ (or $\mathbf{\breve{G})}$ of type (\ref{dmetric}) (or (%
\ref{dmetricvc})),
\begin{equation*}
\Gamma _{\beta \gamma }^{(c)\alpha
}=[L_{jk}^{(c)i},L_{bk}^{(c)a},C_{jc}^{(c)i},C_{bc}^{(c)a}]~(\breve{\Gamma}%
_{\beta \gamma }^{(c)\alpha }=[\breve{L}_{jk}^{(c)i},\breve{L}%
_{~a~.k}^{(c).b},\breve{C}_{~j}^{(c)i\ c},{\breve{C}}_a^{(c)~bc}])
\end{equation*}
with coefficients
\begin{eqnarray}
L_{jk}^{(c)i} &=&L_{.jk}^i,C_{bc}^{(c)a}=C_{.bc}^a~(\breve{L}_{jk}^{(c)i}=%
\breve{L}_{.jk}^i,{\breve{C}}_a^{(c)~bc}={\breve{C}}_a^{~bc}),%
\mbox{ (see
(\ref{lccoef})}  \notag \\
L_{bi}^{(c)a} &=&\frac{\partial N_i^a}{\partial y^b}+\frac 12h^{ac}\left(
\frac{\delta h_{bc}}{\delta x^i}-\frac{\partial N_i^d}{\partial y^b}h_{dc}-%
\frac{\partial N_i^d}{\partial y^c}h_{db}\right)  \notag \\
~(\breve{L}_{~a~.i}^{(c).b} &=&-\frac{\partial {\breve{N}}_i^a}{\partial p_b}%
+\frac 12\breve{h}_{ac}\left( \frac{{\breve{\delta}}\breve{h}^{bc}}{\delta
x^i}+\frac{\partial {\breve{N}}_{id}}{\partial p_b}\breve{h}^{dc}+\frac{%
\partial {\breve{N}}_{id}}{\partial p_c}\breve{h}^{db}\right) ),  \notag \\
~C_{jc}^{(c)i} &=&\frac 12g^{ik}\frac{\partial g_{jk}}{\partial y^c}~(\breve{%
C}_{~j}^{(c)i\ c}=\frac 12\breve{g}^{ik}\frac{\partial \breve{g}_{jk}}{%
\partial p_c}).  \label{inters}
\end{eqnarray}
This is a metric d--connection which satisfies conditions
\begin{eqnarray*}
D_k^{(c)}g_{ij} &=&0,D_c^{(c)}g_{ij}=0,D_k^{(c)}h_{ab}=0,D_c^{(c)}h_{ab}=0 \\
(\breve{D}_k^{(c)}\breve{g}_{jk} &=&0,\breve{D}^{(c)c}\breve{g}_{jk}=0,%
\breve{D}_k^{(c)}\breve{h}^{bc}=0,\breve{D}^{(c)c}\breve{h}^{ab}=0).
\end{eqnarray*}
In physical applications, we shall use the canonical connection and, for
simplicity, we shall omit the index $(c).$ The coefficients (\ref{inters})$%
\, $are to be extended to higher order if we are dealing with derivations of
geometrical objects with ''shell'' indices. In this case the fiber indices
are to be stipulated for every type of shell into consideration.

\item We can consider the N--adapted Christoffel symbols
\begin{equation}
\widetilde{\Gamma }_{\beta \gamma }^{\alpha }=\frac{1}{2}g^{\alpha \tau
}\left( \delta _{\gamma }g_{\tau \beta }+\delta _{\beta }g_{\tau \gamma
}-\delta g_{\beta \gamma }\right) ,  \label{dchrist}
\end{equation}%
which have the components of d--connection $\widetilde{\Gamma }_{\beta
\gamma }^{\alpha }=\left( L_{jk}^{i},0,0,C_{bc}^{a}\right) ,$ with $%
L_{jk}^{i}$ and $C_{bc}^{a}$ as in (\ref{lccoef}) if $g_{\alpha \beta }$ is
taken in the form (\ref{dmetric}).
\end{enumerate}

\subsection{Amost Hermitian anisotropic spaces}

The are possible very interesting particular constructions \cite{ma87,mhss}
on t--bundle $TM$ provided with N--connection which defines a N--adapted
frame structure $\delta _{\alpha }=(\delta _{i},\dot{\partial}_{i})$ (for
the same formulas (\ref{dder}) and (\ref{ddif}) but with identified fiber
and base indices). We are using the 'dot' symbol in order to distinguish the
horizontal and vertical operators because on t--bundles the indices could
take the same values both for the base and fiber objects. This allow us to
define an almost complex structure $\mathbf{J}=\{J_{\alpha }^{\ \beta }\}$
on $TM$ as follows
\begin{equation}
\mathbf{J}(\delta _{i})=-\dot{\partial}_{i},\ \mathbf{J}(\dot{\partial}%
_{i})=\delta _{i}.  \label{alcomp}
\end{equation}%
It is obvious that $\mathbf{J}$ is well--defined and $\mathbf{J}^{2}=-I.$

For d--metrics of type (\ref{dmetric}), on $TM,$ we can consider the case
when\newline
$g_{ij}(x,y)=h_{ab}(x,y),$ i. e.
\begin{equation}  \label{dmetrict}
\mathbf{G}_{(t)}= g_{ij}(x,y)dx^i\otimes dx^j + g_{ij}(x,y)\delta y^i\otimes
\delta y^j,
\end{equation}
where the index $(t)$ denotes that we have geometrical object defined on
tangent space.

An almost complex structure $J_{\alpha }^{\ \beta }$ is compatible with a
d--metric of type (\ref{dmetrict}) and a d--connection $D$ on tangent bundle
$TM$ if the conditions
\begin{equation*}
J_{\alpha }^{\ \beta }J_{\gamma }^{\ \delta }g_{\beta \delta }=g_{\alpha
\gamma }\ \mbox{ and }\ D_{\alpha }J_{\ \beta }^{\gamma }=0
\end{equation*}%
are satisfied.

The pair $(\mathbf{G}_{(t)},\mathbf{J})$ is an almost Hermitian structure on
$TM.$

One can introduce an almost sympletic 2--form associated to the almost
Hermitian structure $(\mathbf{G}_{(t)},\mathbf{J}),$
\begin{equation}  \label{hermit}
\theta = g_{ij}(x,y)\delta y^i\wedge dx^j.
\end{equation}

If the 2--form (\ref{hermit}), defined by the coefficients $g_{ij},$ is
closed, we obtain an almost K\"{a}hlerian structure in $TM.$

\begin{definition}
An almost K\"{a}hler metric connection is a linear connection $D^{(H)}$ on $T%
{\tilde{M}}=TM\setminus \{0\}$ with the properties:

\begin{enumerate}
\item $D^{(H)}$ preserve by parallelism the vertical distribution defined by
the N--connection structure;

\item $D^{(H)}$ is compatible with the almost K\"{a}hler structure $(\mathbf{%
G}_{(t)},\mathbf{J})$, i. e.
\begin{equation*}
D_{X}^{(H)}g=0,\ D_{X}^{(H)}J=0,\ \forall X\in \mathcal{X}\left( T{\tilde{M}}%
\right) .
\end{equation*}
\end{enumerate}
\end{definition}

By straightforward calculation we can prove that a d--connection \newline
$D\Gamma =\left(L_{jk}^i,L_{jk}^i,C_{jc}^i,C_{jc}^i\right) $ with the
coefficients defined by
\begin{equation}  \label{kahlerconm}
D^{(H)}_{\delta _i}\delta _j = L_{jk}^i \delta _i,\ D^{(H)}_{\delta _i}\dot{%
\partial}_j = L_{jk}^i \dot{\partial}_i, \ D^{(H)}_{\delta _i}\delta _j =
C_{jk}^i \delta _i,\ D^{(H)}_{\delta _i}\dot{\partial}_j = C_{jk}^i \dot{%
\partial}_i,
\end{equation}
where $L_{jk}^i$ and $C_{ab}^e \to C_{jk}^i,$\ on $TM$ are defined by the
formulas (\ref{lccoef}), define a torsionless (see the next section on
torsion structures) metric d--connection which satisfy the compatibility
conditions (\ref{comatib}).

Almost complex structures and almost K\"{a}hler models of Finsler, Lagrange,
Hamilton and Cartan geometries (of first an higher orders) are investigated
in details in Refs. \cite{mhss,vbook}.

\subsection{ Torsions and Curvatures}

We outline the basic definitions and formulas for the torsion and curvature
structures in v--bundles and cv--bundles provided with N--connection
structure.

\subsubsection{N--connection curvature}

\begin{enumerate}
\item The curvature $\mathbf{\Omega }$$\,$ of a nonlinear connection $%
\mathbf{N}$ in a v--bundle $\mathcal{E}$ can be defined in local form as $%
\mathbf{\ }$ \cite{ma87}:
\begin{equation*}
\mathbf{\Omega }=\frac{1}{2}\Omega _{ij}^{a}d^{i}\wedge d^{j}\otimes
\partial _{a},
\end{equation*}%
where
\begin{equation}
\Omega _{ij}^{a} =\delta _{j}N_{i}^{a}-\delta _{i}N_{j}^{a} =\partial
_{j}N_{i}^{a}-\partial _{i}N_{j}^{a}+N_{i}^{b}N_{bj}^{a}-N_{j}^{b}N_{bi}^{a},
\label{ncurv}
\end{equation}
$N_{bi}^{a}$ being that from (\ref{nlinearized}).

\item For the curvature $\mathbf{\breve{\Omega}},$ of a nonlinear connection
$\mathbf{\breve{N}}$ in a cv--bundle $\mathcal{\breve{E}}$ we introduce
\begin{equation*}
\mathbf{\breve{\Omega}}\,=\frac{1}{2}\breve{\Omega}_{ija}d^{i}\wedge
d^{j}\otimes \breve{\partial}^{a},
\end{equation*}%
where
\begin{eqnarray}
\breve{\Omega}_{ija} &=&-\breve{\delta}_{j}\breve{N}_{ia}+\breve{\delta}_{i}%
\breve{N}_{ja}=-\partial _{j}\breve{N}_{ia}+\partial _{i}\breve{N}_{ja}+%
\breve{N}_{ib}\breve{N}_{ja}^{\quad b}-\breve{N}_{jb}\breve{N}_{ja}^{\quad
b},  \notag \\
\breve{N}_{ja}^{\quad b} &=&\breve{\partial}^{b}\breve{N}_{ja}=\partial
\breve{N}_{ja}/\partial p_{b}.  \label{ncurvcv}
\end{eqnarray}

\item There were analyzed the curvatures $\mathbf{\tilde{\Omega}}$ of
different type of nonlinear connections $\mathbf{\tilde{N}}$ in higher order
an\-isot\-ro\-pic bundles were analyzed for higher order tangent/dual
tangent bundles and higher order prolongations of generalized Finsler,
Lagrange and Hamilton spaces in Refs. \cite{mhss} and for higher order
anisotropic superspaces and spinor bundles in Refs. \cite%
{vbook,vsp1,vhsp,vstr2}: For every higher order anisotropy shell, we shall
define the coefficients (\ref{ncurv}) or (\ref{ncurvcv}) in dependence of
the fact with type of subfiber we are considering (a vector or covector
fiber).
\end{enumerate}

\subsubsection{d--Torsions in v- and cv--bundles}

The torsion $\mathbf{T}$ of a d--connection $\mathbf{D\ }$ in v--bundle $%
\mathcal{E}$ (cv--bundle $\mathcal{\breve{E})}$ is defined by the equation
\begin{equation}
\mathbf{T\left( X,Y\right) =XY_{\circ }^{\circ }T\doteq }D_X\mathbf{Y-}D_Y%
\mathbf{X\ -\left[ X,Y\right] .}  \label{torsion}
\end{equation}
One holds the following h- and v--decompositions
\begin{equation*}
\mathbf{T\left( X,Y\right) =T\left( hX,hY\right) +T\left( hX,vY\right)
+T\left( vX,hY\right) +T\left( vX,vY\right) .}
\end{equation*}
We consider the projections:
\begin{equation*}
\mathbf{hT\left( X,Y\right) ,vT\left( hX,hY\right) ,hT\left( hX,hY\right)
,...}
\end{equation*}
and say that, for instance, $\mathbf{hT\left( hX,hY\right) }$ is the
h(hh)--torsion of $\mathbf{D}$ , $\mathbf{vT\left( hX,hY\right) }$ is the
v(hh)--torsion of $\mathbf{D}$ and so on.

The torsion (\ref{torsion}) in v-bundle is locally determined by five
d--tensor fields, torsions, defined as
\begin{eqnarray}
T_{jk}^{i} &=&\mathbf{hT}\left( \delta _{k},\delta _{j}\right) \cdot
d^{i},\quad T_{jk}^{a}=\mathbf{vT}\left( \delta _{k},\delta _{j}\right)
\cdot \delta ^{a},  \label{dtorsions} \\
P_{jb}^{i} &=&\mathbf{hT}\left( \partial _{b},\delta _{j}\right) \cdot
d^{i},\quad P_{jb}^{a}=\mathbf{vT}\left( \partial _{b},\delta _{j}\right)
\cdot \delta ^{a},S_{bc}^{a}=\mathbf{vT}\left( \partial _{c},\partial
_{b}\right) \cdot \delta ^{a}.  \notag
\end{eqnarray}%
Using formulas (\ref{dder}), (\ref{ddif}), (\ref{ncurv}) and (\ref{torsion})
we can compute \cite{ma87} in explicit form the components of torsions (\ref%
{dtorsions}) for a d--connection of type (\ref{hgamma}) and (\ref{vgamma}):
\begin{eqnarray}
T_{.jk}^{i} &=&T_{jk}^{i}=L_{jk}^{i}-L_{kj}^{i},\quad
T_{ja}^{i}=C_{.ja}^{i},T_{aj}^{i}=-C_{ja}^{i},T_{.ja}^{i}=0,\quad
T_{.ib}^{a}=-P_{.bi}^{a}  \label{dtorsc} \\
\qquad T_{.bc}^{a} &=&S_{.bc}^{a}=C_{bc}^{a}-C_{cb}^{a},T_{.ij}^{a}=\delta
_{j}N_{i}^{a}-\delta _{j}N_{j}^{a},\quad T_{.bi}^{a}=P_{.bi}^{a}=\partial
_{b}N_{i}^{a}-L_{.bj}^{a}.  \notag
\end{eqnarray}

Formulas similar to (\ref{dtorsions}) and (\ref{dtorsc}) hold for
cv--bundles:
\begin{eqnarray}
\check{T}_{jk}^{i} &=&\mathbf{hT}\left( \delta _{k},\delta _{j}\right) \cdot
d^{i},\quad \check{T}_{jka}=\mathbf{vT}\left( \delta _{k},\delta _{j}\right)
\cdot \check{\delta}_{a},  \label{torsionsa} \\
\check{P}_{j}^{i\quad b} &=&\mathbf{hT}\left( \check{\partial}^{b},\delta
_{j}\right) \cdot d^{i},\quad \check{P}_{aj}^{\quad b}=\mathbf{vT}\left(
\check{\partial}^{b},\delta _{j}\right) \cdot \check{\delta}_{a},\check{S}%
_{a}^{\quad bc}=\mathbf{vT}\left( \check{\partial}^{c},\check{\partial}%
^{b}\right) \cdot \check{\delta}_{a}.  \notag
\end{eqnarray}%
and
\begin{eqnarray}
\check{T}_{.jk}^{i} &=&\check{T}_{jk}^{i}=L_{jk}^{i}-L_{kj}^{i},\quad
\check{T}_{j}^{ia}=\check{C}_{.j}^{i\ a},\check{T}_{~\ j}^{ia}=-\check{C}%
_{j}^{i~a},~\check{T}_{.j}^{i~a}=0,\check{T}_{a~b}^{~j}=-\check{P}_{a\
b}^{~j}  \label{dtorsca} \\
\qquad \check{T}_{a}^{~bc} &=&\check{S}_{a}^{~bc}=\check{C}_{a}^{~bc}-%
\check{C}_{a}^{~cb},\check{T}_{.ija}=-\delta _{j}\check{N}_{ia}+\delta _{j}%
\check{N}_{ja},\quad \check{T}_{a}^{~bi}=\check{P}_{a}^{~bi}=-\check{\partial%
}^{b}\check{N}_{ia}-\check{L}_{a}^{~bi}.  \notag
\end{eqnarray}

The formulas for torsion can be generalized for hvc--bundles (on every shell
we must write (\ref{dtorsc}) or (\ref{dtorsca}) in dependence of the type of
shell, vector or co-vector one, we are dealing).

\subsubsection{d--Curvatures in v- and cv--bundles}

The curvature $\mathbf{R}$ of a d--connection in v--bundle $\mathcal{E}$ is
defined by the equation
\begin{equation*}
\mathbf{R}\left( X,Y\right) Z=XY_{\bullet }^{\bullet }R\bullet Z=D_{X}D_{Y}{%
\ Z}-D_{Y}D_{X}Z-D_{[X,Y]}{Z.}
\end{equation*}%
One holds the next properties for the h- and v--decompositions of curvature:
\begin{eqnarray}
\mathbf{vR}\left( X,Y\right) hZ &=&0,\ \mathbf{hR}\left( X,Y\right) vZ=0,
\label{curvaturehv} \\
\mathbf{R}\left( X,Y\right) Z &=&\mathbf{hR}\left( X,Y\right) hZ+\mathbf{vR}%
\left( X,Y\right) vZ.  \notag
\end{eqnarray}

From (\ref{curvaturehv}) and the equation $\mathbf{R\left( X,Y\right)
=-R\left( Y,X\right) }$ we get that the curvature of a d--con\-nec\-ti\-on $%
\mathbf{D}$ in $\mathcal{E}$ is completely determined by the following six
d--tensor fields:
\begin{eqnarray}
R_{h.jk}^{.i} &=&d^{i}\cdot \mathbf{R}\left( \delta _{k},\delta _{j}\right)
\delta _{h},~R_{b.jk}^{.a}=\delta ^{a}\cdot \mathbf{R}\left( \delta
_{k},\delta _{j}\right) \partial _{b},  \label{rps} \\
P_{j.kc}^{.i} &=&d^{i}\cdot \mathbf{R}\left( \partial _{c},\partial
_{k}\right) \delta _{j},~P_{b.kc}^{.a}=\delta ^{a}\cdot \mathbf{R}\left(
\partial _{c},\partial _{k}\right) \partial _{b},  \notag \\
S_{j.bc}^{.i} &=&d^{i}\cdot \mathbf{R}\left( \partial _{c},\partial
_{b}\right) \delta _{j},~S_{b.cd}^{.a}=\delta ^{a}\cdot \mathbf{R}\left(
\partial _{d},\partial _{c}\right) \partial _{b}.  \notag
\end{eqnarray}%
By a direct computation, using (\ref{dder}),(\ref{ddif}),(\ref{hgamma}),(\ref%
{vgamma}) and (\ref{rps}) we get:
\begin{eqnarray}
R_{h.jk}^{.i} &=&{\delta }_{h}L_{.hj}^{i}-{\delta }%
_{j}L_{.hk}^{i}+L_{.hj}^{m}L_{mk}^{i}-L_{.hk}^{m}L_{mj}^{i}+C_{.ha}^{i}R_{.jk}^{a},
\label{dcurvatures} \\
R_{b.jk}^{.a} &=&{\delta }_{k}L_{.bj}^{a}-{\delta }%
_{j}L_{.bk}^{a}+L_{.bj}^{c}L_{.ck}^{a}-L_{.bk}^{c}L_{.cj}^{a}+C_{.bc}^{a}R_{.jk}^{c},
\notag \\
P_{j.ka}^{.i} &=&{\partial }_{a}L_{.jk}^{i}-\left( {\delta }%
_{k}C_{.ja}^{i}+L_{.lk}^{i}C_{.ja}^{l}-L_{.jk}^{l}C_{.la}^{i}-L_{.ak}^{c}C_{.jc}^{i}\right) +C_{.jb}^{i}P_{.ka}^{b},
\notag \\
P_{b.ka}^{.c} &=&{\partial }_{a}L_{.bk}^{c}-\left( {\delta }%
_{k}C_{.ba}^{c}+L_{.dk}^{c}C_{.ba}^{d}-L_{.bk}^{d}C_{.da}^{c}-L_{.ak}^{d}C_{.bd}^{c}\right) +C_{.bd}^{c}P_{.ka}^{d},
\notag \\
S_{j.bc}^{.i} &=&{\partial }_{c}C_{.jb}^{i}-{\partial }%
_{b}C_{.jc}^{i}+C_{.jb}^{h}C_{.hc}^{i}-C_{.jc}^{h}C_{hb}^{i},  \notag \\
S_{b.cd}^{.a} &=&{\partial }_{d}C_{.bc}^{a}-{\partial }%
_{c}C_{.bd}^{a}+C_{.bc}^{e}C_{.ed}^{a}-C_{.bd}^{e}C_{.ec}^{a}.  \notag
\end{eqnarray}

We note that d--torsions (\ref{dtorsc}) and d--curvatures (\ref{dcurvatures}%
) are computed in explicit form by particular cases of d--connections (\ref%
{berwald}), (\ref{inters}) and (\ref{dchrist}).

For cv--bundles we have
\begin{eqnarray}
\check{R}_{h.jk}^{.i} &=&d^{i}\cdot \mathbf{R}\left( \delta _{k},\delta
_{j}\right) \delta _{h},~\check{R}_{\ a.jk}^{b}=\check{\delta}_{a}\cdot
\mathbf{R}\left( \delta _{k},\delta _{j}\right) \check{\partial}^{b},
\label{rpsa} \\
\check{P}_{j.k}^{.i\quad c} &=&d^{i}\cdot \mathbf{R}\left( \check{\partial}%
^{c},\partial _{k}\right) \delta _{j},~\check{P}_{\quad a.k}^{b\quad c}=%
\check{\delta}_{a}\cdot \mathbf{R}\left( \check{\partial}^{c},\partial
_{k}\right) \check{\partial}^{b},  \notag \\
\check{S}_{j.}^{.ibc} &=&d^{i}\cdot \mathbf{R}\left( \check{\partial}^{c},%
\check{\partial}^{b}\right) \delta _{j},~\check{S}_{.a}^{b.cd}=\check{\delta}%
_{a}\cdot \mathbf{R}\left( \check{\partial}^{d},\check{\partial}^{c}\right)
\check{\partial}^{b}.  \notag
\end{eqnarray}%
and
\begin{eqnarray}
\check{R}_{h.jk}^{.i} &=&{\check{\delta}}_{h}L_{.hj}^{i}-{\check{\delta}}%
_{j}L_{.hk}^{i}+L_{.hj}^{m}L_{mk}^{i}-L_{.hk}^{m}L_{mj}^{i}+C_{.h}^{i~a}%
\check{R}_{.ajk},  \label{dcurvaturesa} \\
\check{R}_{.ajk}^{b.} &=&{\check{\delta}}_{k}\check{L}_{a.j}^{~b}-{\check{%
\delta}}_{j}\check{L}_{~b~k}^{a}+\check{L}_{cj}^{~b}\check{L}_{.ak}^{~c}-%
\check{L}_{ck}^{b}\check{L}_{a.j}^{~c}+\check{C}_{a}^{~bc}\check{R}_{c.jk},
\notag \\
\check{P}_{j.k}^{.i~a} &=&{\check{\partial}}^{a}L_{.jk}^{i}-\left( {\check{%
\delta}}_{k}\check{C}_{.j}^{i~a}+L_{.lk}^{i}\check{C}_{.j}^{l~a}-L_{.jk}^{l}%
\check{C}_{.l}^{i~a}-\check{L}_{ck}^{~a}\check{C}_{.j}^{i~c}\right) +%
\check{C}_{.j}^{i~b}\check{P}_{bk}^{~\quad a},  \notag \\
\check{P}_{ck}^{b~a} &=&{\check{\partial}}^{a}\check{L}_{c.k}^{~b}-({\check{%
\delta}}_{k}\check{C}_{c.}^{~ba}+\check{L}_{c.k}^{bd}\check{C}_{d}^{\ ba}-%
\check{L}_{d.k}^{\quad b}\check{C}_{c.}^{\ ad})-\check{L}_{dk}^{\quad a}%
\check{C}_{c.}^{\ bd})+\check{C}_{c.}^{\ bd}\check{P}_{d.k}^{\quad a},
\notag \\
\check{S}_{j.}^{.ibc} &=&{\check{\partial}}^{c}\check{C}_{.j}^{i\ b}-{\check{%
\partial}}^{b}\check{C}_{.j}^{i\ c}+\check{C}_{.j}^{h\ b}\check{C}_{.h}^{i\
c}-\check{C}_{.j}^{h\ c}\check{C}_{h}^{i\ b},  \notag \\
\check{S}_{\ a.}^{b\ cd} &=&{\check{\partial}}^{d}\check{C}_{a.}^{\ bc}-{%
\check{\partial}}^{c}\check{C}_{a.}^{\ bd}+\check{C}_{e.}^{\ bc}\check{C}%
_{a.}^{\ ed}-\check{C}_{e.}^{\ bd}\check{C}_{.a}^{\ ec}.  \notag
\end{eqnarray}

The formulas for curvature can be also generalized for hvc--bundles (on
every shell we must write (\ref{dtorsc}) or (\ref{torsionsa}) in dependence
of the type of shell, vector or co-vector one, we are dealing).

\section{Generalizations of Finsler Geometry}

We outline the basic definitions and formulas for Finsler, Lagrange and
generalized Lagrange spaces (constructed on tangent bundle) and for Cartan,
Hamilton and generalized Hamilton spaces (constructed on cotangent bundle).
The original results are given in details in the monographs \cite{ma87,mhss}%
, see also developments for suberbundles \cite{vstr2,vbook}.

\subsection{Finsler Spaces}

The Finsler geometry is modeled on tangent bundle $TM.$

\begin{definition}
A Finsler space (manifold) is a pair $F^{n}=\left( M,F(x,y)\right) $ \ where
$M$ is a real $n$--dimensional differentiable manifold and $F:TM\rightarrow
\mathcal{R}$ \ is a scalar function which satisfy the following conditions:

\begin{enumerate}
\item $F$ is a differentiable function on the manifold $\widetilde{TM}$ $%
=TM\backslash \{0\}$ and $F$ is continuous on the null section of the
projection $\pi :TM\rightarrow M;$

\item $F$ is a positive function, homogeneous on the fibers of the $TM,$ i.
e. $F(x,\lambda y)=\lambda F(x,y),\lambda \in \mathcal{R};$

\item The Hessian of $F^{2}$ with elements
\begin{equation}
g_{ij}^{(F)}(x,y)=\frac{1}{2}\frac{\partial ^{2}F^{2}}{\partial
y^{i}\partial y^{j}}  \label{finm}
\end{equation}
is positively defined on $\widetilde{TM}.$
\end{enumerate}
\end{definition}

The function $F(x,y)$ and $g_{ij}(x,y)$ are called respectively the
fundamental function and the fundamental (or metric) tensor of the Finsler
space $F.$

One considers ''anisotropic'' (depending on directions $y^{i})$ Christoffel
symbols, for simplicity we write $g_{ij}^{(F)}=g_{ij},$
\begin{equation*}
\gamma _{~jk}^{i}(x,y)=\frac{1}{2}g^{ir}\left( \frac{\partial g_{rk}}{%
\partial x^{j}}+\frac{\partial g_{jr}}{\partial x^{k}}-\frac{\partial g_{jk}%
}{\partial x^{r}}\right) ,
\end{equation*}
which are used for definition of the Cartan N--connection,
\begin{equation}
N_{(c)~j}^{i}=\frac{1}{2}\frac{\partial }{\partial y^{j}}\left[ \gamma
_{~nk}^{i}(x,y)y^{n}y^{k}\right] .  \label{ncartan}
\end{equation}
This N--connection can be used for definition of an almost complex structure
like in (\ref{alcomp}) and to define on $TM$ a d--metric
\begin{equation}
\mathbf{G}_{(F)}=g_{ij}(x,y)dx^{i}\otimes dx^{j}+g_{ij}(x,y)\delta
y^{i}\otimes \delta y^{j},  \label{dmfin}
\end{equation}
with $g_{ij}(x,y)$ taken as (\ref{finm}).

Using the Cartan N--connection (\ref{ncartan}) and Finsler metric tensor (%
\ref{finm}) (or, equivalently, the d--metric (\ref{dmfin})) we can introduce
the canonical d--connection
\begin{equation*}
D\Gamma \left( N_{(c)}\right) =\Gamma _{(c)\beta \gamma }^{\alpha }=\left(
L_{(c)~jk}^{i},C_{(c)~jk}^{i}\right)
\end{equation*}
with the coefficients computed like in (\ref{kahlerconm}) and (\ref{lccoef})
with $h_{ab}\rightarrow g_{ij}.$ The d--connection $D\Gamma \left(
N_{(c)}\right) $ has the unique property that it is torsionless and
satisfies the metricity conditions both for the horizontal and vertical
components, i. e. $D_{\alpha }g_{\beta \gamma }=0.$

The d--curvatures
\begin{equation*}
\check{R}_{h.jk}^{.i}=\{\check{R}_{h.jk}^{.i},\check{P}_{j.k}^{.i\quad
l},S_{(c)j.kl}^{.i}\}
\end{equation*}
on a Finsler space provided with Cartan N--connection and Finsler metric
structures are computed following the formulas (\ref{dcurvatures}) when the $%
a,b,c...$ indices are identified with $i,j,k,...$ indices. It should be
emphasized that in this case all values $g_{ij,}\Gamma _{(c)\beta \gamma
}^{\alpha }$ and $R_{(c)\beta .\gamma \delta }^{.\alpha }$ are defined by a
fundamental function $F\left( x,y\right) .$

In general, we can consider that a Finsler space is provided with a metric $%
g_{ij}=\partial ^{2}F^{2}/2\partial y^{i}\partial y^{j},$ but the
N--connection and d--connection are be defined in a different manner, even
not be determined by $F.$

\subsection{Lagrange and Generalized Lagrange Spaces}

The notion of Finsler spaces was generalized by J. Kern \cite{ker} and R.
Miron \cite{mironlg}. It is widely developed in monographs \cite{ma87} and
extended to superspaces in Refs. \cite{vlasg,vstr2,vbook}.

The idea of extension was to consider instead of the homogeneous fundamental
function $F(x,y)$ in a Finsler space a more general one, a Lagrangian $%
L\left( x,y\right) $, defined as a differentiable mapping $L:(x,y)\in
TM\rightarrow L(x,y)\in \mathcal{R},$ of class $C^{\infty }$ on manifold $%
\widetilde{TM}$ and continuous on the null section $0:M\rightarrow TM$ of
the projection $\pi :TM\rightarrow M.$ A Lagrangian is regular if it is
differentiable and the Hessian
\begin{equation}
g_{ij}^{(L)}(x,y)=\frac{1}{2}\frac{\partial ^{2}L^{2}}{\partial
y^{i}\partial y^{j}}  \label{lagm}
\end{equation}
is of rank $n$ on $M.$

\begin{definition}
A Lagrange space is a pair $L^{n}=\left( M,L(x,y)\right) $ where $M$ is a
smooth real $n$--dimensional manifold provided with regular Lagrangian \ $%
L(x,y)$ structure $L:TM\rightarrow \mathcal{R}$ $\ $for which $g_{ij}(x,y)$
from (\ref{lagm}) has a constant signature over the manifold $\widetilde{TM}%
. $
\end{definition}

The fundamental Lagrange function $L(x,y)$ defines a canonical
N--con\-nec\-ti\-on
\begin{equation*}
N_{(cL)~j}^{i}=\frac{1}{2}\frac{\partial }{\partial y^{j}}\left[
g^{ik}\left( \frac{\partial ^{2}L^{2}}{\partial y^{k}\partial y^{h}}y^{h}-%
\frac{\partial L}{\partial x^{k}}\right) \right]
\end{equation*}
as well a d-metric
\begin{equation}
\mathbf{G}_{(L)}=g_{ij}(x,y)dx^{i}\otimes dx^{j}+g_{ij}(x,y)\delta
y^{i}\otimes \delta y^{j},  \label{dmlag}
\end{equation}
with $g_{ij}(x,y)$ taken as (\ref{lagm}). As well we can introduce an almost
K\"{a}hlerian structure and an almost Hermitian model of $L^{n},$ denoted as
$H^{2n}$ as in the case of Finsler spaces but with a proper fundamental
Lagange function and metric tensor $g_{ij}.$ The canonical metric
d--connection $D\Gamma \left( N_{(cL)}\right) =\Gamma _{(cL)\beta \gamma
}^{\alpha }=\left( L_{(cL)~jk}^{i},C_{(cL)~jk}^{i}\right) $ is to computed
by the same formulas (\ref{kahlerconm}) and (\ref{lccoef}) with $%
h_{ab}\rightarrow g_{ij}^{(L)},$ for $N_{(cL)~j}^{i}.$ The d--torsions (\ref%
{dtorsc}) and d--curvatures (\ref{dcurvatures}) are defined, in this case,
by $L_{(cL)~jk}^{i}$ and $C_{(cL)~jk}^{i}.$ We also note that instead of $%
N_{(cL)~j}^{i}$ and $\Gamma _{(cL)\beta \gamma }^{\alpha }$ one can consider
on a $L^{n}$--space arbitrary N--connections $N_{~j}^{i},$ d--connections $%
\Gamma _{\beta \gamma }^{\alpha }$ which are not defined only by $L(x,y)$
and $g_{ij}^{(L)}$ but can be metric, or non--metric with respect to the
Lagrange metric.

The next step of generalization is to consider an arbitrary metric $%
g_{ij}\left( x,y\right) $ on $TM$ instead of (\ref{lagm}) which is the
second derivative of ''anisotropic'' coordinates $y^{i}$ of a Lagrangian %
\cite{mironlg}.

\begin{definition}
A generalized Lagrange space is a pair $GL^{n}=\left( M,g_{ij}(x,y)\right) $
where $g_{ij}(x,y)$ is a covariant, symmetric d--tensor field, of rank $n$
and of constant signature on $\widetilde{TM}.$
\end{definition}

\bigskip One can consider different classes of N-- and d--connections on $%
TM, $ which are compatible (metric) or non compatible with (\ref{dmlag}) for
arbitrary $g_{ij}(x,y).$ We can apply all formulas for d--connections,
N--curvatures, d--torsions and d--curvatures as in a v--bundle $\mathcal{E},$
but reconsidering them on $TM,$ by changing \ $h_{ab}\rightarrow g_{ij}(x,y)$
and $N_{i}^{a}\rightarrow N_{~i}^{k}.$

\subsection{Cartan Spaces}

The theory of Cartan spaces (see, for instance, \cite{run,kaw1}) \ was
formulated in a new fashion in R. Miron's works \cite{mironc1} by
considering them as duals to the Finsler spaces (see details and references
in \cite{mhss}). Roughly, a Cartan space is constructed on a cotangent
bundle $T^{\ast }M$ like a Finsler space on the corresponding tangent bundle
$TM.$

Consider a real smooth manifold $M,$ the cotangent bundle $\left( T^{\ast
}M,\pi ^{\ast },M\right) $ and the manifold $\widetilde{T^{\ast }M}=T^{\ast
}M\backslash \{0\}.$

\begin{definition}
A Cartan space is a pair $C^{n}=\left( M,K(x,p)\right) $ \ such that $%
K:T^{\ast }M\rightarrow \mathcal{R}$ is a scalar function which satisfy the
following conditions:

\begin{enumerate}
\item $K$ is a differentiable function on the manifold $\widetilde{T^{\ast }M%
}$ $=T^{\ast }M\backslash \{0\}$ and continuous on the null section of the
projection $\pi ^{\ast }:T^{\ast }M\rightarrow M;$

\item $K$ is a positive function, homogeneous on the fibers of the $T^{\ast
}M,$ i. e. $K(x,\lambda p)=\lambda F(x,p),\lambda \in \mathcal{R};$

\item The Hessian of $K^{2}$ with elements
\begin{equation}
\check{g}_{(K)}^{ij}(x,p)=\frac{1}{2}\frac{\partial ^{2}K^{2}}{\partial
p_{i}\partial p_{j}}  \label{carm}
\end{equation}
is positively defined on $\widetilde{T^{\ast }M}.$
\end{enumerate}
\end{definition}

The function $K(x,y)$ and $\check{g}^{ij}(x,p)$ are called \ respectively
the fundamental function and the fundamental (or metric) tensor of the
Cartan space $C^{n}.$ We use symbols like $"\check{g}"$ as to emphasize that
the geometrical objects are defined on a dual space.

One considers ''anisotropic'' (depending on directions, momenta, $p_{i})$
\newline
Christoffel symbols, for simplicity, we write the inverse to (\ref{carm}) as
$g_{ij}^{(K)}=\check{g}_{ij},$
\begin{equation*}
\check{\gamma}_{~jk}^{i}(x,p)=\frac{1}{2}\check{g}^{ir}\left( \frac{\partial
\check{g}_{rk}}{\partial x^{j}}+\frac{\partial \check{g}_{jr}}{\partial x^{k}%
}-\frac{\partial \check{g}_{jk}}{\partial x^{r}}\right) ,
\end{equation*}
which are used for definition of the canonical N--connection,
\begin{equation}
\check{N}_{ij}=\check{\gamma}_{~ij}^{k}p_{k}-\frac{1}{2}\gamma
_{~nl}^{k}p_{k}p^{l}{\breve{\partial}}^{n}\check{g}_{ij},~{\breve{\partial}}%
^{n}=\frac{\partial }{\partial p_{n}}.  \label{nccartan}
\end{equation}
This N--connection can be used for definition of an almost complex structure
like in (\ref{alcomp}) and to define on $T^{\ast }M$ a d--metric
\begin{equation}
\mathbf{\check{G}}_{(k)}=\check{g}_{ij}(x,p)dx^{i}\otimes dx^{j}+\check{g}%
^{ij}(x,p)\delta p_{i}\otimes \delta p_{j},  \label{dmcar}
\end{equation}
with $\check{g}^{ij}(x,p)$ taken as (\ref{carm}).

Using the canonical N--connection (\ref{nccartan}) and Finsler metric tensor
(\ref{carm}) (or, equivalently, the d--metric (\ref{dmcar}) we can introduce
the canonical d--connection
\begin{equation*}
D\check{\Gamma}\left( \check{N}_{(k)}\right) =\check{\Gamma}_{(k)\beta
\gamma }^{\alpha }=\left( \check{H}_{(k)~jk}^{i},\check{C}_{(k)~i}^{\quad
jk}\right)
\end{equation*}%
with the coefficients \
\begin{equation*}
\check{H}_{(k)~jk}^{i}=\frac{1}{2}\check{g}^{ir}\left( \check{\delta}_{j}%
\check{g}_{rk}+\check{\delta}_{k}\check{g}_{jr}-\check{\delta}_{r}\check{g}%
_{jk}\right) ,\check{C}_{(k)~i}^{\quad jk}=\check{g}_{is}{\breve{\partial}}%
^{s}\check{g}^{jk},
\end{equation*}%
The d--connection $D\check{\Gamma}\left( \check{N}_{(k)}\right) $ has the
unique property that it is torsionless and satisfies the metricity
conditions both for the horizontal and vertical components, i. e. $\check{D}%
_{\alpha }\check{g}_{\beta \gamma }=0.$

The d--curvatures
\begin{equation*}
\check{R}_{(k)\beta .\gamma \delta }^{.\alpha
}=\{R_{(k)h.jk}^{.i},P_{(k)j.km}^{.i},\check{S}_{j.}^{.ikl}\}
\end{equation*}
on a Finsler space provided with Cartan N--connection and Finsler metric
structures are computed following the formulas (\ref{dcurvaturesa}) when the
$a,b,c...$ indices are identified with $i,j,k,...$ indices. It should be
emphasized that in this case all values $\check{g}_{ij,}\check{\Gamma}%
_{(k)\beta \gamma }^{\alpha }$ and $\check{R}_{(k)\beta .\gamma \delta
}^{.\alpha }$ are defined by a fundamental function $K\left( x,p\right) .$

In general, we can consider that a Cartan space is provided with a metric $%
\check{g}^{ij}=\partial ^{2}K^{2}/2\partial p_{i}\partial p_{j},$ but the
N--connection and d--connection could be defined in a different manner, even
not be determined by $K.$

\subsection{ Generalized Hamilton and Hamilton Spaces}

The geometry of Hamilton spaces was defined and investigated by R. Miron in
Refs. \cite{mironh1} (see details and references in \cite{mhss}). It was
developed on the cotangent bundel as a dual geometry to the geometry of
Lagrange spaces. \ Here we start with the definition of generalized Hamilton
spaces and then consider the particular case.

\begin{definition}
A generalized Hamilton space is a pair\newline
$GH^{n}=\left( M,\check{g}^{ij}(x,p)\right) $ where $M$ is a real $n$%
--dimensional manifold and $\check{g}^{ij}(x,p)$ is a contravariant,
symmetric, nondegenerate of rank $n$ and of constant signature on $%
\widetilde{T^{\ast }M}.$
\end{definition}

\bigskip The value $\check{g}^{ij}(x,p)$ is called the fundamental (or
metric) tensor of the space $GH^{n}.$ One can define such values for every
paracompact manifold $M.$ In general, a N--connection on $GH^{n}$ is not
determined by $\check{g}^{ij}.$ Therefore we can consider arbitrary
coefficients $\check{N}_{ij}\left( x,p\right) $ and define on $T^{\ast }M$ a
d--metric like (\ref{dmetricvc})
\begin{equation}
{\breve{\mathbf{G}}}={\breve{g}}_{\alpha \beta }\left( {\breve{u}}\right) {%
\breve{\delta}}^{\alpha }\otimes {\breve{\delta}}^{\beta }={\breve{g}}%
_{ij}\left( {\breve{u}}\right) d^{i}\otimes d^{j}+{\check{g}}^{ij}\left( {%
\breve{u}}\right) {\breve{\delta}}_{i}\otimes {\breve{\delta}}_{j},
\label{dmghs}
\end{equation}
This N--coefficients $\check{N}_{ij}\left( x,p\right) $ and d--metric
structure (\ref{dmghs}) allow to define an almost K\"{a}hler model of
generalized Hamilton spaces and to define canonical d--connections,
d--torsions and d-curvatures (see respectively the formulas (\ref{lccoef}), (%
\ref{inters}), (\ref{dtorsca}) and (\ref{dcurvatures}) with the fiber
coefficients redefined for the cotangent bundle $T^{\ast }M$ ).

A generalized Hamilton space $GH^{n}=\left( M,\check{g}^{ij}(x,p)\right) $
is called reducible to a Hamilton one if there exists a Hamilton function $%
H\left( x,p\right) $ on $T^{\ast }M$ such that
\begin{equation}
\check{g}^{ij}(x,p)=\frac{1}{2}\frac{\partial ^{2}H}{\partial p_{i}\partial
p_{j}}.  \label{hsm}
\end{equation}

\begin{definition}
A Hamilton space is a pair $H^{n}=\left( M,H(x,p)\right) $ \ such that $%
H:T^{\ast }M\rightarrow \mathcal{R}$ is a scalar function which satisfy the
following conditions:

\begin{enumerate}
\item $H$ is a differentiable function on the manifold $\widetilde{T^{\ast }M%
}$ $=T^{\ast }M\backslash \{0\}$ and continuous on the null section of the
projection $\pi ^{\ast }:T^{\ast }M\rightarrow M;$

\item The Hessian of $H$ with elements (\ref{hsm}) is positively defined on $%
\widetilde{T^{\ast }M}$ and $\check{g}^{ij}(x,p)$ is nondegenerate matrix of
rank $n$ and of constant signature.
\end{enumerate}
\end{definition}

For Hamilton spaces the canonical N--connection (defined by $H$ and its
Hessian) exists,
\begin{equation*}
\check{N}_{ij}=\frac{1}{4}\{\check{g}_{ij},H\}-\frac{1}{2}\left( \check{g}%
_{ik}\frac{\partial ^{2}H}{\partial p_{k}\partial x^{j}}+\check{g}_{jk}\frac{%
\partial ^{2}H}{\partial p_{k}\partial x^{i}}\right) ,
\end{equation*}
where the Poisson brackets, for arbitrary functions $f$ and $g$ on $T^{\ast
}M,$ act as
\begin{equation*}
\{f,g\}=\frac{\partial f}{\partial p_{i}}\frac{\partial g}{\partial x^{i}}-%
\frac{\partial g}{\partial p_{i}}\frac{\partial p}{\partial x^{i}}.
\end{equation*}

The canonical d--connection $D\check{\Gamma}\left( \check{N}_{(c)}\right) =%
\check{\Gamma}_{(c)\beta \gamma }^{\alpha }=\left( \check{H}_{(c)~jk}^{i},%
\check{C}_{(c)~i}^{\quad jk}\right) $is defined by the coefficients
\begin{equation*}
\check{H}_{(c)~jk}^{i}=\frac{1}{2}\check{g}^{is}\left( \check{\delta}_{j}%
\check{g}_{sk}+\check{\delta}_{k}\check{g}_{js}-\check{\delta}_{s}\check{g}%
_{jk}\right) ,\check{C}_{(c)~i}^{\quad jk}=-\frac{1}{2}\check{g}_{is}\check{%
\partial}^{j}\check{g}^{sk}.
\end{equation*}%
In result we can compute the d--torsions and d--curvatures like on
cv--bundle \ or on Cartan spaces. On Hamilton spaces all such objects are
defined by the Hamilton function $H(x,p)$ and indices have to be
reconsidered for co--fibers of the co-tangent bundle.

\section{Clifford Bundles and N--Connections}

The theory of anisotropic spinors was extended on higher order anisotropic
(ha) spaces \cite{vhsp,vbook,vp}. In brief, such spinors will be called
ha--spinors which are defined as some Clifford ha--structures defined with
respect to a distinguished quadratic form (\ref{dmetrichcv}) on a
hvc--bundle. For simplicity, the bulk of formulas will be given with respect
to higher order vector bundles. To rewrite such formulas for hvc--bundles is
to consider for the ''dual'' shells of higher order anisotropy some dual
vector spaces and associated dual spinors.

\subsection{Distinguished Clifford Algebras}

The typical fiber of dv--bundle $\xi _{d}\ ,\ \pi _{d}:\ HE\oplus
V_{1}E\oplus ...\oplus V_{z}E\rightarrow E$ is a d-vector space, $\mathcal{F}%
=h\mathcal{F}\oplus v_{1}\mathcal{F\oplus }...\oplus v_{z}\mathcal{F},$
split into horizontal $h\mathcal{F}$ and verticals $v_{p}\mathcal{F}%
,p=1,...,z$ subspaces, with a bilinear quadratic form $G(g,h)$ induced by a
hvc--bundle metric (\ref{dmetrichcv}). Clifford algebras (see, for example,
Refs. \cite{kar,pen}) formulated for d-vector spaces will be called Clifford
d-algebras \cite{vjmp,vdeb}. We shall consider the main properties of
Clifford d--algebras. The proof of theorems will be based on the technique
developed in Ref. \cite{kar,vbook,vp} correspondingly adapted to the
distinguished character of spaces in consideration.

Let $k$ be a number field (for our purposes $k=\mathcal{R}$ or $k=\mathcal{C}%
,\mathcal{R}$ and $\mathcal{C},$ are, respectively real and
complex number fields) and define $\mathcal{F},$ as a d-vector
space on $k$ provided with nondegenerate symmetric quadratic form
(metric)\ $G.$ Let $C$ be an algebra on $k$ (not necessarily
commutative) and $j\ :\ \mathcal{F}$ $\rightarrow C$ a
homomorphism of underlying vector spaces such that
$j(u)^2=\;G(u)\cdot 1\ (1$ is the unity in algebra $C$ and
d-vector $u\in \mathcal{F}).$ We are interested in definition of
the pair $\left( C,j\right) $ satisfying the next universitality
conditions. For every $k$-algebra $A$ and arbitrary homomorphism
$\varphi :\mathcal{F}\rightarrow A$ of the underlying d-vector
spaces, such that $\left( \varphi (u)\right) ^2\rightarrow G\left(
u\right) \cdot 1,$ there is a unique homomorphism of algebras
$\psi \ :\ C\rightarrow A$ transforming the chain of maps into a
commutative diagram.

The algebra solving this problem will be denoted as $C\left( \mathcal{F}%
,A\right) $ [equivalently as $C\left( G\right) $ or $C\left( \mathcal{F}%
\right) ]$ and called as Clifford d--algebra associated with pair $\left(
\mathcal{F},G\right).$

\begin{theorem}
The above-presented chain of maps has a unique solution $\left(
C,j\right) $ up to isomorphism.
\end{theorem}

\textbf{Proof:} See Refs. \cite{vsp1,vbook}.

Now we re--formulate for d--algebras the Chevalley theorem \cite{chev}:

\begin{theorem}
The Clifford d-algebra
\begin{equation*}
C\left( h\mathcal{F}\oplus v_1\mathcal{F}\oplus ...\oplus v_z\mathcal{F}%
,g+h_1+...+h_z\right)
\end{equation*}
is naturally isomorphic to $C(g)\otimes C\left( h_1\right) \otimes
...\otimes C\left( h_z\right) .$
\end{theorem}

\textbf{Proof. }\ See Refs. \cite{vsp1,vbook}.

From the presented Theorems, we conclude that all operations with Clifford
d-algebras can be reduced to calculations for $C\left( h\mathcal{F},g\right)
$ and\newline
$C\left( v_{(p)}\mathcal{F},h_{(p)}\right) $ which are usual Clifford
algebras of dimension $2^{n}$ and, respectively, $2^{m_{p}}$ \cite{kar,ati}.

Of special interest is the case when $k=\mathcal{R}$ and $\mathcal{F}$ is
isomorphic to vector space $\mathcal{R}^{p+q,a+b}$ provided with quadratic
form
\begin{equation*}
-x_1^2-...-x_p^2+x_{p+q}^2-y_1^2-...-y_a^2+...+y_{a+b}^2.
\end{equation*}
In this case, the Clifford algebra, denoted as $\left(
C^{p,q},C^{a,b}\right) ,\,$ is generated by the symbols \newline
$e_1^{(x)},e_2^{(x)},...,e_{p+q}^{(x)},e_1^{(y)},e_2^{(y)},...,e_{a+b}^{(y)}$
satisfying properties
\begin{eqnarray*}
\left( e_i\right) ^2 &=&-1~\left( 1\leq i\leq p\right) ,\left( e_j\right)
^2=-1~\left( 1\leq j\leq a\right) , \\
\left( e_k\right) ^2 &=&1~(p+1\leq k\leq p+q), \\
\left( e_j\right) ^2 &=&1~(n+1\leq s\leq a+b),~e_ie_j=-e_je_i,~i\neq j.\,
\end{eqnarray*}

Explicit calculations of $C^{p,q}$ and $C^{a,b}$ are possible by the using
isomorphisms \cite{kar,pen}
\begin{eqnarray*}
C^{p+n,q+n} &\simeq &C^{p,q}\otimes M_{2}\left( \mathcal{R}\right) \otimes
...\otimes M_{2}\left( \mathcal{R}\right) \\
&\cong &C^{p,q}\otimes M_{2^{n}}\left( \mathcal{R}\right) \cong
M_{2^{n}}\left( C^{p,q}\right) ,
\end{eqnarray*}%
where $M_{s}\left( A\right) $ denotes the ring of quadratic matrices of
order $s$ with coefficients in ring $A.$ Here we write the simplest
isomorphisms $C^{1,0}\simeq \mathcal{C},C^{0,1}\simeq \mathcal{R}\oplus
\mathcal{R}$ and $C^{2,0}=\mathcal{H},$ where by $\mathcal{H}$ is denoted
the body of quaternions.

Now, we emphasize that higher order Lagrange and Finsler spaces, denoted $%
H^{2n}$--spaces, admit locally a structure of Clifford algebra on complex
vector spaces. Really, by using almost \ Hermitian structure $J_\alpha
^{\quad \beta }$ and considering complex space $\mathcal{C}^n$ with
nondegenarate quadratic form $\sum_{a=1}^n\left| z_a\right| ^2,~z_a\in
\mathcal{C}^2$ induced locally by metric (\ref{dmetrichcv}) (rewritten in
complex coordinates $z_a=x_a+iy_a)$ we define Clifford algebra $%
\overleftarrow{C}^n=\underbrace{\overleftarrow{C}^1\otimes ...\otimes
\overleftarrow{C}^1}_n,$ where $\overleftarrow{C}^1=\mathcal{C\otimes }_R%
\mathcal{C=C\oplus C}$ or in consequence, $\overleftarrow{C}^n\simeq
C^{n,0}\otimes _{\mathcal{R}}\mathcal{C}\approx C^{0,n}\otimes _{\mathcal{R}}%
\mathcal{C}.$ Explicit calculations lead to isomorphisms
\begin{equation*}
\overleftarrow{C}^2=C^{0,2}\otimes _{\mathcal{R}}\mathcal{C}\approx
M_2\left( \mathcal{R}\right) \otimes _{\mathcal{R}}\mathcal{C}\approx
M_2\left( \overleftarrow{C}^n\right) ,~C^{2p}\approx M_{2^p}\left( \mathcal{C%
}\right)
\end{equation*}
and
\begin{equation*}
\overleftarrow{C}^{2p+1}\approx M_{2^p}\left( \mathcal{C}\right) \oplus
M_{2^p}\left( \mathcal{C}\right) ,
\end{equation*}
which show that complex Clifford algebras, defined locally for $H^{2n}$%
-spaces, have periodicity 2 on $p.$

Considerations presented in the proof of theorem 2.2 show that map $j:%
\mathcal{F}\rightarrow C\left( \mathcal{F}\right) $ is monomorphic, so we
can identify space $\mathcal{F}$ with its image in $C\left( \mathcal{F}%
,G\right) ,$ denoted as $u\rightarrow \overline{u},$ if $u\in C^{(0)}\left(
\mathcal{F},G\right) ~\left( u\in C^{(1)}\left( \mathcal{F},G\right) \right)
;$ then $u=\overline{u}$ ( respectively, $\overline{u}=-u).$

\begin{definition}
The set of elements $u\in C\left( G\right) ^{\ast },$ where $C\left(
G\right) ^{\ast }$ denotes the multiplicative group of invertible elements
of $C\left( \mathcal{F},G\right) $ satisfying $\overline{u}\mathcal{F}%
u^{-1}\in \mathcal{F},$ is called the twisted Clifford d-group, denoted as $%
\widetilde{\Gamma }\left( \mathcal{F}\right).$
\end{definition}

Let $\widetilde{\rho }:\widetilde{\Gamma }\left( \mathcal{F}\right)
\rightarrow GL\left( \mathcal{F}\right) $ be the homorphism given by $%
u\rightarrow \rho \widetilde{u},$ where $\widetilde{\rho }_u\left( w\right) =%
\overline{u}wu^{-1}.$ We can verify that $\ker \widetilde{\rho }=\mathcal{R}%
^{*}$is a subgroup in $\widetilde{\Gamma }\left( \mathcal{F}\right) .$

The canonical map $j:\mathcal{F}\rightarrow C\left( \mathcal{F}\right) $ can
be interpreted as the linear map $\mathcal{F}\rightarrow C\left( \mathcal{F}%
\right) ^0$ satisfying the universal property of Clifford d-algebras. This
leads to a homomorphism of algebras, $C\left( \mathcal{F}\right) \rightarrow
C\left( \mathcal{F}\right) ^t,$ considered by an anti-involution of $C\left(
\mathcal{F}\right) $ and denoted as $u\rightarrow ~^tu.$ More exactly, if $%
u_1...u_n\in \mathcal{F,}$ then $t_u=u_n...u_1$ and $^t\overline{u}=%
\overline{^tu}=\left( -1\right) ^nu_n...u_1.$

\begin{definition}
The spinor norm of arbitrary $u\in C\left( \mathcal{F}\right) $ is defined as%
\newline
$S\left( u\right) =~^{t}\overline{u}\cdot u\in C\left( \mathcal{F}\right) .$
\end{definition}

It is obvious that if $u,u^{\prime },u^{\prime \prime }\in \widetilde{\Gamma
}\left( \mathcal{F}\right) ,$ then $S(u,u^{\prime })=S\left( u\right)
S\left( u^{\prime }\right) $ and \newline
$S\left( uu^{\prime }u^{\prime \prime }\right) =S\left( u\right) S\left(
u^{\prime }\right) S\left( u^{\prime \prime }\right) .$ For $u,u^{\prime
}\in \mathcal{F} S\left( u\right) =-G\left( u\right) $ and \newline
$S\left( u,u^{\prime }\right) =S\left( u\right) S\left( u^{\prime }\right)
=S\left( uu^{\prime }\right) .$

Let us introduce the orthogonal group $O\left( G\right) \subset GL\left(
G\right) $ defined by metric $G$ on $\mathcal{F}$ and denote sets
\begin{equation*}
SO\left( G\right) =\{u\in O\left( G\right) ,\det \left| u\right|
=1\},~Pin\left( G\right) =\{u\in \widetilde{\Gamma }\left( \mathcal{F}%
\right) ,S\left( u\right) =1\}
\end{equation*}
and $Spin\left( G\right) =Pin\left( G\right) \cap C^0\left( \mathcal{F}%
\right) .$ For ${\mathcal{F}\cong \mathcal{R}}^{n+m}$ we write $Spin\left(
n_E\right) .$ By straightforward calculations (see similar considerations in
Ref. \cite{kar}) we can verify the exactness of these sequences:
\begin{eqnarray*}
1 &\rightarrow &\mathcal{Z}/2\rightarrow Pin\left( G\right) \rightarrow
O\left( G\right) \rightarrow 1, \\
1 &\rightarrow &\mathcal{Z}/2\rightarrow Spin\left( G\right) \rightarrow
SO\left( G\right) \rightarrow 0, \\
1 &\rightarrow &\mathcal{Z}/2\rightarrow Spin\left( n_E\right) \rightarrow
SO\left( n_E\right) \rightarrow 1.
\end{eqnarray*}
We conclude this subsection by emphasizing that the spinor norm was defined
with respect to a quadratic form induced by a metric in dv--bundle $\mathcal{%
E}^{<z>}$. This approach differs from those presented in Refs. \cite{asa88}
and \cite{ono}.

\subsection{Clifford Ha--Bundles}

We shall consider two variants of generalization of spinor constructions
defined for d-vector spaces to the case of distinguished vector bundle
spaces enabled with the structure of N-connection. The first is to use the
extension to the category of vector bundles. The second is to define the
Clifford fibration associated with compatible linear d-connection and metric
$G$ on a dv--bundle. We shall analyze both variants.

\subsubsection{Clifford d--module structure in dv--bundles}

Because functor $\mathcal{F}\rightarrow C(\mathcal{F})$ is smooth we can
extend it to the category of vector bundles of type
\begin{equation*}
\xi ^{<z>}=\{\pi _{d}:HE^{<z>}\oplus V_{1}E^{<z>}\oplus ...\oplus
V_{z}E^{<z>}\rightarrow E^{<z>}\}.
\end{equation*}%
Recall that by $\mathcal{F}$ we denote the typical fiber of such bundles.
For $\xi ^{<z>}$ we obtain a bundle of algebras, denoted as $C\left( \xi
^{<z>}\right) ,\,$ such that $C\left( \xi ^{<z>}\right) _{u}=C\left(
\mathcal{F}_{u}\right) .$ Multiplication in every fiber defines a continuous
map $C\left( \xi ^{<z>}\right) \times C\left( \xi ^{<z>}\right) \rightarrow
C\left( \xi ^{<z>}\right).$ If $\xi ^{<z>}$ is a distinguished vector bundle
on number field $k$, $C\left( \xi ^{<z>}\right) $-module, the d-module, on $%
\xi ^{<z>}$ is given by the continuous map $C\left( \xi ^{<z>}\right) \times
_{E}\xi ^{<z>}\rightarrow \xi ^{<z>}$ with every fiber $\mathcal{F}_{u}$
provided with the structure of the $C\left( \mathcal{F}_{u}\right)$--module,
correlated with its $k$-module structure, Because $\mathcal{F}\subset
C\left( \mathcal{F}\right) ,$ we have a fiber to fiber map $\mathcal{F}%
\times _{E}\xi ^{<z>}\rightarrow \xi ^{<z>},$ inducing on every fiber the
map $\mathcal{F}_{u}\times _{E}\xi _{(u)}^{<z>}\rightarrow \xi _{(u)}^{<z>}$
($\mathcal{R}$-linear on the first factor and $k$-linear on the second one
). Inversely, every such bilinear map defines on $\xi ^{<z>}$ the structure
of the $C\left( \xi ^{<z>}\right) $-module by virtue of universal properties
of Clifford d--algebras. Equivalently, the above--mentioned bilinear map
defines a morphism of v--bundles
\begin{equation*}
m:\xi ^{<z>}\rightarrow HOM\left( \xi ^{<z>},\xi ^{<z>}\right) \quad \lbrack
HOM\left( \xi ^{<z>},\xi ^{<z>}\right)
\end{equation*}%
denotes the bundles of homomorphisms] when $\left( m\left( u\right) \right)
^{2}=G\left( u\right) $ on every point.

Vector bundles $\xi ^{<z>}$ provided with $C\left( \xi ^{<z>}\right) $%
--structures are objects of the category with morphisms being morphisms of
dv-bundles, which induce on every point $u\in \xi ^{<z>}$ morphisms of $%
C\left( \mathcal{F}_u\right) -$modules. This is a Banach category contained
in the category of finite-dimensional d-vector spaces on filed $k.$

Let us denote by $H^s\left( \mathcal{E}^{<z>},GL_{n_E}\left( \mathcal{R}%
\right) \right) ,\,$ where $n_E=n+m_1+...+m_z,\,$ the s-dimensional
cohomology group of the algebraic sheaf of germs of continuous maps of
dv-bundle $\mathcal{E}^{<z>}$ with group $GL_{n_E}\left( \mathcal{R}\right) $
the group of automorphisms of $\mathcal{R}^{n_E}\,$ (for the language of
algebraic topology see, for example, Refs. \cite{kar}. We shall also use the
group $SL_{n_E}\left( \mathcal{R}\right) =\{A\subset GL_{n_E}\left( \mathcal{%
R}\right) ,\det A=1\}.\,$ Here we point out that cohomologies $H^s(M,Gr)$
characterize the class of a principal bundle $\pi :P\rightarrow M $ on $M$
with structural group $Gr.$ Taking into account that we deal with bundles
distinguished by an N-connection we introduce into consideration
cohomologies $H^s\left( \mathcal{E}^{<z>},GL_{n_E}\left( \mathcal{R}\right)
\right) $ as distinguished classes (d-classes) of bundles $\mathcal{E}^{<z>}$
provided with a global N-connection structure.

For a real vector bundle $\xi ^{<z>}$ on compact base $\mathcal{E}^{<z>}$ we
can define the orientation on $\xi ^{<z>}$ as an element $\alpha _d\in
H^1\left( \mathcal{E}^{<z>},GL_{n_E}\left( \mathcal{R}\right) \right) $
whose image on map
\begin{equation*}
H^1\left( \mathcal{E}^{<z>},SL_{n_E}\left( \mathcal{R}\right) \right)
\rightarrow H^1\left( \mathcal{E}^{<z>},GL_{n_E}\left( \mathcal{R}\right)
\right)
\end{equation*}
is the d-class of bundle $\mathcal{E}^{<z>}.$

\begin{definition}
The spinor structure on $\xi ^{<z>}$ is defined as an element\newline
$\beta _d\in H^1\left( \mathcal{E}^{<z>},Spin\left( n_E\right) \right) $
whose image in the composition
\begin{equation*}
H^1\left( \mathcal{E}^{<z>},Spin\left( n_E\right) \right) \rightarrow
H^1\left( \mathcal{E}^{<z>},SO\left( n_E\right) \right) \rightarrow
H^1\left( \mathcal{E}^{<z>},GL_{n_E}\left( \mathcal{R}\right) \right)
\end{equation*}
is the d-class of $\mathcal{E}^{<z>}.$
\end{definition}

The above definition of spinor structures can be re--formulated in terms of
principal bundles. Let $\xi ^{<z>}$ be a real vector bundle of rank n+m on a
compact base $\mathcal{E}^{<z>}.$ If there is a principal bundle $P_d$ with
structural group $SO( n_E ) $ or $Spin( n_E ) ],$ this bundle $\xi ^{<z>}$
can be provided with orientation (or spinor) structure. The bundle $P_d$ is
associated with element\newline
$\alpha _d\in H^1\left(\mathcal{E}^{<z>},SO(n_{<z>})\right) $ [or $\beta
_d\in H^1\left( \mathcal{E}^{<z>}, Spin\left( n_E\right) \right) .$

We remark that a real bundle is oriented if and only if its first
Stiefel--Whitney d--class vanishes,
\begin{equation*}
w_1\left( \xi _d\right) \in H^1\left( \xi ,\mathcal{Z}/2\right) =0,
\end{equation*}
where $H^1\left( \mathcal{E}^{<z>},\mathcal{Z}/2\right) $ is the first group
of Chech cohomology with coefficients in $\mathcal{Z}/2,$ Considering the
second Stiefel--Whitney class $w_2\left( \xi ^{<z>}\right) \in H^2\left(
\mathcal{E}^{<z>},\mathcal{Z}/2\right) $ it is well known that vector bundle
$\xi ^{<z>}$ admits the spinor structure if and only if $w_2\left( \xi
^{<z>}\right) =0.$ Finally, we emphasize that taking into account that base
space $\mathcal{E}^{<z>}$ is also a v-bundle, $p:E^{<z>}\rightarrow M,$ we
have to make explicit calculations in order to express cohomologies $%
H^s\left( \mathcal{E}^{<z>},GL_{n+m}\right) \,$ and $H^s\left( \mathcal{E}%
^{<z>},SO\left( n+m\right) \right) $ through cohomologies
\begin{equation*}
H^s\left( M,GL_n\right) ,H^s\left( M,SO\left( m_1\right) \right)
,...H^s\left( M,SO\left( m_z\right) \right) ,
\end{equation*}
which depends on global topological structures of spaces $M$ and $\mathcal{E}%
^{<z>}$ $.$ For general bundle and base spaces this requires a cumbersome
cohomological calculus.

\subsubsection{Clifford fibration}

Another way of defining the spinor structure is to use Clifford fibrations.
Consider the principal bundle with the structural group $Gr$ being a
subgroup of orthogonal group $O\left( G\right) ,$ where $G$ is a quadratic
nondegenerate form ) defined on the base (also being a bundle space) space $%
\mathcal{E}^{<z>}.$ The fibration associated to principal fibration $P\left(
\mathcal{E}^{<z>},Gr\right) $ with a typical fiber having Clifford algebra $%
C\left( G\right) $ is, by definition, the Clifford fibration $PC\left(
\mathcal{E}^{<z>},Gr\right) .$ We can always define a metric on the Clifford
fibration if every fiber is isometric to $PC\left( \mathcal{E}%
^{<z>},G\right) $ (this result is proved for arbitrary quadratic forms $G$
on pseudo--Riemannian bases. If, additionally, $Gr\subset SO\left( G\right) $
a global section can be defined on $PC\left( G\right) .$

Let $\mathcal{P}\left( \mathcal{E}^{<z>},Gr\right) $ be the set of principal
bundles with differentiable base $\mathcal{E}^{<z>}$ and structural group $%
Gr.$ If $g:Gr\rightarrow Gr^{\prime }$ is an homomorphism of Lie groups and $%
P( \mathcal{E}^{<z>},Gr)$ $\subset \mathcal{P}\left( \mathcal{E}%
^{<z>},Gr\right) $ (for simplicity in this subsection we shall denote
mentioned bundles and sets of bundles as $P,P^{\prime }$ and respectively, $%
\mathcal{P},\mathcal{P}^{\prime }),$ we can always construct a principal
bundle with the property that there is an homomorphism $f:P^{\prime
}\rightarrow P$ of principal bundles which can be projected to the identity
map of $\mathcal{E}^{<z>}$ and corresponds to isomorphism $g:Gr\rightarrow
Gr^{\prime }.$ If the inverse statement also holds, the bundle $P^{\prime }$
is called as the extension of $P$ associated to $g$ and $f$ is called the
extension homomorphism denoted as $\widetilde{g.}$

Now we can define distinguished spinor structures on bundle spaces.

\begin{definition}
Let $P\in \mathcal{P}\left( \mathcal{E}^{<z>},O\left( G\right) \right) $ be
a principal bundle. A distinguished spinor structure of $P,$ equivalently a
ds-structure of $\mathcal{E}^{<z>}$ is an extension $\widetilde{P}$ of $P$
associated to homomorphism $h:PinG\rightarrow O\left( G\right) $ where $%
O\left( G\right) $ is the group of orthogonal rotations, generated by metric
$G,$ in bundle $\mathcal{E}^{<z>}.$
\end{definition}

So, if $\widetilde{P}$ is a spinor structure of the space $\mathcal{E}%
^{<z>}, $ then $\widetilde{P}\in \mathcal{P}\left( \mathcal{E}%
^{<z>},PinG\right).$

The definition of spinor structures on varieties was given in Ref. \cite{cru}
one had been proved that a necessary and sufficient condition for a space
time to be orientable is to admit a global field of orthonormalized frames.
We mention that spinor structures can be also defined on varieties modeled
on Banach spaces \cite{ana77}. As we have shown similar constructions are
possible for the cases when space time has the structure of a v-bundle with
an N-connection.

\begin{definition}
A special distinguished spinor structure, ds-structure, of principal bundle $%
P=P\left( \mathcal{E}^{<z>},SO\left( G\right) \right) $ is a principal bundle%
$\widetilde{P}=\widetilde{P}\left( \mathcal{E}^{<z>},SpinG\right) $ for
which a homomorphism of principal bundles $\widetilde{p}:\widetilde{P}%
\rightarrow P,$ projected on the identity map of $\mathcal{E}^{<z>}$ and
corresponding to representation
\begin{equation*}
R:SpinG\rightarrow SO\left( G\right) ,
\end{equation*}%
is defined.
\end{definition}

In the case when the base space variety is oriented, there is a natural
bijection between tangent spinor structures with a common base. For special
ds--structures we can define, as for any spinor structure, the concepts of
spin tensors, spinor connections, and spinor covariant derivations (see
Refs. \cite{vdeb,vsp1}).

\subsection{Almost Complex Spinor Structures}

Almost complex structures are an important characteristic of $H^{2n}$-spaces
and of osculator bundles $Osc^{k=2k_{1}}(M),$ where $k_{1}=1,2,...$ . For
simplicity in this subsection we restrict our analysis to the case of $%
H^{2n} $-spaces. We can rewrite the almost Hermitian metric \cite{ma87}, $%
H^{2n}$-metric in complex form \cite{vjmp}:

\begin{equation}
G=H_{ab}\left( z,\xi \right) dz^a\otimes dz^b,  \label{2.38a}
\end{equation}
where
\begin{equation*}
z^a=x^a+iy^a,~\overline{z^a}=x^a-iy^a,~H_{ab}\left( z,\overline{z}\right)
=g_{ab}\left( x,y\right) \mid _{y=y\left( z,\overline{z}\right) }^{x=x\left(
z,\overline{z}\right) },
\end{equation*}
and define almost complex spinor structures. For given metric (\ref{2.38a})
on $H^{2n}$-space there is always a principal bundle $P^U$ with unitary
structural group U(n) which allows us to transform $H^{2n}$-space into
v-bundle $\xi ^U\approx P^U\times _{U\left( n\right) }\mathcal{R}^{2n}.$
This statement will be proved after we introduce complex spinor structures
on oriented real vector bundles \cite{kar}.

Let us consider momentarily $k=\mathcal{C}$ and introduce into consideration
[instead of the group $Spin(n)]$ the group $Spin^c\times _{\mathcal{Z}%
/2}U\left( 1\right) $ being the factor group of the product $Spin(n)\times
U\left( 1\right) $ with the respect to equivalence
\begin{equation*}
\left( y,z\right) \sim \left( -y,-a\right) ,\quad y\in Spin(m).
\end{equation*}
This way we define the short exact sequence
\begin{equation}
1\rightarrow U\left( 1\right) \rightarrow Spin^c\left( n\right) \overset{S^c}%
{\to }SO\left( n\right) \rightarrow 1,  \label{2.39a}
\end{equation}
where $\rho ^c\left( y,a\right) =\rho ^c\left( y\right) .$ If $\lambda $ is
oriented , real, and rank $n,$ $\gamma $-bundle $\pi :E_\lambda \rightarrow
M^n,$ with base $M^n,$ the complex spinor structure, spin structure, on $%
\lambda $ is given by the principal bundle $P$ with structural group $%
Spin^c\left( m\right) $ and isomorphism $\lambda \approx P\times
_{Spin^c\left( n\right) }\mathcal{R}^n$ (see (\ref{2.39a})). For such
bundles the categorial equivalence can be defined as
\begin{equation}
\epsilon ^c:\mathcal{E}_{\mathcal{C}}^T\left( M^n\right) \rightarrow
\mathcal{E}_{\mathcal{C}}^\lambda \left( M^n\right) ,  \label{2.40a}
\end{equation}
where $\epsilon ^c\left( E^c\right) =P\bigtriangleup _{Spin^c\left( n\right)
}E^c$ is the category of trivial complex bundles on $M^n,\mathcal{E}_{%
\mathcal{C}}^\lambda \left( M^n\right) $ is the category of complex
v-bundles on $M^n$ with action of Clifford bundle\newline
$C\left( \lambda \right) ,P\bigtriangleup _{Spin^c(n)}$ and $E^c$ is the
factor space of the bundle product $P\times _ME^c$ with respect to the
equivalence $\left( p,e\right) \sim \left( p\widehat{g}^{-1},\widehat{g}%
e\right) ,p\in P,e\in E^c,$ where $\widehat{g}\in Spin^c\left( n\right) $
acts on $E$ by via the imbedding $Spin\left( n\right) \subset C^{0,n}$ and
the natural action $U\left( 1\right) \subset \mathcal{C}$ on complex
v-bundle $\xi ^c,E^c=tot\xi ^c,$ for bundle $\pi ^c:E^c\rightarrow M^n.$

Now we return to the bundle $\xi =\mathcal{E}^{<1>}.$ A real v-bundle (not
being a spinor bundle) admits a complex spinor structure if and only if
there exist a homomorphism $\sigma :U\left( n\right) \rightarrow
Spin^c\left( 2n\right) $ defining a commutative diagram. The explicit
construction of $\sigma $ for arbitrary $\gamma $-bundle is given in Refs. %
\cite{kar} and \cite{ati}. Let $\lambda $ be a complex, rank\thinspace $n,$
spinor bundle with
\begin{equation}
\tau :Spin^c\left( n\right) \times _{\mathcal{Z}/2}U\left( 1\right)
\rightarrow U\left( 1\right)  \label{2.41a}
\end{equation}
the homomorphism defined by formula $\tau \left( \lambda ,\delta \right)
=\delta ^2.$ For $P_s$ being the principal bundle with fiber $Spin^c\left(
n\right) $ we introduce the complex linear bundle $L\left( \lambda ^c\right)
=P_S\times _{Spin^c(n)}\mathcal{C}$ defined as the factor space of $%
P_S\times \mathcal{C}$ on equivalence relation

\begin{equation*}
\left( pt,z\right) \sim \left( p,l\left( t\right) ^{-1}z\right) ,
\end{equation*}
where $t\in Spin^c\left( n\right) .$ This linear bundle is associated to
complex spinor structure on $\lambda ^c.$

If $\lambda ^c$ and $\lambda ^{c^{\prime }}$ are complex spinor bundles, the
Whitney sum $\lambda ^c\oplus \lambda ^{c^{\prime }}$ is naturally provided
with the structure of the complex spinor bundle. This follows from the
holomorphism
\begin{equation}
\omega ^{\prime }:Spin^c\left( n\right) \times Spin^c\left( n^{\prime
}\right) \rightarrow Spin^c\left( n+n^{\prime }\right) ,  \label{2.42a}
\end{equation}
given by formula $\left[ \left( \beta ,z\right) ,\left( \beta
^{\prime },z^{\prime }\right) \right] \rightarrow \left[ \omega
\left( \beta ,\beta ^{\prime }\right) ,zz^{\prime }\right] ,$
where $\omega $ is the homomorphism making the chain of maps into
a commutative diagram. Here, $z,z^{\prime }\in U\left( 1\right) .$
It is obvious that $L\left( \lambda ^c\oplus \lambda ^{c^{\prime
}}\right) $ is isomorphic to $L\left( \lambda ^c\right) \otimes
L\left( \lambda ^{c^{\prime }}\right).$

We conclude this subsection by formulating our main result on complex spinor
structures for $H^{2n}$-spaces:

\begin{theorem}
Let $\lambda ^c$ be a complex spinor bundle of rank $n$ and $H^{2n}$-space
considered as a real vector bundle $\lambda ^c\oplus \lambda ^{c^{\prime }}$
provided with almost complex structure $J_{\quad \beta }^\alpha ;$
multiplication on $i$ is given by $\left(
\begin{array}{cc}
0 & -\delta _j^i \\
\delta _j^i & 0%
\end{array}
\right) $. Then, the chain of maps define a commutative diagram up
to isomorphisms $\epsilon ^c $ and $\widetilde{\epsilon }^c$
defined as in (\ref{2.40a}), $\mathcal{H}$
is functor $E^c\rightarrow E^c\otimes L\left( \lambda ^c\right) $ and $%
\mathcal{E}_{\mathcal{C}}^{0,2n}\left( M^n\right) $ is defined by functor $%
\mathcal{E}_{\mathcal{C}}\left( M^n\right) \rightarrow \mathcal{E}_{\mathcal{%
C}}^{0,2n}\left( M^n\right) $ given as correspondence $E^c\rightarrow
\Lambda \left( \mathcal{C}^n\right) \otimes E^c$ (which is a categorial
equivalence), $\Lambda \left( \mathcal{C}^n\right) $ is the exterior algebra
on $\mathcal{C}^n.$ $W$ is the real bundle $\lambda ^c\oplus \lambda
^{c^{\prime }}$ provided with complex structure.
\end{theorem}

\textbf{Proof: } See Refs. \cite{vjmp,vhsp,vbook,vp}.

Now consider bundle $P\times _{Spin^c\left( n\right) }Spin^c\left( 2n\right)
$ as the principal $Spin^c\left( 2n\right) $-bundle, associated to $M\oplus
M $ being the factor space of the product $P\times Spin^c\left( 2n\right) $
on the equivalence relation $\left( p,t,h\right) \sim \left( p,\mu \left(
t\right) ^{-1}h\right) .$ In this case the categorial equivalence (\ref%
{2.40a}) can be rewritten as
\begin{equation*}
\epsilon ^c\left( E^c\right) =P\times _{Spin^c\left( n\right) }Spin^c\left(
2n\right) \Delta _{Spin^c\left( 2n\right) }E^c
\end{equation*}
and seen as factor space of $P\times Spin^c\left( 2n\right) \times _ME^c$ on
equivalence relation
\begin{equation*}
\left( pt,h,e\right) \sim \left( p,\mu \left( t\right) ^{-1}h,e\right) %
\mbox{and}\left( p,h_1,h_2,e\right) \sim \left( p,h_1,h_2^{-1}e\right)
\end{equation*}
(projections of elements $p$ and $e$ coincides on base $M).$ Every element
of $\epsilon ^c\left( E^c\right) $ can be represented as $P\Delta
_{Spin^c\left( n\right) }E^c,$ i.e., as a factor space $P\Delta E^c$ on
equivalence relation $\left( pt,e\right) \sim \left( p,\mu ^c\left( t\right)
,e\right) ,$ when $t\in Spin^c\left( n\right) .$ The complex line bundle $%
L\left( \lambda ^c\right) $ can be interpreted as the factor space of $%
P\times _{Spin^c\left( n\right) }\mathcal{C}$ on equivalence relation $%
\left( pt,\delta \right) \sim \left( p,r\left( t\right) ^{-1}\delta \right)
. $

Putting $\left( p,e\right) \otimes \left( p,\delta \right) \left( p,\delta
e\right) $ we introduce morphism
\begin{equation*}
\epsilon ^{c}\left( E\right) \times L\left( \lambda ^{c}\right) \rightarrow
\epsilon ^{c}\left( \lambda ^{c}\right)
\end{equation*}%
with properties
\begin{eqnarray*}
\left( pt,e\right) \otimes \left( pt,\delta \right) &\rightarrow &\left(
pt,\delta e\right) =\left( p,\mu ^{c}\left( t\right) ^{-1}\delta e\right) ,
\\
\left( p,\mu ^{c}\left( t\right) ^{-1}e\right) \otimes \left( p,l\left(
t\right) ^{-1}e\right) &\rightarrow &\left( p,\mu ^{c}\left( t\right)
r\left( t\right) ^{-1}\delta e\right)
\end{eqnarray*}%
pointing to the fact that we have defined the isomorphism correctly and that
it is an isomorphism on every fiber. $\Box $

\section{Spinors and N--Connection Geometry}

The purpose of this section is to show how a corresponding abstract spinor
technique entailing notational and calculations advantages can be developed
for arbitrary splits of dimensions of a d-vector space $\mathcal{F}=h%
\mathcal{F}\oplus v_{1}\mathcal{F}\oplus ...\oplus v_{z}\mathcal{F}$, where $%
\dim h\mathcal{F}=n$ and $\dim v_{p}\mathcal{F}=m_{p}.$ For convenience we
shall also present some necessary coordinate expressions.

\subsection{ D--Spinor Techniques}

The problem of a rigorous definition of spinors on locally anisotropic
spaces (d--spinors) was posed and solved \cite{vjmp,vsp1} in the framework
of the formalism of Clifford and spinor structures on v-bundles provided
with compatible nonlinear and distinguished connections and metric. We
introduced d-spinors as corresponding objects of the Clifford d-algebra $%
\mathcal{C}\left( \mathcal{F},G\right) $, defined for a d-vector space $%
\mathcal{F}$ in a standard manner (see, for instance, \cite{kar}) and proved
that operations with $\mathcal{C}\left( \mathcal{F},G\right) \ $ can be
reduced to calculations for $\mathcal{C}\left( h\mathcal{F},g\right) ,%
\mathcal{C}\left( v_{1}\mathcal{F},h_{1}\right) ,...$ and $\mathcal{C}\left(
v_{z}\mathcal{F},h_{z}\right) ,$ which are usual Clifford algebras of
respective dimensions $2^{n},2^{m_{1}},...$ and $2^{m_{z}}$ (if it is
necessary we can use quadratic forms $g$ and $h_{p}$ correspondingly induced
on $h\mathcal{F}$ and $v_{p}\mathcal{F}$ by a metric $\mathbf{G}$ (\ref%
{dmetrichcv})). Considering the orthogonal subgroup $O\mathbf{\left(
G\right) }\subset GL\mathbf{\left( G\right) }$ defined by a metric $\mathbf{G%
}$ we can define the d-spinor norm and parametrize d-spinors by ordered
pairs of elements of Clifford algebras $\mathcal{C}\left( h\mathcal{F}%
,g\right) $ and $\mathcal{C}\left( v_{p}\mathcal{F},h_{p}\right)
,p=1,2,...z. $ We emphasize that the splitting of a Clifford d-algebra
associated to a dv-bundle $\mathcal{E}^{<z>}$ is a straightforward
consequence of the global decomposition defining a N-connection structure in
$\mathcal{E}^{<z>}$.

In this subsection we shall omit detailed proofs which in most cases are
mechanical but rather tedious. We can apply the methods developed in \cite%
{pen,hladik} in a straightforward manner on h- and v-subbundles in order to
verify the correctness of affirmations.

\subsubsection{Clifford d--algebra, d--spinors and d--twistors}

In order to relate the succeeding constructions with Clifford d-algebras %
\cite{vjmp} we consider a la-frame decomposition of the metric (\ref%
{dmetrichcv}):
\begin{equation*}
G_{<\alpha ><\beta >}\left( u\right) =l_{<\alpha >}^{<\widehat{\alpha }%
>}\left( u\right) l_{<\beta >}^{<\widehat{\beta }>}\left( u\right) G_{<%
\widehat{\alpha }><\widehat{\beta }>},
\end{equation*}%
where the frame d-vectors and constant metric matrices are distinguished as

\begin{eqnarray*}
l_{<\alpha >}^{<\widehat{\alpha }>}\left( u\right) &=&\left(
\begin{array}{cccc}
l_j^{\widehat{j}}\left( u\right) & 0 & ... & 0 \\
0 & l_{a_1}^{\widehat{a}_1}\left( u\right) & ... & 0 \\
... & ... & ... & ... \\
0 & 0 & .. & l_{a_z}^{\widehat{a}_z}\left( u\right)%
\end{array}
\right) , \\
G_{<\widehat{\alpha }><\widehat{\beta }>} &=&\left(
\begin{array}{cccc}
g_{\widehat{i}\widehat{j}} & 0 & ... & 0 \\
0 & h_{\widehat{a}_1\widehat{b}_1} & ... & 0 \\
... & ... & ... & ... \\
0 & 0 & 0 & h_{\widehat{a}_z\widehat{b}_z}%
\end{array}
\right) ,
\end{eqnarray*}
$g_{\widehat{i}\widehat{j}}$ and $h_{\widehat{a}_1\widehat{b}_1},...,h_{%
\widehat{a}_z\widehat{b}_z}$ are diagonal matrices with $g_{\widehat{i}%
\widehat{i}}=$ $h_{\widehat{a}_1\widehat{a}_1}=...=h_{\widehat{a}_z\widehat{b%
}_z}=\pm 1.$

To generate Clifford d-algebras we start with matrix equations
\begin{equation}
\sigma _{<\widehat{\alpha }>}\sigma _{<\widehat{\beta }>}+\sigma _{<\widehat{%
\beta }>}\sigma _{<\widehat{\alpha }>}=-G_{<\widehat{\alpha }><\widehat{%
\beta }>}I,  \label{2.43a}
\end{equation}
where $I$ is the identity matrix, matrices $\sigma _{<\widehat{\alpha }%
>}\,(\sigma $-objects) act on a d-vector space $\mathcal{F}=h\mathcal{F}%
\oplus v_1\mathcal{F}\oplus ...\oplus v_z\mathcal{F}$ and theirs components
are distinguished as
\begin{equation}
\sigma _{<\widehat{\alpha }>}\,=\left\{ (\sigma _{<\widehat{\alpha }>})_{%
\underline{\beta }}^{\cdot \underline{\gamma }}=\left(
\begin{array}{cccc}
(\sigma _{\widehat{i}})_{\underline{j}}^{\cdot \underline{k}} & 0 & ... & 0
\\
0 & (\sigma _{\widehat{a}_1})_{\underline{b}_1}^{\cdot \underline{c}_1} & ...
& 0 \\
... & ... & ... & ... \\
0 & 0 & ... & (\sigma _{\widehat{a}_z})_{\underline{b}_z}^{\cdot \underline{c%
}_z}%
\end{array}
\right) \right\} ,  \label{2.44a}
\end{equation}
indices \underline{$\beta $},\underline{$\gamma $},... refer to spin spaces
of type $\mathcal{S}=S_{(h)}\oplus S_{(v_1)}\oplus ...\oplus S_{(v_z)}$ and
underlined Latin indices \underline{$j$},$\underline{k},...$ and $\underline{%
b}_1,\underline{c}_1,...,\underline{b}_z,\underline{c}_z...$ refer
respectively to h-spin space $\mathcal{S}_{(h)}$ and v$_p$-spin space $%
\mathcal{S}_{(v_p)},(p=1,2,...,z)\ $which are correspondingly associated to
a h- and v$_p$-decomposition of a dv-bundle $\mathcal{E}^{<z>}.$ The
irreducible algebra of matrices $\sigma _{<\widehat{\alpha }>}$ of minimal
dimension $N\times N,$ where $N=N_{(n)}+N_{(m_1)}+...+N_{(m_z)},$ $\dim
\mathcal{S}_{(h)}$=$N_{(n)}$ and $\dim \mathcal{S}_{(v_p)}$=$N_{(m_p)},$ has
these dimensions
\begin{equation*}
{N_{(n)}} ={\left\{
\begin{array}{rl}
{\ 2^{(n-1)/2},} & n=2k+1 \\
{2^{n/2},\ } & n=2k;%
\end{array}
\right. }\quad \mbox{and\ } {N}_{(m_p)}{} = {}\left|
\begin{array}{cc}
2^{(m_p-1)/2}, & m_p=2k_p+1 \\
2^{m_p}, & m_p=2k_p%
\end{array}
\right| ,
\end{equation*}
where $k=1,2,...,k_p=1,2,....$

The Clifford d-algebra is generated by sums on $n+1$ elements of form
\begin{equation}
A_1I+B^{\widehat{i}}\sigma _{\widehat{i}}+C^{\widehat{i}\widehat{j}}\sigma _{%
\widehat{i}\widehat{j}}+D^{\widehat{i}\widehat{j}\widehat{k}}\sigma _{%
\widehat{i}\widehat{j}\widehat{k}}+...  \label{2.45a}
\end{equation}
and sums of $m_p+1$ elements of form
\begin{equation*}
A_{2(p)}I+B^{\widehat{a}_p}\sigma _{\widehat{a}_p}+C^{\widehat{a}_p\widehat{b%
}_p}\sigma _{\widehat{a}_p\widehat{b}_p}+D^{\widehat{a}_p\widehat{b}_p%
\widehat{c}_p}\sigma _{\widehat{a}_p\widehat{b}_p\widehat{c}_p}+...
\end{equation*}
with antisymmetric coefficients $C^{\widehat{i}\widehat{j}}=C^{[\widehat{i}%
\widehat{j}]},C^{\widehat{a}_p\widehat{b}_p}=C^{[\widehat{a}_p\widehat{b}%
_p]},D^{\widehat{i}\widehat{j}\widehat{k}}=D^{[\widehat{i}\widehat{j}%
\widehat{k}]},D^{\widehat{a}_p\widehat{b}_p\widehat{c}_p}=D^{[\widehat{a}_p%
\widehat{b}_p\widehat{c}_p]}, $ ... and matrices $\sigma _{\widehat{i}%
\widehat{j}}=\sigma _{[\widehat{i}}\sigma _{\widehat{j}]},\sigma _{\widehat{a%
}_p\widehat{b}_p}=\sigma _{[\widehat{a}_p}\sigma _{\widehat{b}_p]},\sigma _{%
\widehat{i}\widehat{j}\widehat{k}}=\sigma _{[\widehat{i}}\sigma _{\widehat{j}%
}\sigma _{\widehat{k}]},.... $ Really, we have 2$^{n+1}$ coefficients $%
\left( A_1,C^{\widehat{i}\widehat{j}},D^{\widehat{i}\widehat{j}\widehat{k}%
},...\right) $ and 2$^{m_p+1}$ coefficients $(A_{2(p)},C^{\widehat{a}_p%
\widehat{b}_p},D^{\widehat{a}_p\widehat{b}_p\widehat{c}_p},...)$ of the
Clifford algebra on $\mathcal{F}$.

For simplicity, we shall present the necessary geometric constructions only
for h-spin spaces $\mathcal{S}_{(h)}$ of dimension $N_{(n)}.$ Considerations
for a v-spin space $\mathcal{S}_{(v)}$ are similar but with proper
characteristics for a dimension $N_{(m)}.$

In order to define the scalar (spinor) product on $\mathcal{S}_{(h)}$ we
introduce into consideration this finite sum (because of a finite number of
elements $\sigma _{\lbrack \widehat{i}\widehat{j}...\widehat{k}]}$):
\begin{equation}
^{(\pm )}E_{\underline{k}\underline{m}}^{\underline{i}\underline{j}}=\delta
_{\underline{k}}^{\underline{i}}\delta _{\underline{m}}^{\underline{j}}+%
\frac{2}{1!}(\sigma _{\widehat{i}})_{\underline{k}}^{.\underline{i}}(\sigma
^{\widehat{i}})_{\underline{m}}^{.\underline{j}}+\frac{2^{2}}{2!}(\sigma _{%
\widehat{i}\widehat{j}})_{\underline{k}}^{.\underline{i}}(\sigma ^{\widehat{i%
}\widehat{j}})_{\underline{m}}^{.\underline{j}}+\frac{2^{3}}{3!}(\sigma _{%
\widehat{i}\widehat{j}\widehat{k}})_{\underline{k}}^{.\underline{i}}(\sigma
^{\widehat{i}\widehat{j}\widehat{k}})_{\underline{m}}^{.\underline{j}}+...
\label{2.46a}
\end{equation}%
which can be factorized as
\begin{equation}
^{(\pm )}E_{\underline{k}\underline{m}}^{\underline{i}\underline{j}}=N_{(n)}{%
\ }^{(\pm )}\epsilon _{\underline{k}\underline{m}}{\ }^{(\pm )}\epsilon ^{%
\underline{i}\underline{j}}\mbox{ for }n=2k  \label{2.47a}
\end{equation}%
and
\begin{eqnarray}
^{(+)}E_{\underline{k}\underline{m}}^{\underline{i}\underline{j}}
&=&2N_{(n)}\epsilon _{\underline{k}\underline{m}}\epsilon ^{\underline{i}%
\underline{j}},{\ }^{(-)}E_{\underline{k}\underline{m}}^{\underline{i}%
\underline{j}}=0\mbox{ for }n=3(mod4),  \label{2.48a} \\
^{(+)}E_{\underline{k}\underline{m}}^{\underline{i}\underline{j}} &=&0,{\ }%
^{(-)}E_{\underline{k}\underline{m}}^{\underline{i}\underline{j}%
}=2N_{(n)}\epsilon _{\underline{k}\underline{m}}\epsilon ^{\underline{i}%
\underline{j}}\mbox{ for }n=1(mod4).  \notag
\end{eqnarray}

Antisymmetry of $\sigma _{\widehat{i}\widehat{j}\widehat{k}...}$ and the
construction of the objects (\ref{2.45a})--(\ref{2.48a}) define the
properties of $\epsilon $-objects $^{(\pm )}\epsilon _{\underline{k}%
\underline{m}}$ and $\epsilon _{\underline{k}\underline{m}}$ which have an
eight-fold periodicity on $n$ (see details in \cite{pen} and, with respect
to locally anisotropic spaces, \cite{vjmp}).

For even values of $n$ it is possible the decomposition of every h-spin
space $\mathcal{S}_{(h)}$ into irreducible h-spin spaces $\mathbf{S}_{(h)}$
and $\mathbf{S}_{(h)}^{\prime }$ (one considers splitting of h-indices, for
instance, \underline{$l$}$=L\oplus L^{\prime },\underline{m}=M\oplus
M^{\prime },...;$ for v$_p$-indices we shall write $\underline{a}%
_p=A_p\oplus A_p^{\prime },\underline{b}_p=B_p\oplus B_p^{\prime },...)$ and
defines new $\epsilon $-objects
\begin{equation}
\epsilon ^{\underline{l}\underline{m}}=\frac 12\left( ^{(+)}\epsilon ^{%
\underline{l}\underline{m}}+^{(-)}\epsilon ^{\underline{l}\underline{m}%
}\right) \mbox{ and }\widetilde{\epsilon }^{\underline{l}\underline{m}%
}=\frac 12\left( ^{(+)}\epsilon ^{\underline{l}\underline{m}}-^{(-)}\epsilon
^{\underline{l}\underline{m}}\right)  \label{2.49a}
\end{equation}
We shall omit similar formulas for $\epsilon $-objects with lower indices.

In general, the spinor $\epsilon $-objects should be defined for every shell
of an\-iso\-tro\-py where instead of dimension $n$ we shall consider the
dimensions $m_p$, $1\leq p\leq z,$ of shells.

We define a d--spinor space $\mathcal{S}_{(n,m_{1})}\ $ as a direct sum of a
horizontal and a vertical spinor spaces, for instance,
\begin{eqnarray*}
\mathcal{S}_{(8k,8k^{\prime })} &=&\mathbf{S}_{\circ }\oplus \mathbf{S}%
_{\circ }^{\prime }\oplus \mathbf{S}_{|\circ }\oplus \mathbf{S}_{|\circ
}^{\prime },\mathcal{S}_{(8k,8k^{\prime }+1)}\ =\mathbf{S}_{\circ }\oplus
\mathbf{S}_{\circ }^{\prime }\oplus \mathcal{S}_{|\circ }^{(-)},..., \\
\mathcal{S}_{(8k+4,8k^{\prime }+5)} &=&\mathbf{S}_{\triangle }\oplus \mathbf{%
S}_{\triangle }^{\prime }\oplus \mathcal{S}_{|\triangle }^{(-)},...
\end{eqnarray*}%
The scalar product on a $\mathcal{S}_{(n,m_{1})}\ $ is induced by
(corresponding to fixed values of $n$ and $m_{1}$) $\epsilon $-objects
considered for h- and v$_{1}$-components. We present also an example for $%
\mathcal{S}_{(n,m_{1}+...+m_{z})}:$%
\begin{equation}
\mathcal{S}_{(8k+4,8k_{(1)}+5,...,8k_{(p)}+4,...8k_{(z)})}=[\mathbf{S}%
_{\triangle }\oplus \mathbf{S}_{\triangle }^{\prime }\oplus \mathcal{S}%
_{|(1)\triangle }^{(-)}\oplus ...\oplus \mathbf{S}_{|(p)\triangle }\oplus
\mathbf{S}_{|(p)\triangle }^{\prime }\oplus ...\oplus \mathbf{S}_{|(z)\circ
}\oplus \mathbf{S}_{|(z)\circ }^{\prime }.  \notag
\end{equation}

Having introduced d-spinors for dimensions $\left( n,m_1+...+m_z\right) $ we
can write out the generalization for ha--spaces of twistor equations \cite%
{pen} by using the distinguished $\sigma $-objects (\ref{2.44a}):
\begin{equation}
(\sigma _{(<\widehat{\alpha }>})_{|\underline{\beta }|}^{..\underline{\gamma
}}\quad \frac{\delta \omega ^{\underline{\beta }}}{\delta u^{<\widehat{\beta
}>)}}=\frac 1{n+m_1+...+m_z}\quad G_{<\widehat{\alpha }><\widehat{\beta }%
>}(\sigma ^{\widehat{\epsilon }})_{\underline{\beta }}^{..\underline{\gamma }%
}\quad \frac{\delta \omega ^{\underline{\beta }}}{\delta u^{\widehat{%
\epsilon }}},  \label{2.56a}
\end{equation}
where $\left| \underline{\beta }\right| $ denotes that we do not consider
symmetrization on this index. The general solution of (\ref{2.56a}) on the
d-vector space $\mathcal{F}$ looks like as
\begin{equation}
\omega ^{\underline{\beta }}=\Omega ^{\underline{\beta }}+u^{<\widehat{%
\alpha }>}(\sigma _{<\widehat{\alpha }>})_{\underline{\epsilon }}^{..%
\underline{\beta }}\Pi ^{\underline{\epsilon }},  \label{2.57a}
\end{equation}
where $\Omega ^{\underline{\beta }}$ and $\Pi ^{\underline{\epsilon }}$ are
constant d-spinors. For fixed values of dimensions $n$ and $m=m_1+...m_z$ we
mast analyze the reduced and irreducible components of h- and v$_p$-parts of
equations (\ref{2.56a}) and their solutions (\ref{2.57a}) in order to find
the symmetry properties of a d-twistor $\mathbf{Z^\alpha \ }$ defined as a
pair of d-spinors $\mathbf{Z}^\alpha =(\omega ^{\underline{\alpha }},\pi _{%
\underline{\beta }}^{\prime }),$ where $\pi _{\underline{\beta }^{\prime
}}=\pi _{\underline{\beta }^{\prime }}^{(0)}\in {\widetilde{\mathcal{S}}}%
_{(n,m_1,...,m_z)}$ is a constant dual d-spinor. The problem of definition
of spinors and twistors on ha-spaces was firstly considered in \cite{vdeb}
(see also \cite{v87}) in connection with the possibility to extend the
equations (\ref{2.57a}) and theirs solutions (\ref{2.58a}), by using nearly
autoparallel maps, on curved, locally isotropic or anisotropic, spaces. We
note that the definition of twistors have been extended to higher order
anisotropic spaces with trivial N-- and d--connections.

\subsubsection{ Mutual transforms of d-tensors and d-spinors}

The spinor algebra for spaces of higher dimensions can not be considered as
a real alternative to the tensor algebra as for locally isotropic spaces of
dimensions $n=3,4$ \cite{pen}. The same holds true for ha--spaces and we
emphasize that it is not quite convenient to perform a spinor calculus for
dimensions $n,m>>4$. Nevertheless, the concept of spinors is important for
every type of spaces, we can deeply understand the fundamental properties of
geometical objects on ha--spaces, and we shall consider in this subsection
some questions concerning transforms of d-tensor objects into d-spinor ones.

\subsubsection{ Transformation of d-tensors into d-spinors}

In order to pass from d-tensors to d-spinors we must use $\sigma $-objects (%
\ref{2.44a}) written in reduced or irreduced form \quad (in dependence of
fixed values of dimensions $n$ and $m$):

\begin{eqnarray}
&&(\sigma _{<\widehat{\alpha }>})_{\underline{\beta }}^{\cdot \underline{%
\gamma }},~(\sigma ^{<\widehat{\alpha }>})^{\underline{\beta }\underline{%
\gamma }},~(\sigma ^{<\widehat{\alpha }>})_{\underline{\beta }\underline{%
\gamma }},...,(\sigma _{<\widehat{a}>})^{\underline{b}\underline{c}%
},...,(\sigma _{\widehat{i}})_{\underline{j}\underline{k}},...,(\sigma _{<%
\widehat{a}>})^{AA^{\prime }},...,(\sigma ^{\widehat{i}})_{II^{\prime }},....
\label{2.58a} \\
&&  \notag
\end{eqnarray}%
It is obvious that contracting with corresponding $\sigma $-objects (\ref%
{2.58a}) we can introduce instead of d-tensors indices the d-spinor ones,
for instance,
\begin{equation*}
\omega ^{\underline{\beta }\underline{\gamma }}=(\sigma ^{<\widehat{\alpha }%
>})^{\underline{\beta }\underline{\gamma }}\omega _{<\widehat{\alpha }%
>},\quad \omega _{AB^{\prime }}=(\sigma ^{<\widehat{a}>})_{AB^{\prime
}}\omega _{<\widehat{a}>},\quad ...,\zeta _{\cdot \underline{j}}^{\underline{%
i}}=(\sigma ^{\widehat{k}})_{\cdot \underline{j}}^{\underline{i}}\zeta _{%
\widehat{k}},....
\end{equation*}%
For d-tensors containing groups of antisymmetric indices there is a more
simple procedure of theirs transforming into d-spinors because the objects
\begin{equation}
(\sigma _{\widehat{\alpha }\widehat{\beta }...\widehat{\gamma }})^{%
\underline{\delta }\underline{\nu }},\quad (\sigma ^{\widehat{a}\widehat{b}%
...\widehat{c}})^{\underline{d}\underline{e}},\quad ...,(\sigma ^{\widehat{i}%
\widehat{j}...\widehat{k}})_{II^{\prime }},\quad ...  \label{2.59a}
\end{equation}%
can be used for sets of such indices into pairs of d-spinor indices. Let us
enumerate some properties of $\sigma $-objects of type (\ref{2.59a}) (for
simplicity we consider only h-components having q indices $\widehat{i},%
\widehat{j},\widehat{k},...$ taking values from 1 to $n;$ the properties of v%
$_{p}$-components can be written in a similar manner with respect to indices
$\widehat{a}_{p},\widehat{b}_{p},\widehat{c}_{p}...$ taking values from 1 to
$m$):
\begin{equation}
(\sigma _{\widehat{i}...\widehat{j}})^{\underline{k}\underline{l}}%
\mbox{
 is\ }\left\{ \
\begin{array}{c}
\mbox{symmetric on }\underline{k},\underline{l}\mbox{ for
}n-2q\equiv 1,7~(mod~8); \\
\mbox{antisymmetric on }\underline{k},\underline{l}\mbox{ for
}n-2q\equiv 3,5~(mod~8)%
\end{array}%
\right\}  \label{2.60a}
\end{equation}%
for odd values of $n,$ and an object
\begin{equation}
(\sigma _{\widehat{i}...\widehat{j}})^{IJ}~\left( (\sigma _{\widehat{i}...%
\widehat{j}})^{I^{\prime }J^{\prime }}\right) \mbox{ is\ }\left\{
\begin{array}{c}
\mbox{symmetric on }I,J~(I^{\prime },J^{\prime })\mbox{ for
}n-2q\equiv 0~(mod~8); \\
\mbox{antisymmetric on }I,J~(I^{\prime },J^{\prime })\mbox{ for
}n-2q\equiv 4~(mod~8)%
\end{array}%
\right\}  \label{2.61a}
\end{equation}%
or
\begin{equation}
(\sigma _{\widehat{i}...\widehat{j}})^{IJ^{\prime }}=\pm (\sigma _{\widehat{i%
}...\widehat{j}})^{J^{\prime }I}\{%
\begin{array}{c}
n+2q\equiv 6(mod8); \\
n+2q\equiv 2(mod8),%
\end{array}
\label{2.62a}
\end{equation}%
with vanishing of the rest of reduced components of the d-tensor $(\sigma _{%
\widehat{i}...\widehat{j}})^{\underline{k}\underline{l}}$ with prime/
unprime sets of indices.

\subsubsection{ Fundamental d--spinors}

We can transform every d--spinor $\xi ^{\underline{\alpha }}=\left( \xi ^{%
\underline{i}},\xi ^{\underline{a}_{1}},...,\xi ^{\underline{a}_{z}}\right) $
into a corresponding d-tensor. For simplicity, we consider this construction
only for a h-component $\xi ^{\underline{i}}$ on a h-space being of
dimension $n$. The values
\begin{equation}
\xi ^{\underline{\alpha }}\xi ^{\underline{\beta }}(\sigma ^{\widehat{i}...%
\widehat{j}})_{\underline{\alpha }\underline{\beta }}\quad \left( n%
\mbox{ is
odd}\right)  \label{2.63a}
\end{equation}%
or
\begin{equation}
\xi ^{I}\xi ^{J}(\sigma ^{\widehat{i}...\widehat{j}})_{IJ}~\left( \mbox{or }%
\xi ^{I^{\prime }}\xi ^{J^{\prime }}(\sigma ^{\widehat{i}...\widehat{j}%
})_{I^{\prime }J^{\prime }}\right) ~\left( n\mbox{ is even}\right)
\label{2.64a}
\end{equation}%
with a different number of indices $\widehat{i}...\widehat{j},$ taken
together, defines the h-spinor $\xi ^{\underline{i}}\,$ to an accuracy to
the sign. We emphasize that it is necessary to choose only those
h-components of d-tensors (\ref{2.63a}) (or (\ref{2.64a})) which are
symmetric on pairs of indices $\underline{\alpha }\underline{\beta }$ (or $%
IJ\,$ (or $I^{\prime }J^{\prime }$ )) and the number $q$ of indices $%
\widehat{i}...\widehat{j}$ satisfies the condition (as a respective
consequence of the properties (\ref{2.60a}) and/ or (\ref{2.61a}), (\ref%
{2.62a}))
\begin{equation}
n-2q\equiv 0,1,7~(mod~8).  \label{2.65a}
\end{equation}%
Of special interest is the case when
\begin{equation}
q=\frac{1}{2}\left( n\pm 1\right) ~\left( n\mbox{ is odd}\right)
\label{2.66a}
\end{equation}%
or
\begin{equation}
q=\frac{1}{2}n~\left( n\mbox{ is even}\right) .  \label{2.67a}
\end{equation}%
If all expressions (\ref{2.63a}) and/or (\ref{2.64a}) are zero for all
values of $q\,$ with the exception of one or two ones defined by the
conditions (\ref{2.65a}), (\ref{2.66a}) (or (\ref{2.67a})), the value $\xi ^{%
\widehat{i}}$ (or $\xi ^{I}$ ($\xi ^{I^{\prime }}))$ is called a fundamental
h-spinor. Defining in a similar manner the fundamental v-spinors we can
introduce fundamental d-spinors as pairs of fundamental h- and v-spinors.
Here we remark that a h(v$_{p}$)-spinor $\xi ^{\widehat{i}}~(\xi ^{\widehat{a%
}_{p}})\,$ (we can also consider reduced components) is always a fundamental
one for $n(m)<7,$ which is a consequence of (\ref{2.67a})).

\subsection{ Differential Geometry of Ha--Spinors}

This subsection is devoted to the differential geometry of d--spinors in
higher order anisotro\-pic spaces. We shall use denotations of type
\begin{equation*}
v^{<\alpha >}=(v^i,v^{<a>})\in \sigma ^{<\alpha >}=(\sigma ^i,\sigma ^{<a>})
\end{equation*}
and
\begin{equation*}
\zeta ^{\underline{\alpha }_p}=(\zeta ^{\underline{i}_p},\zeta ^{\underline{a%
}_p})\in \sigma ^{\alpha _p}=(\sigma ^{i_p},\sigma ^{a_p})\,
\end{equation*}
for, respectively, elements of modules of d-vector and irreduced d-spinor
fields (see details in \cite{vjmp}). D-tensors and d-spinor tensors
(irreduced or reduced) will be interpreted as elements of corresponding $%
\mathcal{\sigma }$--modules, for instance,
\begin{equation*}
q_{~<\beta >...}^{<\alpha >}\in \mathcal{\sigma ^{<\alpha >}\mathbf{/}}%
^{\prime };[-0\mathcal{_{<\beta >}},\psi _{~\underline{\beta }_p\quad ...}^{%
\underline{\alpha }_p\quad \underline{\gamma }_p}\in \mathcal{\sigma }_{~%
\underline{\beta _p}\quad ...}^{\underline{\alpha }_p\quad \underline{\gamma
}_p}~,\xi _{\quad J_pK_p^{\prime }N_p^{\prime }}^{I_pI_p^{\prime }}\in
\mathcal{\sigma }_{\quad J_pK_p^{\prime }N_p^{\prime }}^{I_pI_p^{\prime
}}~,...
\end{equation*}

We can establish a correspondence between the higher order anisotropic
adapted to the N--connection metric $g_{\alpha \beta }$ (\ref{dmetrichcv})
and d-spinor metric $\epsilon _{\underline{\alpha }\underline{\beta }}$ ($%
\epsilon $-objects for both h- and v$_{p}$-subspaces of $\mathcal{E}^{<z>}$)
of a ha--space $\mathcal{E}^{<z>}$ by using the relation
\begin{equation}
g_{<\alpha ><\beta >}=-\frac{1}{N(n)+N(m_{1})+...+N(m_{z})}\times ((\sigma
_{(<\alpha >}(u))^{\underline{\alpha }\underline{\beta }}(\sigma _{<\beta
>)}(u))^{\underline{\delta }\underline{\gamma }})\epsilon _{\underline{%
\alpha }\underline{\gamma }}\epsilon _{\underline{\beta }\underline{\delta }%
},  \label{2.68a}
\end{equation}%
where
\begin{equation}
(\sigma _{<\alpha >}(u))^{\underline{\nu }\underline{\gamma }}=l_{<\alpha
>}^{<\widehat{\alpha }>}(u)(\sigma _{<\widehat{\alpha }>})^{<\underline{\nu }%
><\underline{\gamma }>},  \label{2.69a}
\end{equation}%
which is a consequence of formulas (\ref{2.43a})--(\ref{2.49a}). In brief we
can write (\ref{2.68a}) as
\begin{equation}
g_{<\alpha ><\beta >}=\epsilon _{\underline{\alpha }_{1}\underline{\alpha }%
_{2}}\epsilon _{\underline{\beta }_{1}\underline{\beta }_{2}}  \label{2.70a}
\end{equation}%
if the $\sigma $-objects are considered as a fixed structure, whereas $%
\epsilon $-objects are treated as caring the metric ''dynamics '' , on
higher order anisotropic space. This variant is used, for instance, in the
so-called 2-spinor geometry \cite{pen} and should be preferred if we have to
make explicit the algebraic symmetry properties of d-spinor objects by using
metric decomposition (\ref{2.70a}). An alternative way is to consider as
fixed the algebraic structure of $\epsilon $-objects and to use variable
components of $\sigma $-objects of type (\ref{2.69a}) for developing a
variational d-spinor approach to gravitational and matter field interactions
on ha-spaces (the spinor Ashtekar variables \cite{ash} are introduced in
this manner).

We note that a d--spinor metric
\begin{equation*}
\epsilon _{\underline{\nu }\underline{\tau }}=\left(
\begin{array}{cccc}
\epsilon _{\underline{i}\underline{j}} & 0 & ... & 0 \\
0 & \epsilon _{\underline{a}_1\underline{b}_1} & ... & 0 \\
... & ... & ... & ... \\
0 & 0 & ... & \epsilon _{\underline{a}_z\underline{b}_z}%
\end{array}
\right)
\end{equation*}
on the d-spinor space $\mathcal{S}=(\mathcal{S}_{(h)},\mathcal{S}%
_{(v_1)},...,\mathcal{S}_{(v_z)})$ can have symmetric or antisymmetric h (v$%
_p$) -components $\epsilon _{\underline{i}\underline{j}}$ ($\epsilon _{%
\underline{a}_p\underline{b}_p})$ , see $\epsilon $-objects. For simplicity,
in order to avoid cumbersome calculations connected with eight-fold
periodicity on dimensions $n$ and $m_p$ of a ha-space $\mathcal{E}^{<z>},$
we shall develop a general d-spinor formalism only by using irreduced spinor
spaces $\mathcal{S}_{(h)}$ and $\mathcal{S}_{(v_p)}.$

\subsubsection{ D-covariant derivation on ha--spaces}

Let $\mathcal{E}^{<z>}$ be a ha-space. We define the action on a d-spinor of
a d-covariant operator
\begin{eqnarray*}
\bigtriangledown _{<\alpha >} &=&\left( \bigtriangledown
_{i},\bigtriangledown _{<a>}\right) =(\sigma _{<\alpha >})^{\underline{%
\alpha }_{1}\underline{\alpha }_{2}}\bigtriangledown _{^{\underline{\alpha }%
_{1}\underline{\alpha }_{2}}}=\left( (\sigma _{i})^{\underline{i}_{1}%
\underline{i}_{2}}\bigtriangledown _{^{\underline{i}_{1}\underline{i}%
_{2}}},~(\sigma _{<a>})^{\underline{a}_{1}\underline{a}_{2}}\bigtriangledown
_{^{\underline{a}_{1}\underline{a}_{2}}}\right) \\
&=&((\sigma _{i})^{\underline{i}_{1}\underline{i}_{2}}\bigtriangledown _{^{%
\underline{i}_{1}\underline{i}_{2}}},~(\sigma _{a_{1}})^{\underline{a}_{1}%
\underline{a}_{2}}\bigtriangledown _{(1)^{\underline{a}_{1}\underline{a}%
_{2}}},...,(\sigma _{a_{p}})^{\underline{a}_{1}\underline{a}%
_{2}}\bigtriangledown _{(p)^{\underline{a}_{1}\underline{a}%
_{2}}},...,(\sigma _{a_{z}})^{\underline{a}_{1}\underline{a}%
_{2}}\bigtriangledown _{(z)^{\underline{a}_{1}\underline{a}_{2}}})
\end{eqnarray*}%
(in brief, we shall write
\begin{equation*}
\bigtriangledown _{<\alpha >}=\bigtriangledown _{^{\underline{\alpha }_{1}%
\underline{\alpha }_{2}}}=\left( \bigtriangledown _{^{\underline{i}_{1}%
\underline{i}_{2}}},~\bigtriangledown _{(1)^{\underline{a}_{1}\underline{a}%
_{2}}},...,\bigtriangledown _{(p)^{\underline{a}_{1}\underline{a}%
_{2}}},...,\bigtriangledown _{(z)^{\underline{a}_{1}\underline{a}%
_{2}}}\right) )
\end{equation*}%
as maps
\begin{equation*}
\bigtriangledown _{{\underline{\alpha }}_{1}{\underline{\alpha }}_{2}}\ :\
\mathcal{\sigma }^{\underline{\beta }}\rightarrow \sigma _{<\alpha >}^{%
\underline{\beta }}=\sigma _{{\underline{\alpha }}_{1}{\underline{\alpha }}%
_{2}}^{\underline{\beta }}=
\end{equation*}%
\begin{equation*}
\left( \sigma _{i}^{\underline{\beta }}=\sigma _{{\underline{i}}_{1}{%
\underline{i}}_{2}}^{\underline{\beta }},\sigma _{(1)a_{1}}^{\underline{%
\beta }}=\sigma _{(1){\underline{\alpha }}_{1}{\underline{\alpha }}_{2}}^{%
\underline{\beta }},...,\sigma _{(p)a_{p}}^{\underline{\beta }}=\sigma _{(p){%
\underline{\alpha }}_{1}{\underline{\alpha }}_{2}}^{\underline{\beta }%
},...,\sigma _{(z)a_{z}}^{\underline{\beta }}=\sigma _{(z){\underline{\alpha
}}_{1}{\underline{\alpha }}_{2}}^{\underline{\beta }}\right)
\end{equation*}%
satisfying conditions
\begin{equation*}
\bigtriangledown _{<\alpha >}(\xi ^{\underline{\beta }}+\eta ^{\underline{%
\beta }})=\bigtriangledown _{<\alpha >}\xi ^{\underline{\beta }%
}+\bigtriangledown _{<\alpha >}\eta ^{\underline{\beta }},
\end{equation*}%
and
\begin{equation*}
\bigtriangledown _{<\alpha >}(f\xi ^{\underline{\beta }})=f\bigtriangledown
_{<\alpha >}\xi ^{\underline{\beta }}+\xi ^{\underline{\beta }%
}\bigtriangledown _{<\alpha >}f
\end{equation*}%
for every $\xi ^{\underline{\beta }},\eta ^{\underline{\beta }}\in \mathcal{%
\sigma ^{\underline{\beta }}}$ and $f$ being a scalar field on $\mathcal{E}%
^{<z>}.\mathcal{\ }$ It is also required that one holds the Leibnitz rule
\begin{equation*}
(\bigtriangledown _{<\alpha >}\zeta _{\underline{\beta }})\eta ^{\underline{%
\beta }}=\bigtriangledown _{<\alpha >}(\zeta _{\underline{\beta }}\eta ^{%
\underline{\beta }})-\zeta _{\underline{\beta }}\bigtriangledown _{<\alpha
>}\eta ^{\underline{\beta }}
\end{equation*}%
and that $\bigtriangledown _{<\alpha >}\,$ is a real operator, i.e. it
commuters with the operation of complex conjugation:
\begin{equation*}
\overline{\bigtriangledown _{<\alpha >}\psi _{\underline{\alpha }\underline{%
\beta }\underline{\gamma }...}}=\bigtriangledown _{<\alpha >}(\overline{\psi
}_{\underline{\alpha }\underline{\beta }\underline{\gamma }...}).
\end{equation*}

Let now analyze the question on uniqueness of action on d--spinors of an
operator $\bigtriangledown _{<\alpha >}$ satisfying necessary conditions .
Denoting by $\bigtriangledown _{<\alpha >}^{(1)}$ and $\bigtriangledown
_{<\alpha >}$ two such d-covariant operators we consider the map
\begin{equation}
(\bigtriangledown _{<\alpha >}^{(1)}-\bigtriangledown _{<\alpha >}):\mathcal{%
\sigma ^{\underline{\beta }}\rightarrow \sigma _{\underline{\alpha }_{1}%
\underline{\alpha }_{2}}^{\underline{\beta }}}.  \label{2.71a}
\end{equation}%
Because the action on a scalar $f$ of both operators $\bigtriangledown
_{\alpha }^{(1)}$ and $\bigtriangledown _{\alpha }$ must be identical, i.e.
\begin{equation*}
\bigtriangledown _{<\alpha >}^{(1)}f=\bigtriangledown _{<\alpha >}f,
\end{equation*}%
the action (\ref{2.71a}) on $f=\omega _{\underline{\beta }}\xi ^{\underline{%
\beta }}$ must be written as
\begin{equation*}
(\bigtriangledown _{<\alpha >}^{(1)}-\bigtriangledown _{<\alpha >})(\omega _{%
\underline{\beta }}\xi ^{\underline{\beta }})=0.
\end{equation*}%
In consequence we conclude that there is an element $\Theta _{\underline{%
\alpha }_{1}\underline{\alpha }_{2}\underline{\beta }}^{\quad \quad
\underline{\gamma }}\in \mathcal{\sigma }_{\underline{\alpha }_{1}\underline{%
\alpha }_{2}\underline{\beta }}^{\quad \quad \underline{\gamma }}$ for which
\begin{equation}
\bigtriangledown _{\underline{\alpha }_{1}\underline{\alpha }_{2}}^{(1)}\xi
^{\underline{\gamma }}=\bigtriangledown _{\underline{\alpha }_{1}\underline{%
\alpha }_{2}}\xi ^{\underline{\gamma }}+\Theta _{\underline{\alpha }_{1}%
\underline{\alpha }_{2}\underline{\beta }}^{\quad \quad \underline{\gamma }%
}\xi ^{\underline{\beta }},\bigtriangledown _{\underline{\alpha }_{1}%
\underline{\alpha }_{2}}^{(1)}\omega _{\underline{\beta }}=\bigtriangledown
_{\underline{\alpha }_{1}\underline{\alpha }_{2}}\omega _{\underline{\beta }%
}-\Theta _{\underline{\alpha }_{1}\underline{\alpha }_{2}\underline{\beta }%
}^{\quad \quad \underline{\gamma }}\omega _{\underline{\gamma }}~.
\label{2.72a}
\end{equation}%
The action of the operator (\ref{2.71a}) on a d-vector $v^{<\beta >}=v^{%
\underline{\beta }_{1}\underline{\beta }_{2}}$ can be written by using
formula (\ref{2.72a}) for both indices $\underline{\beta }_{1}$ and $%
\underline{\beta }_{2}$ :
\begin{eqnarray*}
(\bigtriangledown _{<\alpha >}^{(1)}-\bigtriangledown _{<\alpha >})v^{%
\underline{\beta }_{1}\underline{\beta }_{2}} &=&\Theta _{<\alpha >%
\underline{\gamma }}^{\quad \underline{\beta }_{1}}v^{\underline{\gamma }%
\underline{\beta }_{2}}+\Theta _{<\alpha >\underline{\gamma }}^{\quad
\underline{\beta }_{2}}v^{\underline{\beta }_{1}\underline{\gamma }} \\
&=&(\Theta _{<\alpha >\underline{\gamma }_{1}}^{\quad \underline{\beta }%
_{1}}\delta _{\underline{\gamma }_{2}}^{\quad \underline{\beta }_{2}}+\Theta
_{<\alpha >\underline{\gamma }_{1}}^{\quad \underline{\beta }_{2}}\delta _{%
\underline{\gamma }_{2}}^{\quad \underline{\beta }_{1}})v^{\underline{\gamma
}_{1}\underline{\gamma }_{2}}=Q_{\ <\alpha ><\gamma >}^{<\beta >}v^{<\gamma
>},
\end{eqnarray*}%
where
\begin{equation}
Q_{\ <\alpha ><\gamma >}^{<\beta >}=Q_{\qquad \underline{\alpha }_{1}%
\underline{\alpha }_{2}~\underline{\gamma }_{1}\underline{\gamma }_{2}}^{%
\underline{\beta }_{1}\underline{\beta }_{2}}=\Theta _{<\alpha >\underline{%
\gamma }_{1}}^{\quad \underline{\beta }_{1}}\delta _{\underline{\gamma }%
_{2}}^{\quad \underline{\beta }_{2}}+\Theta _{<\alpha >\underline{\gamma }%
_{1}}^{\quad \underline{\beta }_{2}}\delta _{\underline{\gamma }_{2}}^{\quad
\underline{\beta }_{1}}.  \label{2.73a}
\end{equation}%
The d-commutator $\bigtriangledown _{\lbrack <\alpha >}\bigtriangledown
_{<\beta >]}$ defines the d-torsion. So, applying operators\newline
$\bigtriangledown _{\lbrack <\alpha >}^{(1)}\bigtriangledown _{<\beta
>]}^{(1)}$ and $\bigtriangledown _{\lbrack <\alpha >}\bigtriangledown
_{<\beta >]}$ on $f=\omega _{\underline{\beta }}\xi ^{\underline{\beta }}$
we can write
\begin{equation*}
T_{\quad <\alpha ><\beta >}^{(1)<\gamma >}-T_{~<\alpha ><\beta >}^{<\gamma
>}=Q_{~<\beta ><\alpha >}^{<\gamma >}-Q_{~<\alpha ><\beta >}^{<\gamma >}
\end{equation*}%
with $Q_{~<\alpha ><\beta >}^{<\gamma >}$ from (\ref{2.73a}).

The action of operator $\bigtriangledown _{<\alpha >}^{(1)}$ on d-spinor
tensors of type $\chi _{\underline{\alpha }_1\underline{\alpha }_2\underline{%
\alpha }_3...}^{\qquad \quad \underline{\beta }_1\underline{\beta }_2...}$
must be constructed by using formula (\ref{2.72a}) for every upper index $%
\underline{\beta }_1\underline{\beta }_2...$ and formula (\ref{2.73a}) for
every lower index $\underline{\alpha }_1\underline{\alpha }_2\underline{%
\alpha }_3...$ .

\subsubsection{Infeld--van der Waerden co\-ef\-fi\-ci\-ents}

Let
\begin{equation*}
\delta _{\underline{\mathbf{\alpha }}}^{\quad \underline{\alpha }}=\left(
\delta _{\underline{\mathbf{1}}}^{\quad \underline{i}},\delta _{\underline{%
\mathbf{2}}}^{\quad \underline{i}},...,\delta _{\underline{\mathbf{N(n)}}%
}^{\quad \underline{i}},\delta _{\underline{\mathbf{1}}}^{\quad \underline{a}%
},\delta _{\underline{\mathbf{2}}}^{\quad \underline{a}},...,\delta _{%
\underline{\mathbf{N(m)}}}^{\quad \underline{i}}\right)
\end{equation*}%
be a d--spinor basis. The dual to it basis is denoted as
\begin{equation*}
\delta _{\underline{\alpha }}^{\quad \underline{\mathbf{\alpha }}}=\left(
\delta _{\underline{i}}^{\quad \underline{\mathbf{1}}},\delta _{\underline{i}%
}^{\quad \underline{\mathbf{2}}},...,\delta _{\underline{i}}^{\quad
\underline{\mathbf{N(n)}}},\delta _{\underline{i}}^{\quad \underline{\mathbf{%
1}}},\delta _{\underline{i}}^{\quad \underline{\mathbf{2}}},...,\delta _{%
\underline{i}}^{\quad \underline{\mathbf{N(m)}}}\right) .
\end{equation*}%
A d-spinor $\kappa ^{\underline{\alpha }}\in \mathcal{\sigma }$ $^{%
\underline{\alpha }}$ has components $\kappa ^{\underline{\mathbf{\alpha }}%
}=\kappa ^{\underline{\alpha }}\delta _{\underline{\alpha }}^{\quad
\underline{\mathbf{\alpha }}}.$ Taking into account that $\delta _{%
\underline{\mathbf{\alpha }}}^{\quad \underline{\alpha }}\delta _{\underline{%
\mathbf{\beta }}}^{\quad \underline{\beta }}\bigtriangledown _{\underline{%
\alpha }\underline{\beta }}=\bigtriangledown _{\underline{\mathbf{\alpha }}%
\underline{\mathbf{\beta }}},$ we write out the components $\bigtriangledown
_{\underline{\alpha }\underline{\beta }}$ $\kappa ^{\underline{\gamma }}$ in
the form
\begin{eqnarray}
\delta _{\underline{\mathbf{\alpha }}}^{\quad \underline{\alpha }}~\delta _{%
\underline{\mathbf{\beta }}}^{\quad \underline{\beta }}~\delta _{\underline{%
\gamma }}^{\quad \underline{\mathbf{\gamma }}}~\bigtriangledown _{\underline{%
\alpha }\underline{\beta }}\kappa ^{\underline{\gamma }} &=&  \notag \\
\delta _{\underline{\mathbf{\epsilon }}}^{\quad \underline{\tau }}~\delta _{%
\underline{\tau }}^{\quad \underline{\mathbf{\gamma }}}~\bigtriangledown _{%
\underline{\mathbf{\alpha }}\underline{\mathbf{\beta }}}\kappa ^{\underline{%
\mathbf{\epsilon }}}+\kappa ^{\underline{\mathbf{\epsilon }}}~\delta _{%
\underline{\epsilon }}^{\quad \underline{\mathbf{\gamma }}}~\bigtriangledown
_{\underline{\mathbf{\alpha }}\underline{\mathbf{\beta }}}\delta _{%
\underline{\mathbf{\epsilon }}}^{\quad \underline{\epsilon }} &=&
\bigtriangledown _{\underline{\mathbf{\alpha }}\underline{\mathbf{\beta }}%
}\kappa ^{\underline{\mathbf{\gamma }}}+\kappa ^{\underline{\mathbf{\epsilon
}}}\gamma _{~\underline{\mathbf{\alpha }}\underline{\mathbf{\beta }}%
\underline{\mathbf{\epsilon }}}^{\underline{\mathbf{\gamma }}},
\label{2.74a}
\end{eqnarray}%
where the coordinate components of the d--spinor connection $\gamma _{~%
\underline{\mathbf{\alpha }}\underline{\mathbf{\beta }}\underline{\mathbf{%
\epsilon }}}^{\underline{\mathbf{\gamma }}}$ are defined as
\begin{equation}
\gamma _{~\underline{\mathbf{\alpha }}\underline{\mathbf{\beta }}\underline{%
\mathbf{\epsilon }}}^{\underline{\mathbf{\gamma }}}\doteq \delta _{%
\underline{\tau }}^{\quad \underline{\mathbf{\gamma }}}~\bigtriangledown _{%
\underline{\mathbf{\alpha }}\underline{\mathbf{\beta }}}\delta _{\underline{%
\mathbf{\epsilon }}}^{\quad \underline{\tau }}.  \label{2.75a}
\end{equation}%
We call the Infeld - van der Waerden d-symbols a set of $\sigma $-objects ($%
\sigma _{\mathbf{\alpha }})^{\underline{\mathbf{\alpha }}\underline{\mathbf{%
\beta }}}$ parametri\-zed with respect to a coordinate d-spinor basis.
Defining
\begin{equation*}
\bigtriangledown _{<\mathbf{\alpha >}}=(\sigma _{<\mathbf{\alpha >}})^{%
\underline{\mathbf{\alpha }}\underline{\mathbf{\beta }}}~\bigtriangledown _{%
\underline{\mathbf{\alpha }}\underline{\mathbf{\beta }}},
\end{equation*}%
introducing denotations
\begin{equation*}
\gamma ^{\underline{\mathbf{\gamma }}}{}_{<\mathbf{\alpha >\underline{\tau }}%
}\doteq \gamma ^{\underline{\mathbf{\gamma }}}{}_{\mathbf{\underline{\alpha }%
\underline{\beta }\underline{\tau }}}~(\sigma _{<\mathbf{\alpha >}})^{%
\underline{\mathbf{\alpha }}\underline{\mathbf{\beta }}}
\end{equation*}%
and using properties (\ref{2.74a}), we can write the relations
\begin{eqnarray}
l_{<\mathbf{\alpha >}}^{<\alpha >}~\delta _{\underline{\beta }}^{\quad
\underline{\mathbf{\beta }}}~\bigtriangledown _{<\alpha >}\kappa ^{%
\underline{\beta }} &=&\bigtriangledown _{<\mathbf{\alpha >}}\kappa ^{%
\underline{\mathbf{\beta }}}+\kappa ^{\underline{\mathbf{\delta }}}~\gamma
_{~<\mathbf{\alpha >}\underline{\mathbf{\delta }}}^{\underline{\mathbf{\beta
}}},  \label{2.76a} \\
l_{<\mathbf{\alpha >}}^{<\alpha >}~\delta _{\underline{\mathbf{\beta }}%
}^{\quad \underline{\beta }}~\bigtriangledown _{<\alpha >}~\mu _{\underline{%
\beta }} &=&\bigtriangledown _{<\mathbf{\alpha >}}~\mu _{\underline{\mathbf{%
\beta }}}-\mu _{\underline{\mathbf{\delta }}}\gamma _{~<\mathbf{\alpha >}%
\underline{\mathbf{\beta }}}^{\underline{\mathbf{\delta }}}.  \notag
\end{eqnarray}%
for d-covariant derivations $~\bigtriangledown _{\underline{\alpha }}\kappa
^{\underline{\beta }}$ and $\bigtriangledown _{\underline{\alpha }}~\mu _{%
\underline{\beta }}.$

We can consider expressions similar to (\ref{2.76a}) for values having both
types of d-spinor and d-tensor indices, for instance,
\begin{equation}
l_{<\mathbf{\alpha >}}^{<\alpha >}~l_{<\gamma >}^{<\mathbf{\gamma >}}~\delta
_{\underline{\mathbf{\delta }}}^{\quad \underline{\delta }}~\bigtriangledown
_{<\alpha >}\theta _{\underline{\delta }}^{~<\gamma >}=\bigtriangledown _{<%
\mathbf{\alpha >}}\theta _{\underline{\mathbf{\delta }}}^{~<\mathbf{\gamma >}%
}-\theta _{\underline{\mathbf{\epsilon }}}^{~<\mathbf{\gamma >}}\gamma _{~<%
\mathbf{\alpha >}\underline{\mathbf{\delta }}}^{\underline{\mathbf{\epsilon }%
}}+\theta _{\underline{\mathbf{\delta }}}^{~<\mathbf{\tau >}}~\Gamma _{\quad
<\mathbf{\alpha ><\tau >}}^{~<\mathbf{\gamma >}}  \notag
\end{equation}%
(we can prove this by a straightforward calculation).

Now we shall consider some possible relations between components of
d-connec\-ti\-ons $\gamma _{~<\mathbf{\alpha >}\underline{\mathbf{\delta }}%
}^{\underline{\mathbf{\epsilon }}}$ and $\Gamma _{\quad <\mathbf{\alpha
><\tau >}}^{~<\mathbf{\gamma >}}$ and derivations of $(\sigma _{<\mathbf{%
\alpha >}})^{\underline{\mathbf{\alpha }}\underline{\mathbf{\beta }}}$ . We
can write
\begin{eqnarray*}
\Gamma _{~<\mathbf{\beta ><\gamma >}}^{<\mathbf{\alpha >}} &=&l_{<\alpha
>}^{<\mathbf{\alpha >}}\bigtriangledown _{<\mathbf{\gamma >}}l_{<\mathbf{%
\beta >}}^{<\alpha >}=l_{<\alpha >}^{<\mathbf{\alpha >}}\bigtriangledown _{<%
\mathbf{\gamma >}}(\sigma _{<\mathbf{\beta >}})^{\underline{\epsilon }%
\underline{\tau }}l_{<\alpha >}^{<\mathbf{\alpha >}}\bigtriangledown _{<%
\mathbf{\gamma >}}((\sigma _{<\mathbf{\beta >}})^{\underline{\mathbf{%
\epsilon }}\underline{\mathbf{\tau }}}\delta _{\underline{\mathbf{\epsilon }}%
}^{~\underline{\epsilon }}\delta _{\underline{\mathbf{\tau }}}^{~\underline{%
\tau }}) \\
&=&l_{<\alpha >}^{<\mathbf{\alpha >}}\delta _{\underline{\mathbf{\alpha }}%
}^{~\underline{\alpha }}\delta _{\underline{\mathbf{\epsilon }}}^{~%
\underline{\epsilon }}\bigtriangledown _{<\mathbf{\gamma >}}(\sigma _{<%
\mathbf{\beta >}})^{\underline{\mathbf{\alpha }}\underline{\mathbf{\epsilon }%
}}+l_{<\alpha >}^{<\mathbf{\alpha >}}(\sigma _{<\mathbf{\beta >}})^{%
\underline{\mathbf{\epsilon }}\underline{\mathbf{\tau }}}(\delta _{%
\underline{\mathbf{\tau }}}^{~\underline{\tau }}\bigtriangledown _{<\mathbf{%
\gamma >}}\delta _{\underline{\mathbf{\epsilon }}}^{~\underline{\epsilon }%
}+\delta _{\underline{\mathbf{\epsilon }}}^{~\underline{\epsilon }%
}\bigtriangledown _{<\mathbf{\gamma >}}\delta _{\underline{\mathbf{\tau }}%
}^{~\underline{\tau }}) \\
&=&l_{\underline{\mathbf{\epsilon }}\underline{\mathbf{\tau }}}^{<\mathbf{%
\alpha >}}~\bigtriangledown _{<\mathbf{\gamma >}}(\sigma _{<\mathbf{\beta >}%
})^{\underline{\mathbf{\epsilon }}\underline{\mathbf{\tau }}}+l_{\underline{%
\mathbf{\mu }}\underline{\mathbf{\nu }}}^{<\mathbf{\alpha >}}\delta _{%
\underline{\epsilon }}^{~\underline{\mathbf{\mu }}}\delta _{\underline{\tau }%
}^{~\underline{\mathbf{\nu }}}(\sigma _{<\mathbf{\beta >}})^{\underline{%
\epsilon }\underline{\tau }}(\delta _{\underline{\mathbf{\tau }}}^{~%
\underline{\tau }}\bigtriangledown _{<\mathbf{\gamma >}}\delta _{\underline{%
\mathbf{\epsilon }}}^{~\underline{\epsilon }}+\delta _{\underline{\mathbf{%
\epsilon }}}^{~\underline{\epsilon }}\bigtriangledown _{<\mathbf{\gamma >}%
}\delta _{\underline{\mathbf{\tau }}}^{~\underline{\tau }})
\end{eqnarray*}%
where $l_{<\alpha >}^{<\mathbf{\alpha >}}=(\sigma _{\underline{\mathbf{%
\epsilon }}\underline{\mathbf{\tau }}})^{<\mathbf{\alpha >}}$ , from which
one follows
\begin{equation}
(\sigma _{<\mathbf{\alpha >}})^{\underline{\mathbf{\mu }}\underline{\mathbf{%
\nu }}}(\sigma _{\underline{\mathbf{\alpha }}\underline{\mathbf{\beta }}})^{<%
\mathbf{\beta >}}\Gamma _{~<\mathbf{\gamma ><\beta >}}^{<\mathbf{\alpha >}%
}=(\sigma _{\underline{\mathbf{\alpha }}\underline{\mathbf{\beta }}})^{<%
\mathbf{\beta >}}\bigtriangledown _{<\mathbf{\gamma >}}(\sigma _{<\mathbf{%
\alpha >}})^{\underline{\mathbf{\mu }}\underline{\mathbf{\nu }}}+\delta _{%
\underline{\mathbf{\beta }}}^{~\underline{\mathbf{\nu }}}\gamma _{~<\mathbf{%
\gamma >\underline{\alpha }}}^{\underline{\mathbf{\mu }}}+\delta _{%
\underline{\mathbf{\alpha }}}^{~\underline{\mathbf{\mu }}}\gamma _{~<\mathbf{%
\gamma >\underline{\beta }}}^{\underline{\mathbf{\nu }}}.  \notag
\end{equation}%
Connecting the last expression on \underline{$\mathbf{\beta }$} and
\underline{$\mathbf{\nu }$} and using an orthonormalized d-spinor basis when
$\gamma _{~<\mathbf{\gamma >\underline{\beta }}}^{\underline{\mathbf{\beta }}%
}=0$ (a consequence from (\ref{2.75a})) we have
\begin{equation}
\gamma _{~<\mathbf{\gamma >\underline{\alpha }}}^{\underline{\mathbf{\mu }}}=%
\frac{1}{N(n)+N(m_{1})+...+N(m_{z})}(\Gamma _{\quad <\mathbf{\gamma >~%
\underline{\alpha }\underline{\beta }}}^{\underline{\mathbf{\mu }}\underline{%
\mathbf{\beta }}}-(\sigma _{\underline{\mathbf{\alpha }}\underline{\mathbf{%
\beta }}})^{<\mathbf{\beta >}}\bigtriangledown _{<\mathbf{\gamma >}}(\sigma
_{<\mathbf{\beta >}})^{\underline{\mathbf{\mu }}\underline{\mathbf{\beta }}%
}),  \label{2.78a}
\end{equation}%
where
\begin{equation}
\Gamma _{\quad <\mathbf{\gamma >~\underline{\alpha }\underline{\beta }}}^{%
\underline{\mathbf{\mu }}\underline{\mathbf{\beta }}}=(\sigma _{<\mathbf{%
\alpha >}})^{\underline{\mathbf{\mu }}\underline{\mathbf{\beta }}}(\sigma _{%
\underline{\mathbf{\alpha }}\underline{\mathbf{\beta }}})^{\mathbf{\beta }%
}\Gamma _{~<\mathbf{\gamma ><\beta >}}^{<\mathbf{\alpha >}}.  \label{2.79a}
\end{equation}%
We also note here that, for instance, for the canonical and Berwald
connections and Christoffel d-symbols we can express d-spinor connection (%
\ref{2.79a}) through corresponding locally adapted derivations of components
of metric and N-connection by introducing corresponding coefficients instead
of $\Gamma _{~<\mathbf{\gamma ><\beta >}}^{<\mathbf{\alpha >}}$ in (\ref%
{2.79a}) and than in (\ref{2.78a}).

\subsubsection{ D--spinors of ha--space curvature and torsion}

The d-tensor indices of the commutator $\Delta _{<\alpha ><\beta >}$ can be
transformed into d-spinor ones:
\begin{eqnarray}
\Box _{\underline{\alpha }\underline{\beta }} &=&(\sigma ^{<\alpha ><\beta
>})_{\underline{\alpha }\underline{\beta }}\Delta _{\alpha \beta }=(\Box _{%
\underline{i}\underline{j}},\Box _{\underline{a}\underline{b}})=(\Box _{%
\underline{i}\underline{j}},\Box _{\underline{a}_{1}\underline{b}%
_{1}},...,\Box _{\underline{a}_{p}\underline{b}_{p}},...,\Box _{\underline{a}%
_{z}\underline{b}_{z}}),  \label{2.80a} \\
&&  \notag
\end{eqnarray}%
with h- and v$_{p}$-components,
\begin{equation*}
\Box _{\underline{i}\underline{j}}=(\sigma ^{<\alpha ><\beta >})_{\underline{%
i}\underline{j}}\Delta _{<\alpha ><\beta >}\mbox{ and }\Box _{\underline{a}%
\underline{b}}=(\sigma ^{<\alpha ><\beta >})_{\underline{a}\underline{b}%
}\Delta _{<\alpha ><\beta >},
\end{equation*}%
being symmetric or antisymmetric in dependence of corresponding values of
dimensions $n\,$ and $m_{p}$ (see eight-fold parametizations. Considering
the actions of operator (\ref{2.80a}) on d-spinors $\pi ^{\underline{\gamma }%
}$ and $\mu _{\underline{\gamma }}$ we introduce the d-spinor curvature $X_{%
\underline{\delta }\quad \underline{\alpha }\underline{\beta }}^{\quad
\underline{\gamma }}\,$ as to satisfy equations
\begin{equation}
\Box _{\underline{\alpha }\underline{\beta }}\ \pi ^{\underline{\gamma }}=X_{%
\underline{\delta }\quad \underline{\alpha }\underline{\beta }}^{\quad
\underline{\gamma }}\pi ^{\underline{\delta }}\mbox{ and }\Box _{\underline{%
\alpha }\underline{\beta }}\ \mu _{\underline{\gamma }}=X_{\underline{\gamma
}\quad \underline{\alpha }\underline{\beta }}^{\quad \underline{\delta }}\mu
_{\underline{\delta }}.  \label{2.81a}
\end{equation}%
The gravitational d-spinor $\Psi _{\underline{\alpha }\underline{\beta }%
\underline{\gamma }\underline{\delta }}$ is defined by a corresponding
symmetrization of d-spinor indices:
\begin{equation}
\Psi _{\underline{\alpha }\underline{\beta }\underline{\gamma }\underline{%
\delta }}=X_{(\underline{\alpha }|\underline{\beta }|\underline{\gamma }%
\underline{\delta })}.  \label{2.82a}
\end{equation}%
We note that d-spinor tensors $X_{\underline{\delta }\quad \underline{\alpha
}\underline{\beta }}^{\quad \underline{\gamma }}$ and $\Psi _{\underline{%
\alpha }\underline{\beta }\underline{\gamma }\underline{\delta }}\,$ are
transformed into similar 2-spinor objects on locally isotropic spaces \cite%
{pen} if we consider vanishing of the N-connection structure and a limit to
a locally isotropic space.

Putting $\delta _{\underline{\gamma }}^{\quad \mathbf{\underline{\gamma }}}$
instead of $\mu _{\underline{\gamma }}$ in (\ref{2.81a}) and using (\ref%
{2.82a}) we can express respectively the curvature and gravitational
d-spinors as
\begin{equation*}
X_{\underline{\gamma }\underline{\delta }\underline{\alpha }\underline{\beta
}}=\delta _{\underline{\delta }\underline{\mathbf{\tau }}}\Box _{\underline{%
\alpha }\underline{\beta }}\delta _{\underline{\gamma }}^{\quad \mathbf{%
\underline{\tau }}}\mbox{ and }\Psi _{\underline{\gamma }\underline{\delta }%
\underline{\alpha }\underline{\beta }}=\delta _{\underline{\delta }%
\underline{\mathbf{\tau }}}\Box _{(\underline{\alpha }\underline{\beta }%
}\delta _{\underline{\gamma })}^{\quad \mathbf{\underline{\tau }}}.
\end{equation*}

The d-spinor torsion $T_{\qquad \underline{\alpha }\underline{\beta }}^{%
\underline{\gamma }_1\underline{\gamma }_2}$ is defined similarly as for
d-tensors by using the d-spinor commutator (\ref{2.80a}) and equations
\begin{equation*}
\Box _{\underline{\alpha }\underline{\beta }}f=T_{\qquad \underline{\alpha }%
\underline{\beta }}^{\underline{\gamma }_1\underline{\gamma }%
_2}\bigtriangledown _{\underline{\gamma }_1\underline{\gamma }_2}f.
\end{equation*}

The d-spinor components $R_{\underline{\gamma }_1\underline{\gamma }_2\qquad
\underline{\alpha }\underline{\beta }}^{\qquad \underline{\delta }_1%
\underline{\delta }_2}$ of the curvature d-tensor $R_{\gamma \quad \alpha
\beta }^{\quad \delta }$ can be computed by using relations (\ref{2.79a}),
and (\ref{2.80a}) and (\ref{2.82a}) as to satisfy the equations
\begin{equation*}
(\Box _{\underline{\alpha }\underline{\beta }}-T_{\qquad \underline{\alpha }%
\underline{\beta }}^{\underline{\gamma }_1\underline{\gamma }%
_2}\bigtriangledown _{\underline{\gamma }_1\underline{\gamma }_2})V^{%
\underline{\delta }_1\underline{\delta }_2}=R_{\underline{\gamma }_1%
\underline{\gamma }_2\qquad \underline{\alpha }\underline{\beta }}^{\qquad
\underline{\delta }_1\underline{\delta }_2}V^{\underline{\gamma }_1%
\underline{\gamma }_2},
\end{equation*}
here d-vector $V^{\underline{\gamma }_1\underline{\gamma }_2}$ is considered
as a product of d-spinors, i.e. $V^{\underline{\gamma }_1\underline{\gamma }%
_2}=\nu ^{\underline{\gamma }_1}\mu ^{\underline{\gamma }_2}$. We find
\begin{eqnarray}
R_{\underline{\gamma }_{1}\underline{\gamma }_{2}\qquad \underline{\alpha }%
\underline{\beta }}^{\qquad \underline{\delta }_{1}\underline{\delta }_{2}}
&=&\left( X_{\underline{\gamma }_{1}~\underline{\alpha }\underline{\beta }%
}^{\quad \underline{\delta }_{1}}+T_{\qquad \underline{\alpha }\underline{%
\beta }}^{\underline{\tau }_{1}\underline{\tau }_{2}}\quad \gamma _{\quad
\underline{\tau }_{1}\underline{\tau }_{2}\underline{\gamma }_{1}}^{%
\underline{\delta }_{1}}\right) \delta _{\underline{\gamma }_{2}}^{\quad
\underline{\delta }_{2}}  \notag \\
&& +\left( X_{\underline{\gamma }_{2}~\underline{\alpha }\underline{\beta }%
}^{\quad \underline{\delta }_{2}}+T_{\qquad \underline{\alpha }\underline{%
\beta }}^{\underline{\tau }_{1}\underline{\tau }_{2}}\quad \gamma _{\quad
\underline{\tau }_{1}\underline{\tau }_{2}\underline{\gamma }_{2}}^{%
\underline{\delta }_{2}}\right) \delta _{\underline{\gamma }_{1}}^{\quad
\underline{\delta }_{1}}.  \label{2.83a}
\end{eqnarray}

It is convenient to use this d--spinor expression for the curvature d-tensor
\begin{eqnarray*}
R_{\underline{\gamma }_1\underline{\gamma }_2\qquad \underline{\alpha }_1%
\underline{\alpha }_2\underline{\beta }_1\underline{\beta }_2}^{\qquad
\underline{\delta }_1\underline{\delta }_2} &=&\left( X_{\underline{\gamma }%
_1~\underline{\alpha }_1\underline{\alpha }_2\underline{\beta }_1\underline{%
\beta }_2}^{\quad \underline{\delta }_1}+T_{\qquad \underline{\alpha }_1%
\underline{\alpha }_2\underline{\beta }_1\underline{\beta }_2}^{\underline{%
\tau }_1\underline{\tau }_2}~\gamma _{\quad \underline{\tau }_1\underline{%
\tau }_2\underline{\gamma }_1}^{\underline{\delta }_1}\right) \delta _{%
\underline{\gamma }_2}^{\quad \underline{\delta }_2} \\
&&+\left( X_{\underline{\gamma }_2~\underline{\alpha }_1\underline{\alpha }_2%
\underline{\beta }_1\underline{\beta }_2}^{\quad \underline{\delta }%
_2}+T_{\qquad \underline{\alpha }_1\underline{\alpha }_2\underline{\beta }_1%
\underline{\beta }_2~}^{\underline{\tau }_1\underline{\tau }_2}\gamma
_{\quad \underline{\tau }_1\underline{\tau }_2\underline{\gamma }_2}^{%
\underline{\delta }_2}\right) \delta _{\underline{\gamma }_1}^{\quad
\underline{\delta }_1}
\end{eqnarray*}
in order to get the d--spinor components of the Ricci d-tensor
\begin{eqnarray}
& &R_{\underline{\gamma }_1\underline{\gamma }_2\underline{\alpha }_1%
\underline{\alpha }_2} = R_{\underline{\gamma }_1\underline{\gamma }_2\qquad
\underline{\alpha }_1\underline{\alpha }_2\underline{\delta }_1\underline{%
\delta }_2}^{\qquad \underline{\delta }_1\underline{\delta }_2}=X_{%
\underline{\gamma }_1~\underline{\alpha }_1\underline{\alpha }_2\underline{%
\delta }_1\underline{\gamma }_2}^{\quad \underline{\delta }_1}+
\label{2.84a} \\
&&T_{\qquad \underline{\alpha }_1\underline{\alpha }_2\underline{\delta }_1%
\underline{\gamma }_2}^{\underline{\tau }_1\underline{\tau }_2}~\gamma
_{\quad \underline{\tau }_1\underline{\tau }_2\underline{\gamma }_1}^{%
\underline{\delta }_1}+X_{\underline{\gamma }_2~\underline{\alpha }_1%
\underline{\alpha }_2\underline{\delta }_1\underline{\gamma }_2}^{\quad
\underline{\delta }_2}+T_{\qquad \underline{\alpha }_1\underline{\alpha }_2%
\underline{\gamma }_1\underline{\delta }_2~}^{\underline{\tau }_1\underline{%
\tau }_2}\gamma _{\quad \underline{\tau }_1\underline{\tau }_2\underline{%
\gamma }_2}^{\underline{\delta }_2}  \notag
\end{eqnarray}
and this d-spinor decomposition of the scalar curvature:
\begin{eqnarray}
q\overleftarrow{R} &=&R_{\qquad \underline{\alpha }_1\underline{\alpha }_2}^{%
\underline{\alpha }_1\underline{\alpha }_2}=X_{\quad ~\underline{~\alpha }%
_1\quad \underline{\delta }_1\underline{\alpha }_2}^{\underline{\alpha }_1%
\underline{\delta }_1~~\underline{\alpha }_2}+T_{\qquad ~~\underline{\alpha }%
_2\underline{\delta }_1}^{\underline{\tau }_1\underline{\tau }_2\underline{%
\alpha }_1\quad ~\underline{\alpha }_2}~\gamma _{\quad \underline{\tau }_1%
\underline{\tau }_2\underline{\alpha }_1}^{\underline{\delta }_1}
\label{2.85a} \\
&&+X_{\qquad \quad \underline{\alpha }_2\underline{\delta }_2\underline{%
\alpha }_1}^{\underline{\alpha }_2\underline{\delta }_2\underline{\alpha }%
_1}+T_{\qquad \underline{\alpha }_1\quad ~\underline{\delta }_2~}^{%
\underline{\tau }_1\underline{\tau }_2~~\underline{\alpha }_2\underline{%
\alpha }_1}\gamma _{\quad \underline{\tau }_1\underline{\tau }_2\underline{%
\alpha }_2}^{\underline{\delta }_2}.  \notag
\end{eqnarray}

Using (\ref{2.84a}) and (\ref{2.85a}), see details in Refs. \cite{pen}, we
define the d--spinor components of the Einstein and $\Phi _{<\alpha><\beta>}$
d--tensors:
\begin{eqnarray}
\overleftarrow{G}_{<\gamma ><\alpha >} &=&\overleftarrow{G}_{\underline{%
\gamma }_{1}\underline{\gamma }_{2}\underline{\alpha }_{1}\underline{\alpha }%
_{2}}=X_{\underline{\gamma }_{1}~\underline{\alpha }_{1}\underline{\alpha }%
_{2}\underline{\delta }_{1}\underline{\gamma }_{2}}^{\quad \underline{\delta
}_{1}}+T_{\qquad \underline{\alpha }_{1}\underline{\alpha }_{2}\underline{%
\delta }_{1}\underline{\gamma }_{2}}^{\underline{\tau }_{1}\underline{\tau }%
_{2}}~\gamma _{\quad \underline{\tau }_{1}\underline{\tau }_{2}\underline{%
\gamma }_{1}}^{\underline{\delta }_{1}}  \notag \\
&&+X_{\underline{\gamma }_{2}~\underline{\alpha }_{1}\underline{\alpha }_{2}%
\underline{\delta }_{1}\underline{\gamma }_{2}}^{\quad \underline{\delta }%
_{2}}+T_{\qquad \underline{\alpha }_{1}\underline{\alpha }_{2}\underline{%
\gamma }_{1}\underline{\delta }_{2}~}^{\underline{\tau }_{1}\underline{\tau }%
_{2}}\gamma _{\quad \underline{\tau }_{1}\underline{\tau }_{2}\underline{%
\gamma }_{2}}^{\underline{\delta }_{2}}-\frac{1}{2}\varepsilon _{\underline{%
\gamma }_{1}\underline{\alpha }_{1}}\varepsilon _{\underline{\gamma }_{2}%
\underline{\alpha }_{2}}[X_{\quad ~\underline{~\beta }_{1}\quad \underline{%
\mu }_{1}\underline{\beta }_{2}}^{\underline{\beta }_{1}\underline{\mu }%
_{1}~~\underline{\beta }_{2}}  \notag \\
&&+T_{\qquad ~~\underline{\beta }_{2}\underline{\mu }_{1}}^{\underline{\tau }%
_{1}\underline{\tau }_{2}\underline{\beta }_{1}\quad ~\underline{\beta }%
_{2}}~\gamma _{\quad \underline{\tau }_{1}\underline{\tau }_{2}\underline{%
\beta }_{1}}^{\underline{\mu }_{1}}+X_{\qquad \quad \underline{\beta }_{2}%
\underline{\mu }_{2}\underline{\nu }_{1}}^{\underline{\beta }_{2}\underline{%
\mu }_{2}\underline{\nu }_{1}}+T_{\qquad \underline{\beta }_{1}\quad ~%
\underline{\delta }_{2}~}^{\underline{\tau }_{1}\underline{\tau }_{2}~~%
\underline{\beta }_{2}\underline{\beta }_{1}}\gamma _{\quad \underline{\tau }%
_{1}\underline{\tau }_{2}\underline{\beta }_{2}}^{\underline{\delta }_{2}}]
\label{2.86a}
\end{eqnarray}%
and%
\begin{eqnarray}
\Phi _{<\gamma ><\alpha >} &=&\Phi _{\underline{\gamma }_{1}\underline{%
\gamma }_{2}\underline{\alpha }_{1}\underline{\alpha }_{2}}=\frac{1}{%
2(n+m_{1}+...+m_{z})}\varepsilon _{\underline{\gamma }_{1}\underline{\alpha }%
_{1}}\varepsilon _{\underline{\gamma }_{2}\underline{\alpha }_{2}}[X_{\quad ~%
\underline{~\beta }_{1}\quad \underline{\mu }_{1}\underline{\beta }_{2}}^{%
\underline{\beta }_{1}\underline{\mu }_{1}~~\underline{\beta }_{2}}+  \notag
\\
&&T_{\qquad ~~\underline{\beta }_{2}\underline{\mu }_{1}}^{\underline{\tau }%
_{1}\underline{\tau }_{2}\underline{\beta }_{1}\quad ~\underline{\beta }%
_{2}}~\gamma _{\quad \underline{\tau }_{1}\underline{\tau }_{2}\underline{%
\beta }_{1}}^{\underline{\mu }_{1}}+X_{\qquad \quad \underline{\beta }_{2}%
\underline{\mu }_{2}\underline{\nu }_{1}}^{\underline{\beta }_{2}\underline{%
\mu }_{2}\underline{\nu }_{1}}+T_{\qquad \underline{\beta }_{1}\quad ~%
\underline{\delta }_{2}~}^{\underline{\tau }_{1}\underline{\tau }_{2}~~%
\underline{\beta }_{2}\underline{\beta }_{1}}\gamma _{\quad \underline{\tau }%
_{1}\underline{\tau }_{2}\underline{\beta }_{2}}^{\underline{\delta }_{2}}]
\notag \\
&&-\frac{1}{2}[X_{\underline{\gamma }_{1}~\underline{\alpha }_{1}\underline{%
\alpha }_{2}\underline{\delta }_{1}\underline{\gamma }_{2}}^{\quad
\underline{\delta }_{1}}+T_{\qquad \underline{\alpha }_{1}\underline{\alpha }%
_{2}\underline{\delta }_{1}\underline{\gamma }_{2}}^{\underline{\tau }_{1}%
\underline{\tau }_{2}}~\gamma _{\quad \underline{\tau }_{1}\underline{\tau }%
_{2}\underline{\gamma }_{1}}^{\underline{\delta }_{1}}  \notag \\
&& + X_{\underline{\gamma }_{2}~\underline{\alpha }_{1}\underline{\alpha }%
_{2}\underline{\delta }_{1}\underline{\gamma }_{2}}^{\quad \underline{\delta
}_{2}}+T_{\qquad \underline{\alpha }_{1}\underline{\alpha }_{2}\underline{%
\gamma }_{1}\underline{\delta }_{2}~}^{\underline{\tau }_{1}\underline{\tau }%
_{2}}\gamma _{\quad \underline{\tau }_{1}\underline{\tau }_{2}\underline{%
\gamma }_{2}}^{\underline{\delta }_{2}}].  \label{2.87a}
\end{eqnarray}

We omit this calculus in this work.

\subsection*{Acknowledgements}

~~ The S. V. work is supported by a NATO/Portugal fellowship at CENTRA,
Instituto Superior Tecnico, Lisbon. The author is grateful to R. Ablamowicz,
John Ryan and B. Fauser for collaboration and support of his participation
at ''The 6th International Conference on Clifford Algebras'', Cookeville,
Tennessee, USA (May, 20-25, 2002). He would like to thank J. P. S. Lemos, R.
Miron, M. Anastasiei and P. Stavrinos for hospitality and support.


\end{document}